
\documentstyle{amsppt}
\TagsOnRight
\magnification=1100

\newif\ifFIRST\newdimen\MAXright\MAXright0pt 
\def\sdynkin{\bgroup\eightpoint\dynkin} 
\def\endsdynkin{\enddynkin\egroup} 
\def\dynkin{\bgroup\FIRSTtrue\hskip.5em\setbox1\hbox{$\diagup$}%
\setbox2\hbox{$\diagdown$}%
\setbox0\hbox to2\wd1{\hrulefill}%
\setbox3\hbox{$\ \bullet\ $}%
\setbox4\hbox{$\ \circ\ $}%
\setbox7\hbox{$\smallmatrix\otimes\endsmallmatrix$}%
\def\root##1{\ifFIRST\setbox5\hbox{$##1$}\ifdim\wd5>1.3em%
\hskip-.5em\hskip.5\wd5\fi\fi\FIRSTfalse%
\hskip-.25em\raise1.0\wd3\hbox to0pt{\hss\hskip.45em$%
\ssize##1$\hss}\copy3\hskip-.25em\setbox6\hbox{$##1$}%
\MAXright\wd6}%
\def\xroot##1{\ifFIRST\setbox5\hbox{$##1$}\ifdim\wd5>1.3em%
\hskip-.5em\hskip.5\wd5\fi\fi\FIRSTfalse%
\hskip-.25em\raise1.0\wd3\hbox to0pt{\hss\hskip.45em$%
\ssize##1$\hss}\copy7\hskip-.25em\setbox6\hbox{$##1$}%
\MAXright\wd6}%
\def\droot##1{\ifFIRST\setbox5\hbox{$##1$}\ifdim\wd5>1.3em%
\hskip-.5em\hskip.5\wd5\fi\fi\FIRSTfalse%
\hskip-.25em\lower1.8\wd3\hbox to0pt{\hss\hskip.45em$%
\ssize##1$\hss}\copy3\hskip-.25em\setbox6\hbox{$##1$}%
\MAXright\wd6}%
\def\rroot##1{\hskip-.25em\copy3\hbox to0pt{\hskip.3em$\ssize##1$\hss}%
\hskip-.25em\setbox6\hbox{\hskip.6em$##1##1$}%
\MAXright\wd6}%
\def\wrroot##1{\hskip-.36em\copy4\hbox to0pt{\hskip.3em$\ssize##1$\hss}%
\hskip-.48em\setbox6\hbox{\hskip.6em$##1##1$}%
\MAXright\wd6}%
\def\wroot##1{\ifFIRST\setbox5\hbox{$##1$}\ifdim\wd5>1.3em%
\hskip-.5em\hskip.5\wd5\fi\fi\FIRSTfalse%
\hskip-.36em\raise1.0\wd3\hbox to0pt{\hss\hskip.6em$%
\ssize##1$\hss}\copy4\hskip-.38em\setbox6\hbox{$##1$}%
\MAXright\wd6}%
\def\wdroot##1{\ifFIRST\setbox5\hbox{$##1$}\ifdim\wd5>1.3em%
\hskip-.5em\hskip.5\wd5\fi\fi\FIRSTfalse%
\hskip-.36em\lower1.8\wd3\hbox to0pt{\hss\hskip.6em$%
\ssize##1$\hss}\copy4\hskip-.38em\setbox6\hbox{$##1$}%
\MAXright\wd6}%
\def\link{\raise.22em\copy0}%
\def\llink##1{\raise.32em\copy0\hskip-\wd0%
\raise.12em\copy0\hskip-.5\wd0\hbox to0pt{\hss$##1$\hss}\hskip.5\wd0}%
\def\lllink##1{\raise.22pt\copy0\hskip-\wd0\raise.32em\copy0\hskip-\wd0%
\raise.12em\copy0\hskip-.5\wd0\hbox to0pt{\hss$##1$\hss}\hskip.5\wd0}%
\def\rootupright##1{\hbox to0pt{\raise.45em\copy1\hskip-.25em\raise1.3\ht1%
\hbox{\copy3\hskip.3em$\ssize##1$}\hss}%
\setbox6\hbox{\hskip.6em\copy1\copy1$##1##1$}%
\ifdim\MAXright<\wd6\MAXright\wd6\fi}%
\def\wrootupright##1{\hbox to0pt{\raise.45em\copy1\hskip-.36em\raise1.3\ht1%
\hbox{\copy4\hskip.1em\raise.20em\hbox{$\ssize##1$}}\hss}}%
\def\xrootupright##1{\hbox to0pt{\raise.45em\copy1\hskip-.25em\raise1.3\ht1%
\hbox{\copy7\hskip.2em\raise.20em\hbox{$\ssize##1$}}\hss}}%
\def\rootdownright##1{\hbox to0pt{\raise-.5em\copy2\hskip-.25em\raise-1.35\ht1 
\hbox{\copy3\hskip.3em$\ssize##1$}\hss}%
\setbox6\hbox{\hskip.6em\copy2\copy2$##1##1$}%
\ifdim\MAXright<\wd6\MAXright\wd6\fi}%
\def\wrootdownright##1{\hbox to0pt{\raise-.5em\copy2\hskip-.36em\raise-1.35\ht1 
\hbox{\copy4\hskip.1em\raise.1em\hbox{$\ssize##1$}}\hss}%
\setbox6\hbox{\hskip.6em\copy1\copy1$##1##1$}%
\ifdim\MAXright<\wd6\MAXright\wd6\fi}%
\def\xrootdownright##1{\hbox to0pt{\raise-.5em\copy2\hskip-.36em\raise-1.35\ht1 
\hbox{\copy7\hskip.1em\raise.1em\hbox{$\ssize##1$}}\hss}%
\setbox6\hbox{\hskip.6em\copy1\copy1$##1##1$}%
\ifdim\MAXright<\wd6\MAXright\wd6\fi}%
\def\wrootdown##1{\hbox to0pt{\hskip-.24em\vrule height.06em depth.67em%
\hskip-.59em\raise-1.15em\hbox{\copy4\hskip.3em$\ssize##1$}\hss}%
\setbox6\hbox{$##1$}%
\ifdim\MAXright<\wd6\MAXright\wd6\fi}%
\def\dots{\hskip.5em\cdots\hskip.5em}}%
\def\enddynkin{\ifdim\MAXright>1em\hskip.5\MAXright\else\hskip.5em\fi\egroup}%

\define\g{{\goth g}}
\define\h{{\goth h}}
\define\k{{\goth k}}
\redefine\l{{\goth l}}
\define\p{{\goth p}}
\define\z{{\goth z}}
\define\n{{\goth n}}
\define\m{{\goth m}}
\define\r{{\goth r}}
\define\q{{\goth q}}
\define\a{{\goth a}}
\define\={\overset\text{def}\to=}
\redefine\B{{\Cal B}}
\define\C{{\Bbb C}}
\define\R{{\Bbb R}}
\redefine\D{{\Cal D}}
\define\e{{\goth e}}

\documentstyle{amsppt}
\TagsOnRight

\topmatter

\title
Invariant CR structures on compact homogeneous manifolds\\
\endtitle

\author
Dmitry V. Alekseevsky and Andrea F. Spiro\\
\phantom{aa}
\endauthor

\address
\phantom{ }\newline
D. V. Alekseevsky\newline
Center Sophus Lie\newline
117279 Moscow\newline
gen. Antonova 2-99\newline
RUSSIA\newline 
\phantom{ }\newline
\phantom{ }
\endaddress

\email
fort\@slip.rsuh.ru
\endemail

\address
\phantom{ }\newline
A. F. Spiro\newline
Dipartimento di Matematica e Fisica\newline
Universit\`a di Camerino\newline
Via Madonna delle Carceri\newline
62032 Camerino (Macerata)\newline
ITALY\newline
\phantom{ }\newline
\phantom{ }
\endaddress

\email
spiro\@campus.unicam.it
\endemail


\leftheadtext{D. V. Alekseevsky and A. F. Spiro}
\rightheadtext{Invariant CR structures on Compact homogeneous
manifolds}

\abstract An explicit classification of  simply connected
compact homogeneous CR manifolds  $G/L$
of codimension one, with non-degenerate Levi form, is given.
There are three classes of such manifolds:\par
 a) the standard CR homogeneous 
manifolds which are  homogeneous $S^1$-bundles over a flag manifold $F$,
with CR structure induced by an invariant complex structure on $F$; \par
b) the Morimoto-Nagano spaces, i.e.  
sphere bundles $S(N)\subset TN$ of a compact
rank one symmetric space $N = G/H$, with the CR structure induced by the 
natural complex structure of $TN = G^\C/H^\C$; \par
c) the following manifolds: $SU_n/T^1\cdot SU_{n-2}$, 
$SU_p\times SU_q/T^1 \cdot U_{p-2}\cdot U_{q-2}$, 
$SU_n/T^1\cdot SU_2\cdot SU_2\cdot SU_{n-4}$, 
$SO_{10}/T^1\cdot SO_6$, $E_6/T^1\cdot SO_8$;
these manifolds admit canonical holomorphic fibrations over 
a flag manifold $(F,J_F)$ 
with typical fiber $S(S^k)$, where $k = 2, 3, 5, 7$ or $9$, respectively;
the CR structure  is determined by the invariant complex
structure  $J_F$ on $F$ and by 
an invariant CR structure on the typical fiber,
depending on one complex parameter.
\endabstract


\endtopmatter
\document
\subhead 1. Introduction
\endsubhead
\bigskip
An  {\it almost CR structure\/} on a manifold $M$ is 
a pair $(\D, J)$, where $\D\subset TM$ is a distribution and $J$ is
a complex
structure on $\D$. The complexification $\D^\C$ can be decomposed as
$\D^\C = \D^{10} + \D^{01}$ into sum of complex eigendistributions of
$J$, with eigenvalues $i$ and $-i$.\par
 An almost  CR structure is called {\it
integrable\/} or, shortly, {\it CR structure\/}
if the distribution $\D^{01}$ (and hence also
 $\D^{10}$) is involutive, i.e. with space of sections  closed
under Lie bracket.
This is equivalent to the following conditions:
$$J([JX, Y] + [X, JY])\in \D\ ,$$
$$[JX, JY] - [X, Y] - J([JX, Y] + [X, JY]) = 0\ ,$$
for any two fields $X, Y$ in $\D$.\par
A map $\varphi\: (M, \D, J) \to (M', \D', J')$ between
two CR
manifolds
is  called {\it holomorphic map\/} if
$\varphi_*(\D) \subset \D'$ and  $\varphi_*(J X) = J'\varphi_*(X)$.
\par
Two CR structures $(\D,J)$ and
$(\D',J')$  are called {\it
equivalent\/} if there exists a diffeomorphism such that $\phi_*(\D)
= \D'$ and $\phi_*J  = J'$.\par
The
codimension of
$\D$ is called {\it 
codimension of the CR structure\/}. Note that a 
CR structure of codimension zero 
is the same as  a complex structure.\par
A codimension one CR structure
$(\D,J)$ on a $2n+1$-dimensional 
manifold $M$ is called {\it Levi non-degenerate\/}
if $\D$ is a contact
distribution. This means that any  local (contact) 1-form
$\theta$, which defines the distribution (i.e. such that $ker \theta =
\D$)
is maximally non-degenerate, that is $(d\theta)^n\wedge\theta \neq 0$.
\par 
Note that any real hypersurface $M$ of a complex manifold $N$
has a natural codimension one CR structure $(\D,J_\D)$ induced by
the complex structure $J$ of $N$, where  
$$\D = \{\ X \in TM\ , \ JX \in TM\ \}\ ,\qquad
J_\D = J|_\D\ .$$
\par
In the
following,
if the opposite is not
stated,  by {\it CR structure\/}
we will mean {\it integrable codimension one Levi non-degenerate
CR structure\/}. Sometimes, if the contact distribution $\D$ is given,
we will identify  a CR structure with the associated  complex structure
$J$.\par
\medskip
A CR manifold, that is a manifold $M$ with a CR structure $(\D,J)$,
is called {\it homogeneous\/} if it admits a transitive Lie group of
holomorphic  transformations.\par
If the opposite is not stated, we will always assume that
the homogeneous CR manifold $(M,\D,J)$ is simply connected.\par 
 The aim of this
paper is to give a complete
classification of  simply connected homogeneous CR manifolds
$M = G/L$ of a compact Lie group $G$. This gives a classification of all 
simply connected homogeneous CR manifolds, since any compact homogeneous
CR manifold admits a compact transitive Lie group of 
holomorphic transformations
(see [12]).\par
\medskip 
The simplest example of compact homogeneous CR manifold is   
 the standard sphere $S^{2n-1}\subset \C^n$ 
with the induced CR structure.\par
\medskip
More elaborated examples are provided by the following construction
of A. Morimoto and T. Nagano ([9]). 
Let $N = G/H$ be a compact
rank one symmetric space
(shortly 'CROSS'). The tangent space $TN$ can be identified
with the homogeneous space $G^\C/H^\C$.
Hence,  it  admits
a natural   
$G^\C$-invariant complex structure $J$.  
Any regular orbit $G\cdot p = S(N) \simeq G/L$ in $TN = G^\C/H^\C$ is
a sphere bundle; in particular it is 
a real hypersurface with (Levi non-degenerate) $G$-invariant
CR structure. \par
Moreover, these examples together with the 
standard sphere $S^{2n-1}\subset \C^n$ exhaust the class of   
CR structures induced on a codimension one orbit 
$M = G\cdot x \subset C$
of a compact Lie group $G$ of holomorphic transformations
of a  Stein  manifold $C$. We call the homogeneous CR manifolds
which are equivalent to such  orbits $G\cdot p = S(N)$ 
in the tangent space of a CROSS
{\it Morimoto-Nagano spaces\/}.\par
\medskip
In the fundamental paper [1],
 H. Azad, A. Huckleberry and W. Richthofer showed that these manifolds
play a basic role in the description of  compact homogeneous CR 
manifolds (see also [8] and [11]). \par
More precisely, for any compact homogeneous
CR manifold $M = G/L$ they define a  
holomorphic map (called
{\it anticanonical map\/}) $\phi\: M = G/L \to \C P^N$. This map is 
$G$-equivariant with respect to some explicitly defined projective
action of $G$ on $\C P^N$. For any compact homogeneous
CR manifold $M$ only two possibilities may occur: 
the orbit  $\phi(M) = G\cdot p$, $p\in \phi(M)$, is either a flag manifold
with the  complex structure  induced by the complex structure
$J_P$ of $\C P^N$ and in this case
$\phi\: M\to \phi(M)$ is an $S^1$-fibering, or it is 
a CR manifold with CR structure induced by $J_P$ and 
in this case $\phi\: M\to \phi(M)$ is a finite covering. \par
 This reduces the description 
of  CR homogeneous manifolds of the second type to the description
of compact 
orbits $G\cdot p \subset \C P^N$ of a real subgroup 
$G\subset Aut(\C P^N)$ of projective
transformations, on which $J_P$ induces a CR structure. \par
A simple argument shows that an orbit $G\cdot p\subset \C P^N$ 
of a connected Lie subgroup $G\subset Aut(\C P^N)$
carries a  (possibly Levi degenerate)  CR structure induced by $\C P^N$
if and only if $G\cdot p $ is a real hypersurface of $G^\C\cdot p$.
Moreover, if the orbit is compact,
one may assume that $G$ is a compact semisimple Lie group.\par
The following important result in [1] describes the structure of such 
orbits. \par
\medskip
\proclaim{Theorem} Let $G^\C \subset Aut(\C P^N)$ be a connected
complex semisimple
group of projective transformations and  $G$  its compact form.
Assume that the orbit $M = G\cdot p = G/L$ carries a Levi non-degenerate
CR structure induced by $J_P$ and hence it is a real 
hypersurface in $B = G^\C\cdot p = G^\C/H$. Denote by $P$ a minimal 
parabolic subgroup of $G^\C$ which properly contains $H$. Then 
the fiber $C = P/H$ of the 
$G^\C$-equivariant fibration 
$$\pi\: B = G^\C/H \to F = G^\C/P$$
over the flag manifold $F = G^\C/P$
is a homogeneous Stein manifold biholomorphic to $\C^*$, $\C^n$
or to the tangent space of a CROSS.
\endproclaim
\medskip
This fibration is called {\it Stein-rational fibration\/}. 
Note that   $P$ not necessarily acts  effectively on $C$.\par
The Stein-rational fibration induces a 
$G$-equivariant holomorphic fibration of the homogeneous
CR manifold $M = G/L$ over the flag manifold $F$
$$\pi'\: M = G/L \to F = G^\C/P$$
(it is a {\it CRF fibration\/} according to our definitions, see
below).
Moreover, in correspondence to a fiber of $\pi$, a fiber of $\pi'$
is either $S^1$, $S^{2n-1}$ or a Morimoto-Nagano space.\par
This Theorem gives necessary conditions
for the induced CR structure
on  $M = G\cdot p\subset \C P^N$ being Levi non-degenerate. 
Our classification gives
necessary and sufficient conditions. In particular, we show
that only the sphere
bundles $S(S^k)$ with $k = 2,3,5,7, 9$ and $11$
occur as fibers of the 
fibration $\pi'$.\par
\bigskip
Now we  describe the main results of this paper.
Section \S 2 collects the basics facts on homogeneous CR manifolds.\par
Section \S 3 is devoted to the infinitesimal description of
homogeneous contact manifolds $M = G/L$ of a compact Lie group.\par
We prove that the center of $G$ is at most one
dimensional and
we establish a natural one to one correspondence between 
simply connected homogeneous manifolds $M = G/L$ with an invariant
contact distribution $\D$ and an element $Z\in \g = Lie(G)$
(defined up to scaling)
such that:\par
\roster
\item"a)"  the centralizer of $Z$ has the following
orthogonal decomposition 
$$C_\g(Z) = \l \oplus \R Z\ ,\qquad \l = Lie(L)$$
 w.r.t. the Cartan-Killing form
$\B$;
\item"b)" the 1-parametric subgroup generated by $Z$ is closed.
\endroster
\par
This element $Z$ (called {\it contact element\/}) defines an
orthogonal
decomposition 
$$\g = \l + \R Z + \m\ .$$
The subspace $\m$ is 
$\operatorname{Ad}_L$-invariant and defines the contact distribution
$\D$ on $M=G/L$, while the $\operatorname{Ad}_L$-invariant
1-form $\theta = \B\circ Z \in \g^*$ is extended to a 
$G$-invariant contact form $\theta$ on $G/L$.\par
We associate to $Z$
a flag manifold $F_Z$ , which  is the
adjoint orbit 
$$F_Z = \operatorname{Ad}_G(Z) = G/K\ ,$$
where $K = C_G(Z)$ is the centralizer of $Z$.
There is a  natural principal $S^1$-fibration 
$$\pi\: M = G/L \longrightarrow F_Z = G/K\ .$$
In general, a homogeneous manifold $G/L$ admits no more
then one invariant contact structure. If it admits more then one then 
it is called {\it special contact manifold\/}. \par
The main examples of 
such manifolds can be described as follows. \par
Let $G$ be a simple compact Lie group without center and let 
$Q = G/Sp_1\cdot H'$ be the associated Wolf space, that is the
homogeneous quaternionic K\"ahler manifold, where 
$Sp_1 \cdot H'$ is the normalizer in $G$ of the 3-dimensional
subalgebra $\goth{sp}_1(\mu)$ of $\g$ associated with the maximal 
root $\mu$. Then the associated 3-Sasakian homogeneous
manifold
$M = G/H'$ is a special contact manifold.\par
Any $0 \neq Z \in
\goth{sp}_1(\mu)$ is a contact element. 
Furthermore, any two invariant contact 
structures on $M$ are equivalent under a 
transformation, which commutes with $G$, defined by the right action  
of an element from $Sp_1$.\par
We prove the following theorem.
\medskip 
\proclaim{Theorem 1.1} Any special contact manifold
$M = G/L$ is either the 3-Sasakian homogeneous manifold
$G/H'$ of a simple Lie group $G$,  described above, or
$M = G_2/Sp_1$, where $Sp_1$ is the 3-dimensional subgroup
of the exceptional Lie group $G_2$, with Lie algebra 
$\goth{sp}_1(\mu)$, where $\mu$ is the maximal root of $G_2$.
\endproclaim
\medskip
In section \S 4 we establish some general properties of compact
homogeneous CR manifolds. Let 
$(M = G/L, \D)$ be a homogeneous contact manifold and  
$$\g = \l + \R Z + \m$$
 the associated decomposition of 
 $\g$. Then any invariant (integrable) CR structure $J$ is
defined by the $\operatorname{Ad}_L$-invariant decomposition 
$$\m^\C = \m^{10} + \m^{01} \tag1.1$$
of the complexified tangent space $\m^\C = T^\C_{eL} M$
into holomorphic and antiholomorphic subspaces; this decomposition is such that
$$\l^\C +\m^{10}\quad \text{is a subalgebra of}
\quad \g^\C\ .\tag1.2$$
 The subspace
$\m$ is naturally identified with the tangent space of the
associated flag manifold $F_Z = G/K$, $\k = \l + \R Z = Lie(K)$.
It is known that any invariant complex structure on $F_Z$ 
is defined by an  $\operatorname{Ad}_K$-invariant decomposition
(1.1), where $\m^{10}$ is a subalgebra (in fact it is the 
nilradical of a parabolic subalgebra $\k^\C + \m^{10}$).
Hence any invariant complex structure $J_F$ on $F_Z$ defines an 
invariant CR structure $J_M$ on $M = G/L$. It is called
{\it standard CR structure induced by $J_F$\/}.\par
The natural
$S^1$-fibration $\pi\: M = G/L \to F_Z = G/K$ is holomorphic 
with respect to the CR structure $J_M$ and the complex structure
$J_F$. \par
Since the description of all invariant complex structures
on a flag manifold is known (see e.g. [10],  [4], [5], [3]), it is
sufficient to classify the non-standard CR structures.\par
\medskip
The following notion is important for such classification.\par
A compact homogeneous CR manifold $(M = G/L, \D, J)$ is called
{\it non-primitive\/} if it admits a holomorphic
$G$-equivariant fibration $\pi$ (called {\it CRF-fibration\/})
$$\pi\: M = G/L \longrightarrow F = G/Q\ ,$$
where   $F = G/Q$ is a flag manifold of positive dimension,
equipped with an invariant complex
structure $J_F$. \par
Note that a fiber of $\pi$ is a homogeneous compact CR
manifold $Q/L$  and that any  standard CR manifold is
non-primitive.\par
The classification of primitive CR structures given in \S 5 and \S 
6 is an important
step for classification  of  non-standard CR structures.\par
\medskip
A basic tool for studying the homogeneous  CR manifolds is 
the anticanonical map $\phi$ defined in [1]. \par
Let $(M= G/L, \D_Z, J)$ be a homogeneous CR manifold of 
a compact Lie group $G$ and 
$$\g^\C = \l^\C + \C Z + \m^{10} + \m^{01}$$
the corresponding decomposition of $\g^\C$.  Then the
anticanonical map $\phi$ is the holomorphic map 
of $M$ into the Grassmanian of $k$-planes, 
$k = \dim_\C(\l^\C + \m^{01})$,  given by
$$\phi\: M = G/L \longrightarrow Gr_k(\g^\C) \subset \C P^N$$
$$\phi\: gL \mapsto  \operatorname{Ad}_g([\l^\C + \m^{01}])\ .$$
Note that $\phi$ is a $G$-equivariant map onto the orbit $G\cdot p$ of 
$p = [\l^\C + \m^{01}]\in Gr_k(\g^\C)$ under the natural adjoint action of 
$G$ on  $Gr_k(\g^\C)$.\par
We obtain the following characterization of standard 
CR structures (see Theorems 4.9 and 4.11):
\medskip
\proclaim{Theorem 1.2} Let $(M=G/L, \D_Z, J)$ be a homogeneous CR
manifold. 
\roster
\item If it is standard, then  the image $\phi(M) = G\cdot
p$ of the anticanonical map is the flag manifold $F_Z = G/K$, 
associated with the contact structure $\D_Z$. Hence $\phi\: M \to
\phi(M) = F_Z$ is the natural $S^1$-fibration.
\item If it is non-standard, then $\phi\: M \to \phi(M) = G\cdot p$
is a finite holomorphic covering, with respect to the CR structure
of $G\cdot p \subset Gr_k(\g^\C)$ induced by the complex
structure of $Gr_k(\g^\C)$.
\endroster
\endproclaim
 \medskip
In section \S 5, we classify all invariant CR structures on special
contact manifolds $G/L$. The result is the following:
\medskip
\proclaim{Theorem 1.3} Let $M = G/L$ be a special contact manifold
with an invariant contact structure $\D_Z$. 
\roster
\item if $G \neq SU_n$, then there exists (up to
sign of $J$) only one invariant
CR structure $(\D_Z, J)$, which is standard;
\item if $G = SU_2$ and hence $M = SU_2$, then 
any CR structure is non-primitive
and it admits a CRF fibration with typical fiber $S^1$;   
\item if $G = SU_n$, $n>2$, and hence $M = SU_n/U_{n-2}$, then there 
exist (up to sign of $J$) three standard CR structures and three 
1-parameter families of non-standard CR structures; 
one of such families consists of 
primitive   CR structures; 
any other 
non-standard CR structure
is non-primitive and it admits a 
CRF fibration 
$$\pi_1: M = SU_n/U_{n-2} \longrightarrow Q_1=SU_n/T^1\cdot U_{n-2}$$ 
with fiber $S^1$ over the flag manifold $Q_1 = SU_n/T^1\cdot U_{n-2}$; 
moreover, any non-standard, non-primitive CR structure
admits also a CRF fibration
$$\pi\: M = SU_n/U_{n-2} \longrightarrow Q_2=SU_n/S(U_2\cdot U_{n-2})$$
with  fiber $SO_3$ and base given by the Wolf space 
$Q_2 = SU_n/S(U_2\cdot U_{n-2})$, equipped with its (unique up to sign)
complex structure. 
\endroster
\endproclaim
The explicit description of all non-standard CR structures 
on $SU_2$ and $SU_n/U_{n-2}$ is given in Proposition 5.1.
\medskip
In section \S 6, we obtain  the classification of non-standard
invariant CR structures on non-special homogeneous
contact manifolds.\par
 Together with above results it leads to the following classification
of  primitive CR structures.\par
\medskip
\proclaim{Theorem 1.4} Let $(M= G/L, \D_Z, J)$ be a simply
connected, primitive,
homogeneous CR manifold and let
$\theta = \B\circ Z|_{\goth t}$ be the dual form of the contact element
$Z$ restricted
to a Cartan subalgebra $\goth t$ of $\k = C_\g(Z) = \l + \R Z$.    
Then $G/L$ is  isomorphic to the universal covering space of
a sphere bundle $S(N)\subset T(N)$ of 
a CROSS $N$.  The groups
$G$, $K = C_G(Z)$, the form $\vartheta = -i \theta$
and the CROSS $N$  are listed
in
the
following table.  
\bigskip
\moveright 0.2 cm
\vbox{\offinterlineskip
\halign {\strut\vrule\hfil\ $#$\ \hfil
 &\vrule\hfil\ $#$\
\hfil&\vrule\hfil\ $#$\
\hfil&\vrule\hfil\  $#$\
\hfil & \vrule\hfil\ $#$
\hfil
\vrule\cr
\noalign{\hrule} n^o &
\phantom{\frac{\frac{1}{1}}{\frac{1}{1}}}G
\ \ &
K = C_G(Z)
&
\vartheta
&  
N = G/H
\cr \noalign{\hrule}
1 &
SU_2\times SU'_2
&
T^1\times T^1{}'
&
(\varepsilon_1  - \varepsilon_2) + (\varepsilon_1'
-\varepsilon_2')
&   
S^3 = \frac{SO_4}{SO_3}^{\phantom{\sum^A_B}}_{\phantom{\sum^A_B}}
\cr \noalign{\hrule}
2 &
Spin_7
&
T^1\cdot SU_3
&
\varepsilon_1 +\varepsilon_2 + \varepsilon_3
&
S^{7} = \frac{Spin_7}{G_2}^{\phantom
{\sum^A_B}}_{\phantom{\sum^A_B}}
\cr \noalign{\hrule}
3 &
F_4
&
T^1\cdot SO_7
&
\varepsilon_1
&
\Bbb O P^2 = 
\frac{F_4}{Spin_9}^{\phantom{\sum^A_B}}_{\phantom{\sum^A_B}}
\cr \noalign{\hrule}
4 &
SO_{2n+1}\ 
\smallmatrix{n > 1}
\endsmallmatrix
&
T^1\cdot SO_{2n-1}
&
\varepsilon_1 
&
S^{2n} = \frac{SO_{2n+1}}{SO_{2n}}^{\phantom
{\sum^A_B}}_{\phantom{\sum^A_B}}
\cr \noalign{\hrule}
5 &
SO_{2n}\ 
\smallmatrix{n > 2}
\endsmallmatrix
&
T^1\cdot SO_{2n-2}
&
\varepsilon_1
&
S^{2n-1} = \frac{SO_{2n}}{SO_{2n-1}}^{\phantom{\sum^A_B}}
_{\phantom{\sum^A_B}}
\cr \noalign{\hrule}
6 &
SU_{n+1}
\ 
\smallmatrix{n > 1}
\endsmallmatrix
&
T^1 \cdot U_{n-1}
&
\varepsilon_1- \varepsilon_2 
&
\C P^n = \frac{SU_{n+1}}{U_{n}}^{\phantom
{\sum^A_B}}_{\phantom{\sum^A_B}}
\cr \noalign{\hrule}
7 &
Sp_n
&
T^1\cdot Sp_1\cdot Sp_{n-2}
&
\varepsilon_1 + \varepsilon_2{}
&
\Bbb H P^{n-1} = \underset{\phantom{B}}\to
{\overset{\phantom{A}}\to{\frac{Sp_n}{
Sp_1\cdot Sp_{n-1}}}}
\cr \noalign{\hrule}
}}
\noindent 
In each of these cases, the set of all invariant
 CR structures (considered up to 
sign) on $M = G/L$ is
parameterized by the points of the unit disc $D$ in $\R^2$. The center of $D$
corresponds to the (unique) standard CR structure of $M$
and all other points
correspond to  primitive CR structures.
\endproclaim
The explicit description of all non-standard CR structures on these manifolds
is given in Corollary 5.2 and Propositions 6.3  and 6.4.\par
For what concerns the non-primitive and non-standard CR structures,
we have the following theorem (for an explicit description
of the CR structures,
see also Theorem 5.1 and  Proposition 6.5).\par
\proclaim{Theorem 1.5} Let  $(M = G/L, \D_Z, J)$ be a simply connected
homogeneous CR manifold
with a non-standard non-primitive CR structure. \par
Then, either $M = SU_2$ or
there exists a unique CRF 
fibration
$$\pi\: M = G/L \longrightarrow F = G/Q$$
over a flag manifold $F$ with an invariant complex structure $J_F$, 
such that the fiber $C = Q/L$ is either 
a primitive CR manifold or is equal to $SO_3 = S(S^2)$. 
Moreover
the  groups $G$, $L$, the primitive fiber $C = Q/L$ and the flag manifold
$F = G/Q$ are  as in the 
following table
(in  n.2, the subgroups $U_{p-2}$ and $U'_{q-2}$ of $L$ are 
subgroups of the  factors $SU_p$ and $SU'_q$ of $G$, respectively):\par
\medskip
\moveleft 0.3cm
\vbox{\offinterlineskip
\halign {\strut\vrule\hfil\ $#$\ \hfil
 &\vrule\hfil\ $#$\
\hfil&\vrule\hfil\ $#$\
\hfil&\vrule\hfil\  $#$
\hfil
&\vrule\hfil\  $#$
\hfil
\vrule\cr
\noalign{\hrule} n^o &
\phantom{\frac{\frac{1}{1}}{\frac{1}{1}}}G
\ \ &
L
&
C = Q/L
& F = G/Q
\cr \noalign{\hrule}
1 &
SU_{n}\ \smallmatrix n >2 \endsmallmatrix
& T^1 \cdot SU_{n-2}
& 
SO_3 = S(S^2)
&
\underset{\phantom{A}}\to
{\overset{\phantom{B}}\to{
\frac{SU_n}{S(U_2 \cdot U_{n-2})}
}}
\cr \noalign{\hrule}
2 &
\underset{\phantom{B}}\to{
\overset{\phantom{A}}\to{\matrix SU_p\times SU'_q\\
\smallmatrix p + q > 4\endsmallmatrix
\endmatrix}}
&
T^1\cdot U_{p-2} \cdot U'_{q-2}
&
\frac{SO_4}{SO_2} = S(S^3)
& \frac{SU_p}{S(U_2\cdot U_{p-2})}\times \frac{SU_q} 
{S(U_2\cdot U_{q-2})}
\cr \noalign{\hrule}
3 &
SU_n\ \smallmatrix n >4 \endsmallmatrix
&
T^1\cdot (SU_2 \times SU_2) \cdot SU_{n-4} 
&
\frac{SO_6}{SO_4} = S(S^5)
&
\underset{\phantom{A}}\to
{\overset{\phantom{B}}\to{
\frac{SU_n}{S(U_4\times U_{n-4})}
}}
\cr \noalign{\hrule}
4 &
SO_{10}
&
T^1\cdot SO_6
&
\frac{SO_8}{SO_6} = 
S(S^{7})
&
\underset{\phantom{A}}\to
{\overset{\phantom{B}}\to{
\frac{SO_{10}}{T^1 \cdot SO_8}
}}
\cr \noalign{\hrule}
5 &
E_6
&
T^1\cdot SO_8
&
\frac{SO_{10}}{SO_8} = S(S^9)
&
\underset{\phantom{A}}\to
{\overset{\phantom{B}}\to{
\frac{E_6}{T^1\cdot SO_{10}}
}}
\cr \noalign{\hrule}
}}
 \par
\medskip
\noindent In particular,  the  fiber $C$ is  
a  sphere 
bundle $S(S^r) \subset TS^r$ where $r = 2, 3, 5, 7$ or $9$.
 The CR manifolds in $n.1$ admit also a CRF fibration with fiber $S^1$.
\endproclaim
\medskip
\proclaim{Corollary 1.6} Let $\pi\: M= G/L \to F = G/Q$ be the 
CRF fibration of a non-primitive non-standard CR manifold $(G/L, \D, J_0)$
onto the flag manifold $F = G/Q$ with a fixed
invariant complex structure $J_F$.
Then the set of all invariant CR structures $(\D, J)$ on $G/L$
(up to sign of $J$), such that the fibering $\pi\: M = G/L \to 
F = G/Q$ is holomorphic, is  
parameterized by the points of the unit disc $D$ in $\R^2$. The center of $D$
corresponds to the unique standard CR structure $J_s$ of this 
family and all other points
correspond to non-standard CR structures.
\endproclaim
\medskip
The unique standard CR structure $J_s$ on $M= G/L$ such that the 
fibration  $\pi\: M= G/L \to F = G/Q$ is holomorphic w.r.t. $J_s$
and $J_F$ is called  {\it the standard CR structure associated with 
the non-standard CR structure $J_0$\/}.\par
\medskip
Finally we give the description of all  non-primitive CR manifolds $G/L$
of a given compact Lie group $G$ 
 in terms of {\it painted Dynkin graphs of $\g = Lie(G)$\/}, 
that is of  Dynkin graphs of the   Lie algebra $\g$
with nodes painted in three colors: white ($\circ$), 
black ($\bullet$) and 'grey' ($\otimes$). \par
Recall that any flag manifold
$F = G/Q$ with an invariant complex structure $J_F$ is  
defined (up to 
equivalence)  by a  black-white
Dynkin graph, where the subalgebra ${\goth q} = Lie(Q)$ is
generated by the Cartan subalgebra and the root
vectors associated with the white nodes. The complex structure $J_F$
is determined by the decomposition
$$\g^\C  = {\goth q}^\C + \m^{10} + \m^{01}$$
where $\m^{10}$ is the nilpotent subalgebra generated by the root
vectors  associated to  black nodes (see e.g. [3], [4]). \par
With a painted Dynkin graph $\Gamma$ (equipped by simple
roots in a standard way), we associate two flag manifolds
$F_1(\Gamma) = G/K$  and $F_2(\Gamma) = G/Q$ and 
two invariant complex structure
$J_1(\Gamma)$ and $J_2(\Gamma)$
on $F_1(\Gamma)$ and $F_2(\Gamma)$, respectively, as follows. 
The pairs  $(F_1(\Gamma) = G/K, J_1(\Gamma))$
and $(F_2(\Gamma) = G/Q, J_2(\Gamma))$
are the flag manifolds with invariant complex structures 
defined by the black-white graphs obtained from $\Gamma$
by considering the grey nodes as black and, respectively, white.\par 
Note that $Q$ contains $K$ and that the natural fibration 
$$\varpi\: F_1(\Gamma)  = G/K \to F_2(\Gamma) = G/Q$$
is holomorphic and a fiber $Q/K$ is a flag manifold
with an induced invariant complex structure $J'$. Moreover, 
$J_1(\Gamma)$ is canonically defined by $J_2(\Gamma)$ and $J'$.\par
Conversely, if $F_1= G/K$ and $F_2 = G/Q$ are two flag manifolds
with invariant complex structures $J_1$ and $J_2$ such that 
$Q\supset K$  and the equivariant fibration $\varpi\: F_1 \to F_2$
is holomorphic, then we may associate with $F_1$ and $F_2$ a painted Dynkin
graph in an obvious way.\par
\medskip
\definition{Definition 1.7} A {\it  CR-graph\/} is 
 a pair $(\Gamma, \vartheta(\Gamma))$, formed by a 
painted Dynkin graph $\Gamma$  and
 a linear combination  $\vartheta(\Gamma)$ of simple
roots,  given in the following table:
\par
\medskip
\moveleft 0.2 cm
\vbox{\offinterlineskip
\halign {\strut\vrule\hfil\ $#$
\ \hfil
 &\vrule\hfil\ $#$\ \hfil
 &\vrule\hfil\ $#$\
\hfil&\vrule\hfil\ $#$\
\hfil
\vrule\cr
\noalign{\hrule} type & \g &
\phantom{\frac{\frac{1}{1}}{\frac{1}{1}}}\Gamma
\ \ &
\vartheta(\Gamma)
\cr \noalign{\hrule}
I &
A_n{}^{\phantom{a}}_{\phantom{a}}\ \ 
\smallmatrix (n >1) \endsmallmatrix&
\underset{\phantom{a}}\to
{\overset{\phantom{a}}\to{
\dynkin\xroot{}\link\root{}\link
\wroot{}\link\dots\link\wroot{}
\enddynkin
}}
 &
\varepsilon_1 - \varepsilon_2
\cr \noalign{\hrule}
II & \matrix
A_p + A'_q{}^{\phantom{a}}_{\phantom{a}}\\ 
\smallmatrix (p + q >2)\endsmallmatrix
\endmatrix
 &
\underset{\phantom{a}}\to
{\overset{\phantom{a}}\to{
\matrix
\dynkin\xroot{}\link\root{}\link\wroot{}\link\dots \link\wroot{}\link\wroot{}
\enddynkin \\
\dynkin\xroot{}\link\root{}\link\wroot{}\link\dots\link\wroot{}\link\wroot{}
\enddynkin 
\endmatrix
}}
&
(\varepsilon_1 - \varepsilon_2) 
- (\varepsilon'_1 - \varepsilon'_2)
\cr \noalign{\hrule}
III & A_n{}^{\phantom{a}}_{\phantom{a}}\ \ 
\smallmatrix (n >3) \endsmallmatrix
&
\underset{\phantom{a}}\to
{\overset{\phantom{a}}\to{
\dynkin\wroot{}
\link\xroot{}\link\wroot{}
\link\root{}\link\wroot{}\link\dots\link\wroot{}\enddynkin
}}
& \varepsilon_1 + \varepsilon_2 - \varepsilon_3 - 
\varepsilon_4 \cr \noalign{\hrule}
IV & D_5{}^{\phantom{a}}_{\phantom{a}}\ \ 
\smallmatrix  \endsmallmatrix
&
\underset{\phantom{a}}\to
{\overset{\phantom{a}}\to{\dynkin
\root{}\link\wroot{}
\link\wroot{}\xrootupright{}
\wrootdownright{}\enddynkin}} &
\varepsilon_{n-3}  + 
\varepsilon_{n-2}  + \varepsilon_{n-1}  
+\varepsilon_n
\cr \noalign{\hrule}
V & E_6 & 
\underset{\phantom{a}}\to
{\overset{\phantom{a}}\to{
\dynkin\xroot{}\link\wroot{}\link\wroot{}\wrootdown{}\link\wroot{}
\link\root{}
\enddynkin
}}
&
2 \varepsilon_1 
 +\varepsilon_6
+\varepsilon
\cr \noalign{\hrule}
}}
\enddefinition
The correspondence between nodes and simple roots is 
as in Table 4 of the Appendix.\par
\medskip
The CR-graphs of type I are called {\it special CR-graph\/}. All the others are 
called {\it non-special CR-graphs\/}.\par
Let $(\Gamma, \vartheta(\Gamma))$ be a CR-graph. We fix
a Cartan subalgebra $\h$ of the associated compact Lie algebra $\g$
and define the element 
$Z(\Gamma) = i\B^{-1}\circ \vartheta(\Gamma)\in \h$. Then $Z(\Gamma)$
is a contact element and we call the corresponding contact 
manifold $(M(\Gamma) =
G/L, \D_{Z(\Gamma)})$  {\it the contact 
manifold associated with the CR-graph $(\Gamma, \vartheta(\Gamma))$\/}.
Note that $M(\Gamma)$ is special if and only 
if the CR-graph is special.\par
\medskip
Using the concept of CR-graph, the results of 
our classification may be stated as follows.\par
\bigskip
\proclaim{Theorem 1.8} Let $M =G/L$ be a simply connected,
homogeneous CR manifold
with a non-primitive, non-standard CR structure $(\D_Z, J)$. Suppose also
that $M \neq SU_2$.\par
 Denote by  
$\pi\: G/L \to F_Z = G/K$ the natural (non-holomorphic)
fibration associated with the contact structure $\D_Z$
and by $\pi'\: G/L \to F_2 = G/Q$ 
the unique CRF fibration over a flag manifold
$F_2 = G/Q$ with invariant complex structure $J_2$, with non-standard 
fiber $Q/L$ of minimal dimension, 
which is either primitive or admitting a CRF fibration with fiber
$S^1$. \par
Then $Q\supset K$ and the sequence of fibering
$$M = G/L \longrightarrow F_Z = G/K \longrightarrow F_2 = G/Q$$
is holomorphic with respect to the standard CR structure $(\D, J_s)$ on $M$,
associated to $(\D,J)$, the 
corresponding complex structure $J_s$ on $F_Z$ and the complex structure
$J_2$ on $F_2$. \par
Moreover, the painted Dynkin graph $\Gamma$  associated to the flag manifolds
$F_1 = F_Z$, $F_2$ with  complex structures 
$J_1 = J_s$ and $J_2$, respectively, is a CR graph
and (up to a transformation from the Weyl group) 
$Z$ is proportional to $Z(\Gamma)$.\par
Conversely, if $\Gamma$ is a CR-graph, 
then there exists a unique homogeneous contact manifold
$(M= G/L, \D_Z)$ such that $Z = Z(\Gamma)$ and 
$F_Z = F_1(\Gamma) = G/K$. The complex structure
$J_1(\Gamma)$ defines the unique standard CR structure $(\D_Z, J_1(\Gamma))$
on $M$ such that the sequence of fibrations 
$$M = G/L \longrightarrow F_Z = F_1(\Gamma) = G/K \longrightarrow
 F_2(\Gamma) = G/Q$$
is holomorphic w.r.t. $(\D_Z, J_1(\Gamma))$, $J_1(\Gamma)$ and
$J_2(\Gamma)$. The space of the invariant CR structures $(\D_Z, J)$ on 
$M$ such that the projection $\pi'\: M \to F_2(\Gamma)$ is holomorphic,
is parameterized by the points of the unit disc $D\in \R^2$.
The center of $D$ corresponds
to the CR structure $(\D_Z,J_1(\Gamma))$ and the other points correspond
to the non-standard CR structures. Moreover a CR structure
is non-standard if and only if it induces a non-standard CR structure
on the fiber $Q/L$; such induced CR structure is always primitive, 
with the  exceptions of the cases in which $\Gamma$ is a special CR-graph.
\endproclaim
\bigskip
As final remark, we would like to mention that our classification 
of compact homogeneous CR manifolds have several important corollaries
concerning  compact cohomogeneity
one K\"ahler manifolds. In particular, such corollaries
are an essential tool towards   
the  classification of   K\"ahler-Einstein 
manifolds in the above class. 
They will be discussed in a forthcoming paper. 
\bigskip
\subhead 2. Basic facts about CR structures
\endsubhead
\bigskip
\definition{Definition 2.1}
\roster
\item A  CR structure on a manifold $M$ is a
pair $(\Cal D , J)$, where $\Cal D \subset TM$ is a distribution on
$M$ and $J \in \operatorname{End}\Cal D , \, J^2 = -1$, is a complex
 structure on $\Cal D$.
\item  A CR structure $(\Cal D, J)$ is called {\it integrable\/} if $J$
satisfies the following integrability condition:
$$J([JX, Y] + [X, JY])\in \D\ ,$$
$$[JX, JY] - [X, Y] - J([JX, Y] + [X, JY])  = 0\ \tag2.1$$
for any pair of vector fields $X$, $Y$ in $\D$.\par
\endroster
\enddefinition
  In the sequel, by CR manifold  we will  understand  a manifold $M$
with {\it integrable} CR structure.\par
If $(\Cal D , J)$ is a CR structure then the complexification
$\Cal D^{\Bbb C} \subset T^{\Bbb C}M$ of the distribution $\Cal D$ is
decomposed into a sum
$\Cal D^{\Bbb C} = \Cal D^{10} + \Cal D^{01}$  of two mutually 
conjugated $( \D^{10} =\bar \D^{01} ) $ 
$J$-eigendistributions with eigenvalues $i$ and $-i$. The 
integrability condition (2.1) means that these eigendistributions are
involutive (i.e. closed under the Lie bracket).\par
The {\it codimension} of a CR structure $(\Cal D, J)$
is defined as the codimension of the distribution $\Cal D$.
Remark that a codimension zero CR structure is the same as a complex
structure on a manifold.
 A codimension one CR structure $(\Cal D, J)$ is also called   a
CR structure  of {\it hypersurface type\/}, because such is the structure 
which is
induced on a real hypersurface of a complex manifold. In this case
the distribution $\Cal D$ can be described locally as the kernel of
a 1-form $\theta$. The form $\theta$ defines an Hermitian symmetric 
bilinear form 
$${\Cal L}^\theta_q\: \D_q\times\D_q \to \R$$
given by
$${\Cal L}^\theta( v, w) = (d\theta)(v, Jw) $$
for any $v,w \in \Cal D$. It is called the {\it Levi form\/}.
Remark that the 1-form $\theta$ is defined up to  multiplication 
by a function $f$ everywhere different from zero and
that ${\Cal L}^{f\theta} = f{\Cal L}^{\theta}$. In particular, the conformal
class
of a Levi form depends only on the CR structure. 
\par
A  CR structure $(\D, J)$ of hypersurface type is called
{\it non-degenerate} if it has non-degenerate Levi form or, in other
words, if $\D$ is a contact distribution. In this case a 1-form
$\theta$ with $\ker \theta = \D$ is called {\it contact form}.
\par
  A smooth map
$\varphi\: M \to M'$ of one  CR manifold $(M, \D, J)$
into another  one  $(M', \D', J')$ is called
 {\it  holomorphic map  \/} if \par
a) $\varphi_*(\D) \subset \D'$;\par
b) $\varphi_*(Jv) = J'\varphi_*(v)$ for all $v\in \D$.
\par
 In particular, we may speak about CR {\it  transformation} of a CR 
manifold $(M, \D, J)$ as a transformation $\varphi$ such that
$\varphi$ and $\varphi^{-1}$ are CR maps. In general, the group of all
CR transformations is not a Lie group, but it is a Lie group
when $(\D, J)$ is of hypersurface type and it is Levi non-degenerate.
\par
\definition{Definition 2.2}
  A CR manifold $(M,\D,J)$ is called {\it homogeneous\/} if it admits a 
  transitive Lie group $G$ of CR transformations. 
\enddefinition  
 \par
 Our aim is to classify compact homogeneous codimension one
non-degenerate CR manifolds. The following theorem, which is 
indeed a consequence of the results in [1], 
shows that we 
may identify any such manifold with a quotient space $G/L$ of a compact
Lie group $G$.  
\par         
\medskip
\proclaim{Theorem 2.3} [12]  Let $(M,\D,J) $ be a compact non-degenerate  
CR manifold of hypersurface type.
Assume that it is homogeneous, i.e. that there exists a transitive Lie group
$A$ of CR transformations. Then a maximal compact connected subgroup
$G$ of $A$ acts on $M$ transitively  and one may identify $M$ with 
the quotient space $G/L$ where $L$ is the stabilizer of a point $p 
\in M$.
\endproclaim
\medskip
 \par
Now we fix some notations.
  If the opposite is not stated, we will assume that a CR structure 
  is of hypersurface type, integrable and Levi non-degenerate.\par
The Lie algebra of a Lie group is denoted by the corresponding
gothic letter.\par
  For any subset $A$ of a Lie group $G$ or of its Lie algebra $\frak g$,
we denote by $C_G(A)$ and $C_{\g}(A)$ its centralizer in $G$ and 
$\g$,
respectively.
$Z(G)$ and $Z(\g)$ denote the center of a Lie group $G$ and a Lie 
algebra $\g$.
 By homogeneous manifold $M= G/L$ we mean a homogeneous manifold
of a compact connected Lie group $G$ with {\it connected stability 
subgroup $L$\/} and such that the action of $G$ on $M$ is {\it effective\/}.
\bigskip
\bigskip
\subhead 3. Compact Homogeneous Contact Manifold
\endsubhead
\bigskip
\subsubhead 3.1 Homogeneous contact manifolds of a compact
Lie group $G$
\endsubsubhead
\medskip
 Let $M=G/L$ be a homogeneous manifold of a compact Lie group $G$ 
 with connected stabilizer $L$. 
\par
 A 1-form $\theta \in \g^*$ on the Lie algebra  $\g$ of $G$ is 
 called {\it contact form} if it is $\operatorname {Ad}_{\l}$-invariant and 
 vanishes on $\l = \operatorname { Lie}L$.  Such form defines a global invariant
1-form $\theta$ on the manifold $M$ which is a contact form of the
contact distribution $\D = \ker \theta$. This establishes a 1-1
correspondence between invariant contact structures $\D$ on $M$
and contact 1-forms $\theta \in \g^*$ up to a scaling (see e.g.[2]).
\par 
  Fix now an $\operatorname {Ad}_G $-invariant Euclidean metric $\B$ on $\g$
and denote by $\l^{\perp}$ the orthogonal complement to $\l$ in $\g$.
\par
 The vector $Z= \B^{-1}\circ \theta$ which corresponds to a contact 
 form $\theta$ is called a {\it contact element} of the manifold 
$M= G/L$.
\par
 It is characterized by the properties that: 
\roster
\item $Z \in \l^{\perp}$ and
\item  the centralizer $C_{\g}(Z) = \l \oplus \Bbb R Z$.
\endroster
  Hence, we have the following
\medskip
\proclaim{Proposition 3.1} There exists a natural bijection between  
invariant contact structures on a homogeneous manifold $M =G/L$ and 
contact elements $Z$ defined up to a scaling. 
\endproclaim 
\medskip
 We will denote by $\D_Z$ the contact structure on $M$  defined 
 by a contact element $Z$. A homogeneous manifold $M=G/L$ with an
invariant contact structure $\D$ is called {\it homogeneous contact
manifold\/}.
\par
\bigskip 
 Proposition 3.1  implies the following 
\proclaim{Corollary 3.2}
Let $G/L$ be a homogeneous contact  manifold of a compact Lie group $G$ which 
acts effectively.  Then the the center $Z(G)$ of $G$ has dimension 0 
or 1.
\par
Moreover, if  $Z(G)$ is one dimensional, 
then any contact element $Z$ has not zero orthogonal projections 
$ Z_{Z(\g)} ,\, Z_{\g'}$ on $Z(\g)$ and $\g' = [\g,\g]$, and
the stability subalgebra $\l$ can be written as
$$\l = [C_{\g'}(Z_{\g'})]_\varphi \=
 \{ X = Y + \varphi(Y)\ ,\ Y \in C_{\g'}(Z_{\g'})\}$$
where $\varphi\: C_{\g'}(Z_{\g'}) \to Z(\g)\approx \Bbb R $ is a non-trivial
 Lie algebra homomorphism. 
\endproclaim
\demo{Proof} Clearly  $C_\g(Z) \supset
Z(\g)$. If $\dim Z(\g) \geq 2$
then  $\l\cap Z(\g) \neq \{0\}$ and this contradicts
the fact that $G$ acts effectively. 
The other claims follow immediately.\qed
\enddemo
\bigskip
  Now we associate with a homogeneous contact manifold $(M=G/L,D_Z)$
a flag manifold  
$$F_Z \= \operatorname{Ad}_G Z = \operatorname{Ad}_{G'}(Z_{\g'})$$
where $K= C_G(Z)$ is the centralizer of the contact element $Z$.
We will call
$F_Z$ {\it the flag manifold associated to a contact element\/}
$Z$.\par
 Note that the contact
form $\theta = \B\circ Z$ 
is a connection (form) in the $S^1$ bundle $\pi: 
G/L \to F_Z$ and that the corresponding 
contact structure  
$\D =\operatorname{ker}\theta$ is the horizontal distribution  of  
this connection.\par
 \bigskip
 We describe now all homogeneous contact manifolds $(G/L, \D_Z)$
with given  associated flag manifold $F = G/K$ of a semisimple Lie group 
$G$.
\par
 Consider the orthogonal reductive  decomposition
$$\g = \k + \m 
$$
associated with the flag manifold $F=G/K$.
\par
 We say that an element $Z$ of the center $Z(\k)$ is {\it $\k$-regular} if
it generates a closed 1-parametric subgroup of $G$ and  the 
centralizer $C_G(Z)= K$. 
\par
One can check that if $Z$ is $\k$-regular, then the 
subalgebra 
$$\l_Z = \k \cap (Z)^\perp$$
generates a closed subgroup, which we denote by $L_Z$. 
 Therefore
 \par
\proclaim{Proposition 3.3} Let $F=G/K$ be a flag manifold of a semisimple 
Lie group $G$.
 There is a natural 1-1 correspondence \par
$$ Z \Longleftrightarrow (G/L_Z, \D_Z) $$
between the $\k$-regular elements
$Z \in \g$ (determined up to a scaling) 
and the homogeneous contact
manifolds $(G/L, D)$  
associated  flag manifold $F = G/K$.
\endproclaim
\demo{Proof} The proof is straightforward.
\qed
\enddemo
\bigskip
\subsubhead 3.2 Invariant contact structures on a contact 
manifold $M = G/L$
\endsubsubhead
\medskip
 Now we describe  all invariant contact structures  on a given
homogeneous manifold $M=G/L$. We will show that generically there is
no more then one such structure.
\par
\medskip 
\definition {Definition 3.4} A homogeneous  manifold $G/L$ is called 
homogeneous contact manifold of {\it non-special type} (respectively, 
  of \,{\it  special type } or, shortly, {\it special})
   if it admits a unique (respectively, more then one)
invariant contact structure.
\enddefinition 
\bigskip
\subsubhead 3.2.1 Main  examples of 
special homogeneous contact manifolds
\endsubsubhead
\medskip
 Let $\g$ be a compact semisimple Lie algebra, $\h $  a Cartan
subalgebra of $\g$ and  $R$  the root system  of the  pair 
 $(\g^{\C}, \h^\C)$.
\par
Recall that a root $\alpha \in R$  
defines a 3-dimensional regular subalgebra $\g^{\C}(\alpha)$ with  
standard basis given by the root vectors
$E_\alpha, E_{-\alpha}$ and 
$$H_\alpha = [E_{\alpha}, E_{-\alpha}] = \frac{2}{|\alpha|^2} \Cal B^{-1}\circ
\alpha
\tag 3.1$$ 
verifying
 the relation $[H_{\alpha}, E_{\pm\alpha}] = \pm 2 E_{\pm\alpha}$.
Its intersection with $\g$ is a 3-dimensional compact 
subalgebra denoted by $\g(\alpha)$. We will call  $\g(\alpha)$ the
{\it subalgebra associated 
with the root $\alpha$\/} and denote by $G(\alpha)$
the 3-dimensional subgroup of the adjoint group 
$G = \operatorname {Int}(\g)= \operatorname {Aut}(\g)^0$ generated by
$\g(\alpha)$. \par
Note that any two such subalgebras 
are conjugated by an inner 
automorphism of $\g$ if and only if  the corresponding roots have the same  
length.\par
\medskip
 Fix a system $R^+$ of positive roots of $R$ and put $R^- =-R^+$.
 The highest root  $\mu$ of $R^+$ 
defines the following gradation of the complex Lie algebra  $\g^{\C}$:
$$\g^\C = \g_{-2} + \g_{-1} + \g_0 + 
\g_1 + \g_2\ ,\tag3.2$$
where 
$$\g_{-2} = \C E_{-\mu}\qquad
\g_2 = \C E_{\mu}\qquad \g_0 = \C H_\mu + \g'_0\qquad
\g_0' = C_{\g^\C}(\g(\mu))\tag3.3$$
$$\g_{-1} = \sum_{\beta \in R^-\setminus (\{ - \mu\}\cup R_o)} 
\C E_\beta\qquad
\g_{1} = \sum_{\beta \in R^+\setminus (\{\mu\}\cup R_o)} \C E_\beta$$
and $R_o = \{ \alpha \in R, \, \alpha \perp \mu  \}$ is the root 
system of the subalgebra $\g_0 = C_{\g}(H_{\mu})$.
\par
(3.2) is called the {\it gradation associated with the highest root\/}.
\par
The explicit decomposition (3.2) for any simple complex Lie algebra
is given in Table 1 of the Appendix. 
\bigskip
  Denote by 
$\l = C_{\g}(\g(\mu)) = \g_0' \cap \g$ the centralizer of $\g(\mu)$ 
in $\g$ and by $L$ the corresponding connected subgroup of $G$. It is
easy to check that $L= C_G(\g(\mu))$.
\par
\proclaim{Lemma 3.5}
  Let $G$ be a compact simple Lie group without center and 
let $L = C_G(\g(\mu))$ be as defined  above.
 Then any non zero vector $Z \in \g(\mu)$ is a contact element of
the manifold $G/L$. In particular, $G/L$ is a homogeneous contact 
manifold of special type.
\endproclaim
\demo{Proof}
Observe that $Z\in \g(\mu)$ is a contact element if and only if
$C_\g(Z) = \l + \R Z$  and then   $g\cdot Z$ is contact for any
 $g\in G(\mu)$.  Since $G(\mu)$
acts transitively on the unit sphere 
of $\g(\mu)$, the Lemma follows from the fact that 
$$C_\g(iH_\mu) = \g_0 \cap \g = \l + \R(iH_\mu)$$
and hence that $i H_\mu$ is a contact element. \qed
\enddemo
 Remark that the contact manifolds $M= G/L = G/C_G(\g(\mu)) $, with  
 $G$ simple, carry invariant 3-Sasakian structure and they exhaust all
homogeneous 3-Sasakian manifolds (see [6]).   
\bigskip
\subsubhead 3.2.2  Classification of special homogeneous contact 
manifolds
\endsubsubhead
\medskip
The previous examples almost  exhaust the class of special 
homogeneous contact  manifolds. In fact, we have the following
classification theorem.\par
\bigskip
\proclaim{Theorem 3.6} Let $M = G/L$ be a special homogeneous 
contact manifold of a compact Lie group $G$. Then
the group $G$ is simple and either $L$ is the centralizer of
the subalgebra $\g(\mu)$ associated with the highest root and $M$ is
 a homogeneous 3-Sasakian manifold or $G=G_2$ and $L$ is the 
 centralizer of the subalgebra $\g(\nu)$ associated with a short root 
$\nu$. 
\endproclaim 
\demo{Proof}
 We prove first that if $G$ is not semisimple and, hence, 
 $\dim Z(\g)=1$, then a contact element $Z$ is unique up to a scaling
and $M$ is not special. Indeed, we have the decomposition
$$ \k = C_{\g}(Z) = \l \oplus \Bbb R Z = \l + Z(\g) $$
since $Z(\g)\cap \l =0$, by effectivity.
The line $\Bbb R Z$ is determined uniquely as the orthogonal 
complement to $\l$ in $\frak k = \l + Z(\g)$.
\par
  Now we may assume that $\g$ is semisimple. We need the following: 
\par
\bigskip
\proclaim{Lemma 3.7}
Let $\g$ be compact semisimple
and let $\l\subset \g$ be a closed subalgebra,
which contains no ideal of $\g$. If there exist two not 
proportional vectors $Z$, $Z'\in \l^\perp$
such that
$$C_{\g}(Z)=\l + \Bbb RZ, \quad \l +\Bbb RZ'\subseteq C_{\g}(Z')\ ,$$ 
then $\g$ is simple and
there exists a root $\alpha\in R$ such that:
\roster
\item $\l = C_\g(\g(\alpha))$;
\item $Z, Z' \in \g(\alpha)$ and 
$C_\g(Z') = C_\g(\g(\alpha)) + \R Z'$;
\item $C_\g(\l) = Z(\l) + \g(\alpha)$; 
\item for any root $\beta$ which is orthogonal to 
$\alpha$, $\alpha \pm \beta$ is not a root.
\endroster
\endproclaim
\demo{Proof}
We put  $\k = C_{\g}(Z)$  and consider the orthogonal decomposition
$$ \g = \k +\m =(\l +\Bbb RZ)+ \m\ . $$ 
 Denote by  $R$ the root system of the complex Lie algebra
$\g^{\Bbb C}$ with respect to a Cartan subalgebra $\h^{\Bbb C}$
which is the complexification of a Cartan subalgebra $\h$ of $\k$.
Then the element $Z'$ can be written as 
$$ Z' = cZ + \sum_{i=1}^k c_iE_{\alpha_i} $$
for some root vectors $E_{\alpha_i}$ and constants $c, c_i$. 
The condition $[\l, Z'] = 0$ implies
$\alpha_i(\h\cap \l)= 0 $ if $c_i \neq 0$. Since $\h\cap \l$ is of codimension
one in $\h$, there exist exactly  two (proportional) roots with this
 property, say
$\alpha$ and $-\alpha$. This shows that 
$\l \subset C_\g(\g(\alpha))$. Moreover,
 since $Z\in \h\cap \l^\perp$, we obtain also that $Z$
is proportional to $H_\alpha= [E_{\alpha}, E_{-\alpha}]$
and (1) follows. In particular, $\g$ must be simple
and now (2) is  clear.  (3) follows from (2).\par 
To prove (4), 
assume that there is a root $\beta$ which is orthogonal to $\alpha$ 
and such that $\alpha + \beta $ is a root. Then the vector 
$E_{\beta} + E_{-\beta} \in \g^\C $ does not belong to  
$\l^\C =C_{\g^\C}(\g(\alpha))$, 
but it is orthogonal to $Z$ 
(since  $Z$ is proportional to  
$H_\alpha$) and belongs to the 
centralizer of $Z$: contradiction.
\qed 
\enddemo 
 \bigskip
 Now we conclude the proof of Theorem 3.6. Let $G$ be a compact
semisimple Lie group and let $Z$, $Z'$  two  non-proportional contact 
elements for $G/L$. By Lemma 3.7, $G$ is simple
and $L = C_G(\g(\alpha))$.
 By direct inspection of the root
systems of  simple Lie groups,  a root $\alpha$
verifies the condition (4) of Lemma 3.7 if and only if it is 
a long root or if  it is a short root in the $G_2$ type system.
This concludes the proof. \qed
\enddemo
\bigskip
\bigskip
\subsubhead 3.3 Isotropy representation of a homogeneous
contact manifold
\endsubsubhead
\medskip
Let $M= G/L$ be a homogeneous contact manifold with invariant 
contact structure $\D$ associated to a contact element $Z$. Let
$\g = \l + \R Z + \m$ be the corresponding orthogonal decomposition.
Fix a Cartan subalgebra $\h$ of $\g$ which belongs to
$\k = \l + \R Z = Z(\k) + \k'$
(where 
$\k' = [\k, \k]$ is the semisimple
part of $\k$). Then
$$\h = Z(\k) + \h' = Z(\l) + \R Z + \h'\ ,$$
where we denote by $\h'$ a Cartan subalgebra of $\k'$. 
Remark that $\h(\l) = Z(\l) + \h'$ is a Cartan
subalgebra of $\l$. \par
Denote by $R$ (resp. $R_o$) the root system of $\g^\C$
(resp. $\k^\C$) w.r.t. the Cartan subalgebra $\h^\C$ and let
$R' = R \setminus R_o$. We will denote by $\h(\R)$ the standard
real form of $\h$, spanned by $R$, that is
$$\h(\R) = \h \cap \B^{-1}(<R>)\ .$$
We put  $\goth t = \z(\k) \cap \h(\R)$. Then $Z \in i \goth t$ and 
we may identify 
$$\vartheta = -i\theta = -i \B(Z, \cdot)$$
 with the corresponding
element in $\goth t^* \subset \h(\R)^* = \text{span}_\R R$.\par
\bigskip
Consider the decomposition of the $\k^\C$-module $\m^\C$
into sum of irreducible $\k^\C$-modules
$$\m^\C = \sum \m(\gamma)\ .\tag3.4$$
Here, $\m(\gamma)$ stands for the irreducible $\k^\C$-module
with highest weight $\gamma \in R'$.\par
\bigskip
The following Lemma states a well known property of flag manifolds
(see e.g. [4] or [3]).
\medskip
\proclaim{Lemma 3.8} The $\k^\C$-modules $\m(\gamma)$ are pairwise
not equivalent and, in particular, the decomposition (3.4) is unique.
The moduli $\m(\gamma)$ are irreducible also as $\l^\C$-modules. 
\endproclaim
\demo{Proof} We only need to check that a module $\m(\gamma)$
is irreducible also as an $\l^\C$-module. But it is sufficient to observe
that the semisimple parts of $\l^\C$ and of $\k^\C$ coincide and to
recall that, 
whenever 
 $\dim_\C \m(\gamma) > 1$, the semisimple part of 
 $\k^\C$ acts non-trivially and irreducibly on $\m(\gamma)$.\qed
\enddemo 
\bigskip
 From Lemma 3.8 we derive the following technical 
proposition, which will be 
useful in the following sections.\par
\medskip
\proclaim{Proposition 3.9} Let $M = G/L$ be a homogeneous
contact manifold and let $Z$ be a contact element for $M$. Assume
that $G\neq G_2$ or that $G = G_2$ and $\vartheta = -i \B \circ Z$
is not proportional to a short root of $R$. \par
Then for any irreducible $\k^\C$-module $\m(\gamma)$
there exists at most one distinct $\k^\C$-module $\m(\gamma')$
which is isomorphic to $\m(\gamma)$  
as $\l^\C$-module.\par
This is the case if and only if the 
highest weights $\gamma$ and $\gamma'$
are $\vartheta$-congruent, i.e. $\gamma' = \gamma + \lambda \vartheta$
for some real number $\lambda$.
\endproclaim
\proclaim{Corollary 3.10} Let $M$ and $Z$ as in the 
Proposition 3.9. Then:
\roster
\item"a)" if the modules $\m(\gamma)$,
$\m(\gamma')$ are equivalent as $\l^\C$-modules, then for any weight
$\alpha \in R'$ of $\m(\gamma)$, there exists exactly one
weight $\alpha' \in R'$ of $\m(\gamma')$ which is $\vartheta$-congruent
to $\alpha$;
\item"b)" for any root $\alpha \in R'$ there exists at most
one  root $\alpha' \in R'$ which is $\vartheta$-congruent to 
$\alpha$, i.e. such that $\alpha' = \alpha + \lambda \vartheta$
for some real number $\lambda\neq 0$.
\endroster
\endproclaim
\demo{Proof of Proposition 3.9} Observe that two irreducible
$\l^\C$-modules $\m(\gamma)$ and $\m(\gamma')$ are isomorphic if and
only if their highest weights
$\gamma|_{\h(\l)}$ and $\gamma'|_{\h(\l)}$ coincide. This occurs if and
only if $\gamma' = \gamma + \lambda \vartheta$ for some 
$\lambda \in \R$.\par
Assume now that there exist three distinct isomorphic $\l^\C$-modules
$\m(\gamma)$, $\m(\gamma')$ and $\m(\gamma'')$. Then
$\tilde R = \text{span}_\R(\gamma, \gamma', \gamma'')\cap R$
is a 2-dimensional root system and $\gamma$, $\gamma'$ and
$\gamma''$ belong to the straight line $\gamma + \R \vartheta$.
Checking all 2-dimensional root systems, $2 A_1$, 
$A_2$, $B_2$, $G_2$, we conclude that this is possible only if
$\tilde R$ is of type $B_2$ or $G_2$ and $\vartheta$ is proportional 
to a short root. We claim that both these cases cannot occur.\par
If $\tilde R$ has type  $G_2$, then $\tilde R = R$ which contradicts
to the assumptions. 
\par
If  $\tilde R$ has type  $B_2$, one of the roots $\gamma$, $\gamma'$,
$\gamma''$  is orthogonal to $\vartheta$ and this is impossible
because
$$\vartheta^\perp \cap R = R_o = R\setminus R'$$
while $\gamma$, $\gamma'$, $\gamma'' \in R'$.\qed
\enddemo   
\bigskip
\bigskip
\subhead 4. General Properties of Compact Homogeneous CR manifolds 
\endsubhead
\bigskip
\subsubhead 4.1 Infinitesimal description of invariant CR structures
\endsubsubhead
\medskip
Let $(M = G/L, \D_Z)$ be a homogeneous contact manifold  of a 
connected compact Lie group $G$ with connected stabilizer $L$
and let
$$\g = \l + \R Z + \m \tag4.1$$
 the associated orthogonal decomposition where 
$\k = C_\g(Z) = \l + \R Z$.\par
\medskip
\definition{Definition 4.1} A complex subspace
$\m^{10}$ of $\m^\C$ is called {\it holomorphic\/} 
if 
\roster
\item"i)" $\m^{10}\cap \m^{01} = \{0\}$, where $\m^{01} = \overline{\m^{10}}$
and 'bar' denotes 
the complex conjugation with respect to the real subspace $\g$;
\item"ii)" $\m^\C = \m^{10} + \m^{01}$;
\item"iii)" $\l^\C + \m^{10}$ is a complex subalgebra of $\g^\C$.\par
\endroster
In the following we will referee to condition iii) as the {\it
integrability condition\/}.
\enddefinition
Note that if the integrability
condition holds,   also $\l^\C + \m^{01}$ is  a  subalgebra.
Furthermore, any holomorphic subspace $\m^{10}$ defines an 
$\operatorname{ad}_\l$-invariant complex structure $J$ on $\m$, whose $(+i)-$
and $(-i)$-eigenspaces are exactly $\m^{10}$ and $\m^{01}$.
\bigskip
\proclaim{Proposition 4.2}  Let $(M=G/L, \D_Z)$ be a compact
 homogeneous contact
manifold and  $\g = \l + \Bbb R Z + \m$ be the associated 
decomposition. Then there exists a 
natural one to one correspondence   between the set of  invariant
CR structures $(\D_Z, J)$ on $M$ 
and the set of holomorphic subspaces $\m^{10}$ of $\m^\C$.
\endproclaim
\demo{Proof} Recall that,  under the natural 
identification of $\R Z + \m$ with  the 
tangent space $T_{eL}M$, we have that
 $\m = \D_Z|_{eL}$.  Moreover, any invariant CR structure
$(\D_Z, J)$   defines a decomposition $\D^\C_Z =  \D^{10} + \D^{01}$
into two mutually conjugated invariant integrable distributions. 
Then one can easily check that the complex subspace 
 $\m^{10} = \D^{10}_{eL}\subset \m^\C$ is a holomorphic 
subspace. \par
Conversely an holomorphic subspace $\m^{10}$ and its conjugate
subspace $\m^{01} = \overline{\m^{10}}$ are $\operatorname{ad}_\l$-invariant
and also $\operatorname{Ad}_L$-invariant since $L$ is connected. 
Then they can be extended to two invariant integrable complex distributions
$\D^{10}$ and $\D^{01}$ such that $\D^\C = \D^{10} + \D^{01}$ with 
$\D^{10} \cap  \D^{01} = 0$. Hence they may be considered as eigendistributions
of an invariant CR structure $(\D_Z, J)$ on $M$.\qed
\enddemo
 \bigskip
\subsubhead 4.2 Standard CR structures
\endsubsubhead
\medskip
We want to show  how to construct an invariant CR structure $(\D_Z,J)$
on a homogeneous contact manifold $(M = G/L, \D_Z)$ starting from
an invariant complex structure $J$ on the associated flag manifold
$F_Z$.\par
Let $F= G/K$ be a flag manifold and let $\g = \k + \m$ the associated reductive
decomposition. Recall that an invariant complex structure $J_F$ on $F$ is 
associated with a decomposition $\m^\C = \m^{10} + \m^{01}$ such that
$$a)\ \ \m^{01} = \overline{\m^{10}}\ \ ;\quad b)\ \ \p = \k^\C + \m^{10}\ 
\text{is a subalgebra of}\ \g^\C\ .\tag 4.2$$
We say that $\m^{10}$ is {\it the holomorphic subspace associated with $J_F$\/}.
\par
It is known that $\p$ is a parabolic subalgebra, with reductive part
$\k^\C$ and nilradical $\m^{10}$. Moreover, we can always choose a system
of positive roots $R^+$ for $\g^\C$, such that $\m^{10}$ is generated
by root vectors $E_\alpha$, with $\alpha\in R^+$. We say that such system $R^+$
is {\it compatible with the complex structure $J_F$\/}.\par
\smallskip
Let  $(M = G/L, \D_Z)$  be a homogeneous contact manifold,
  $\g = (\l + \R Z) + \m = \k + \m$  the corresponding decomposition 
and   $F_Z = G/K$  the associated flag manifold. Any
invariant complex structure $J_F$ on $F_Z$ 
induces an invariant CR structure 
$(\D_Z, J)$, which is the one corresponding to the same 
holomorphic
subspace  $\m^{10}\subset \m^\C$ as $J_F$.
\medskip
\definition{Definition 4.3}
 An invariant CR structure  $(\D, J)$ on a homogeneous contact 
 manifold $(M=G/L, \D)$, which is induced by an invariant complex
structure $J_F$ on the associated flag manifold $F= G/K$, is called
{\it standard CR structure\/}.
\enddefinition
\medskip
\remark{Remark 4.4} Since any flag manifold admits at least one 
invariant complex structure, we may conclude that 
{\it any homogeneous contact manifold $(G/L, \D)$, with $G$
compact, admits an invariant CR structure $(\D, J)$}.
\endremark
\bigskip
\medskip
The following Lemma gives an algebraic characterization of the 
standard CR structures.
\medskip
\proclaim{Lemma 4.5} An invariant CR structure $(\D, J)$
on a   homogeneous contact   manifold   $(M = G/L, \D)$ 
is standard if and only if
the corresponding complex structure $J$ on $\m$ is 
$\operatorname{Ad}(K)$-invariant.
\endproclaim
\demo{Proof} The proof is straightforward.
 \qed 
\enddemo
\bigskip
Since the description of all invariant complex structures 
on flag manifolds is well known (see  [10], [4],
[5], [3]),  the problem of classification
of the invariant   CR structures  on compact homogeneous
spaces reduces to the description of {\it non-standard\/} invariant
CR structures.\par
\medskip
The following proposition reduces the problem to the case of $G$
semisimple.\par
\medskip
\proclaim{Proposition 4.6} Let $(M = G/L, \D)$ be a contact manifold of a 
compact Lie group $G$
with $\dim Z(G) = 1$. Then any invariant CR structure
with underlying distribution $\D$ is standard.
\endproclaim 
\demo{Proof} It follows immediately from
the fact that any $\operatorname{Ad}(L)$-invariant decomposition 
$\m^\C = \m^{10} + \m^{01}$ is
clearly also  $\operatorname{Ad}(K)$-invariant, since $ K = L\cdot Z(G)$.
\qed
\enddemo 
\bigskip
\subsubhead 4.3 Holomorphic fibering of homogeneous CR manifolds
\endsubsubhead
\medskip
Let $(M = G/L, \D, J)$ be a homogeneous CR manifold with a {\it standard\/} CR
structure $J$
associated to a complex structure $J_F$ on the associated
flag manifold $F = G/K$. Then the natural projection
$$\pi\: G/L \longrightarrow F = G/K$$
is a  $G$-equivariant   holomorphic fibration.
\par
More generally we give the following definition.\par
\bigskip
\definition{Definition 4.7} Let $M = G/L$ be a homogeneous manifold with
invariant CR structure
$(\D, J)$. 
\roster
\item  Any $G$-equivariant holomorphic fibering 
$$\pi\: M = G/L \longrightarrow F = G/Q$$
of $(M, \D, J)$ over a flag manifold $F = G/Q$ equipped with an invariant
 complex structure $J_F$
is called {\it CRF fibration\/};
\item  we say that a homogeneous CR manifold $(M = G/L, \D, J)$ is 
 {\it  primitive \/} if it doesn't admit a non-trivial CRF fibration;
\item a non-primitive  homogeneous CR manifold $(M = G/L, \D, J)$,
admitting a CRF fibration with typical fiber $S^1$, is called
{\it circular\/}.
\endroster
Remark that {\it any standard CR structure is circular\/} and that {\it 
the typical fiber $Q/L$ of a CRF fibration 
carries a natural invariant CR structure\/}.
\enddefinition
\bigskip
The following Lemma  gives a characterization of 
primitive CR structures.\par
\bigskip
\proclaim{Lemma 4.8} A homogeneous CR manifold $(G/L, \D, J)$
admits a non-trivial CRF fibration
 if and only if there exists a proper parabolic
subalgebra $\p = \goth r + \n \subsetneq \g^\C$ (here $\goth
r$ is a   reductive part and $\n$ the nilpotent part)
such that 
$$a)\ \goth r = (\p \cap \g)^\C\ ;\qquad\quad
b)\ \l^\C + \m^{10} \subset \p\ ;\qquad\quad c)\ 
\l^\C\subsetneq \goth r\ .
$$
In this case, $G/L$ admits a CRF fibration with basis $G/Q$,
where $Q$ is the connected subgroup generated by ${\goth  q} = \goth r 
\cap \g$.
\endproclaim
\demo{Proof} Suppose that $(M=G/L, \D, J)$ is non-primitive and 
let $\pi\: G/L \to G/Q$
be a CRF fibration over a flag manifold $F = G/Q$ with invariant complex
 structure $J_F$. Consider the  decompositions associated to $J$ and $J_F$
$$\g = \l + \R Z + \m\ ,\qquad \m^\C = \m^{10} + \m^{01}\ ,$$
$$\g = {\goth  q} + \m'\ , \qquad \m'{}^\C = \m'{}^{10} + \m'{}^{01}\ .$$
Since $\pi$ is holomorphic and non-trivial, the subalgebra 
$\l^\C + \m^{10}$ is properly contained in the parabolic subalgebra
$\p = {\goth  q}^\C + \m'{}^{10}$, with reductive part
${\goth  q}^\C = (\g\cap \p)^\C$. Furthermore, since the fiber has 
positive dimension, 
 $\l\subsetneq
 {\goth  q}$.
\par
  Conversely, if $\p = \goth r + \n \subset \g^\C$ is a parabolic
subalgebra with reductive subalgebra $\goth r = {\goth  q}^\C$, where ${\goth  q}
= \p\cap \g$,
then  we may consider the orthogonal decompositions
$$\g = {\goth  q} + \m'\ ,
\qquad \g^\C = \goth r + \m'{}^\C = \goth r + \n + \n'\ ,$$
where $\n' = \n^\perp \cap \m'{}^\C$. 
By the remarks at the beginning of \S 4.2,
there exists a unique
invariant complex structure $J_F$ with associated 
holomorphic space  $\m'{}^{10} = \n$. Therefore
if  $\l^\C + \m^{10} \subset \p$, $\l\subsetneq {\goth  q}$ and $Q$ is the 
reductive subgroup generated by ${\goth  q}$,
it is clear that $\pi\: G/L \to G/Q$ is  a non-trivial 
CRF fibration.\qed
\enddemo 
\bigskip
\subsubhead 4.4 The anticanonical map of a homogeneous CR manifold 
\endsubsubhead
\bigskip
Let $(M = G/L, \D_Z, J)$ be a homogeneous CR manifolds of a compact
Lie group $G$  and
$$\g = \l + \R Z + \m\quad, \quad \m^\C = \m^{1 0} + \m^{0 1}$$
 the associated  decompositions of $\g$ and of $\m^\C$.
\par
To  characterize the non-standard invariant CR structures we recall
 the definition
of {\it anticanonical map} of a homogeneous CR manifold
introduced for the first time in [1]. It is a $G$-equivariant
holomorphic map 
$$\phi: M = G/L \longrightarrow \text{Gr}_k(\g^\C)$$
into the Grassmanian of complex $k$-planes,  $k = \dim_\C(\l^\C + \m^{01})$,
of $\g^\C$ given by 
$$\phi\: g L \mapsto \operatorname{Ad}_g(\l^\C + \m^{01})\ .$$
Due to the existence of standard holomorphic $G$-equivariant
embedding 
$$\imath\: \text{Gr}_k(\g^\C) \longrightarrow \C P^N\ ,
\quad N = \left(
\smallmatrix \dim \g^\C\\
k \endsmallmatrix
\right) - 1\ ,$$
$$V = \text{span}(e_1, \dots, e_k) \overset{\imath}\to\mapsto
[V] = \C (e_1\wedge \dots \wedge e_k)\ , $$
we may consider $\phi$ as a $G$-equivariant  map into $\Bbb C P^N$.
 To prove that the map $\phi$ is holomorphic it is sufficient to check 
that the linear map 
$$\phi_*\: \D_0 = \ker\theta|_{T_0M}= \m \longrightarrow
T_{[\l^\C + \m^{01}]}\text{Gr}_k(\g^\C)$$
commutes with the complex structure. \par
Let $v = X + \bar X \in \m$, where $X \in \m^{10}$. Then 
$$\phi_*(v) = \operatorname{ad}_{(X+ \bar X)}([\l^\C + \m^{01}])= 
\operatorname{ad}_X([\l^\C + \m^{01}])\ .$$
Therefore
$$ \phi_*(Jv) = \phi_*(iX - i \bar X) = 
\operatorname{ad}_{iX}([\l^\C + \m^{01}]) = 
i \operatorname{ad}_{X}([\l^\C + \m^{01}]) = i\phi_*(v)\ .$$
This shows that the map $\phi$ is holomorphic.\par
   Remark that the stabilizer $Q$ of the point $[\l^\C +\m^{01}]$ in
$\phi(M)= G/Q$ is the normalizer $Q= N_G(\l^\C +\m^{01})$.
\medskip
Now, the following theorem establishes some important
properties of the anticanonical  map.\par
\bigskip
\proclaim{Theorem 4.9} Let 
$$\phi: M = G/L \longrightarrow \text{Gr}_k(\g^\C)$$
be the anticanonical map of a homogeneous CR manifold $(M= G/L,
\D_Z , J)$. \par
\roster
\item If the CR structure is standard, then the image
$\phi(M)$ is $G$-equivariantly biholomorphic to the 
associated flag manifold $F_Z = G/K = Ad_G Z$ endowed with the 
complex structure $J_F$ which induces the CR sructure
$(\D_Z, J)$. \par
In this case, 
$\phi$ is a CRF 
fibration with fiber $S^1$ and the normalizer in $\g$ of $\l^\C + \m^{01}$ is 
$${\goth  k} = N_\g(\l^\C + \m^{01}) =
\l + \R Z\ $$
and it is equal to
the stabilizer of the point $[\l^\C + \m^{01}] \in \phi(M)$
in $G$.
\item If the CR structure is not standard, then the image
$\phi(M) = G/Q$ is a homogeneous CR manifold with CR structure
induced by the complex structure of 
$\text{Gr}_k(\g^\C)$ and $\phi\: M \to \phi(M)$ is a finite covering.
\endroster
\endproclaim
\demo{Proof} We first need the following Lemma, which in fact was
proved in [1].\par
\bigskip
\proclaim{Lemma 4.10} Let $G/Q = \phi(G/L)$ be the image of the 
anticanonical map. Then $\dim Q/L \leq 1$. 
\endproclaim
\demo{Proof}
We need to prove that $\dim {\goth  q}/\l \leq 1$, where 
${\goth  q} = N_{\g}(\l^\C + \m^{01})$ is the stability subalgebra of the 
flag manifold $G/Q$. Since
$\g = \l +\Bbb R Z +\m$,
it is sufficient to check that ${\goth  q}\cap \m =0$. 
Let $ v \in {\goth  q}\cap \m$. Then
$$ \B(Z, [v, \l^\C +  \m^{01}]) \subset  \B(Z,  \l^\C + \m^{01}) = \{0\}
$$
and in particular 
$$ \{0\} = \B(Z, [v,\l + \m]) =  - \B([v,Z], \l + \m)\ .$$ 
 This means that $v \in N_\g(Z) = \k =\l + \R Z$ 
and hence that $v \in \k \cap \m = \{0\}$.
\qed
\enddemo
\medskip
 Let us prove (1).  Notice  that, by Lemma  4.5, if $(\D_Z, J)$ is standard 
then $N_\g(\l^\C + \m^{01}) \supset \l + \R Z$. Therefore, from Lemma
4.10, we get that $N_\g(\l^\C + \m^{01})  = \l + \R Z = \k$ and 
the image $\phi(G/L)$ of
the anticanonical map  coincides with the flag manifold $F=G/K$.
\par
\medskip
For the claim (2), we first show that if the CR structure is
non-standard, then the anti-canonical map
$\varphi: G/L \to \phi(G/L)$ is a finite covering. In fact, if
the CR structure is 
non-circular,  the fiber of the anticanonical map is
not 1-dimensional (otherwise it would give a CRF fibration
with $S^1$-fiber) and by Lemma 4.10  this implies that
 $\varphi: G/L \to \phi(G/L)$ is a finite covering. If the CR structure
is circular and non-standard, then $(M = G/L, \D_Z, J)$
is one of the CR manifolds described in the next subsection \S 4.5. 
In particular, for all these cases
 $\varphi: G/L \to \phi(G/L)$ is  a finite covering (see
the following  Theorem 4.11 and Corollary 5.2). The other
part of the claim  follows immediately by the holomorphicity
and the $G$-equivariance of $\phi$.\qed
\enddemo
\bigskip
\subsubhead 4.5 Circular CR structures which are non-standard 
\endsubsubhead
\bigskip
As we already pointed out, any standard
CR structure is circular. Now we describe 
 the   circular CR structures, which are not standard.\par
 Let 
$(\D, J)$ be a circular CR structure on 
$G/L$ and let   $Z_\D$ be a contact element
associated to $\D$. Let also $\pi: G/L \to G/Q$ be
the CRF fibration onto the flag manifold
$G/Q$ with fiber $S^1 = Q/L$. Notice that, since $\q$ is
the isotropy subalgebra of a flag manifold, $\q$ is of the form
$\q = \l + \R Z_J$
for some $Z_J \in C_\g(\l) \cap (\l)^\perp$.\par
Since $\pi$ is holomorphic, $\l^\C + \m^{01} \subset 
\q^\C + \m^{01}$ and $\q^\C + \m^{01}$ is a subalgebra
with nilradical $\m^{01}$. This 
implies that $\q = \l + \R Z_J \subset N_\g(\l^\C + \m^{01})$.\par 
By Lemma 4.10, $\dim N_\g(\l^\C + \m^{01}) \leq 
\dim \l + 1$ and therefore $\q = N_\g(\l^\C + \m^{01})$. In 
particular, the CR structure is standard 
if and only if $\q = \k$, i.e. if and 
only if $\R Z_J = \R Z_\D$. \par
If $G/L$ is a contact manifold of non-special type,
 then $\dim C_\g(\l) \cap (\l)^\perp = 1$ and hence
$\R Z_J = \R Z_\D$; in particular
any circular CR structure is standard. \par
If $G/L$ is a contact manifold of  special type,
the class of all invariant CR structures 
is explicitly classified  in \S 5. 
From that classification, the following description of all
 circular CR structures is immediately obtained.\par
\bigskip
\proclaim{Theorem 4.11} Let $M = G/L$ be a homogeneous 
contact manifold of a compact  Lie group $G$. Then  $M = G/L$
admits an invariant non-standard  circular CR structure  $(\D,J)$
if and only if  $M = SU_\ell/U_{\ell -2}$ for  $\ell \geq 2$. \par
Moreover, any invariant non-standard CR structure $(\D,J)$ on $M= SU_\ell/U_{\ell -2}$
is circular.
\endproclaim
\medskip
We refer to Theorem 5.1 for an explicit description of the non-standard circular CR structures
on  $M= SU_\ell/U_{\ell -2}$. 
\bigskip
\bigskip
\subhead 5. Classification of CR structures on special contact manifolds
\endsubhead
\bigskip
We  describe here all  
invariant CR structures $(\D_Z, J)$ on a special
contact manifold $G/L$. Recall that in this case $G$ is simple and
$L = C_G(\g(\alpha))$, by Theorem 3.6,
where either $\alpha =\mu$ is the highest 
root   or
$G = G_2$ and $\alpha = \nu$ is 
a short root. \par
We have the following orthogonal
decomposition of $\g$ 
$$\g = \l + \R Z + \m = \l + \a + \n\ ,\tag5.1$$ 
where   $\a = \g(\alpha)$ is the 3-dimensional
subalgebra associated with the root
$\alpha$, $Z = i H_\alpha \in \a$ and $\l = C_\g(\a)$ is its centralizer.\par
 Let
$(\D, J)$ be an invariant CR structure on $G/L$ which is determined
by the contact element $Z = iH_\alpha$ and 
by the decompositions 
$$\m^\C = \m^{10} + \m^{01} = \a^{10} + \n^{10} + \a^{01} + \n^{01}
\ ,\tag5.2$$
where $\a^{10} = \a^\C \cap \m^{10}$, $\n^{10} = \n^\C \cap \m^{10}$
and $\m^{01} = \a^{01} + \n^{01} = \overline
{\m^{10}}$.\par
Since $\a^\C \simeq \goth{sl}_2(\C)$ and $\a^{10} + \a^{01}$ is the orthogonal
complement to $\C Z$ in $\a^\C$,
we can write $\a^{10} = \C Z'$, for some $Z' \in \m^\C \cap \a^\C$.\par
\bigskip
Note that  a regular element $X$ of $\a^\C$ (up to rescaling) 
can  be  always identified with $i H_\alpha$, where $\alpha$ is a root
of $\g^\C$ with respect to some Cartan
subalgebra $\h$ of $\g^\C$ and such that $\a = \g(\alpha)$. In particular,
since any contact element $Z$ of $\g$
is a regular element for $\a^\C$, it can  be always  
identified with $i H_\alpha$.
\par
\medskip
If  $\alpha= \mu$ is the highest root, 
the eigenspace decomposition of $\operatorname{ad}_{H_\alpha}$ gives the
gradation   
$$\g^\C = \g_{-2} + \g_{-1} + \g_0 + \g_1 + \g_2\ ,\tag5.3$$
which is described in (3.3) and Table 1.  
Table 1 shows that for
$\g^\C \neq A_{\ell}$,   the $\g_0$-moduli $\g_{\pm1}$ are  irreducible,
their dimension is  $\dim_\C \g_{\pm1} = 1/2 \dim_\C \n^\C$ and 
$$[\g_{\pm1}, \g_{\pm1}] = \g_{\pm2}\ .\tag5.4$$
If $\g^\C = A_{\ell}$,
 each $\g_0$-module $\g_{\pm1}$ decomposes into two non-equivalent
irreducible 
$\g_0$-moduli:
 $\g_{\pm 1} = \g^{(1)}_{\pm 1} + \g^{(2)}_{\pm 1}$. Moreover,
 the following relations hold:
$$[\g^{(i)}_{1}, \g^{(i)}_{1}] = \{0\} = [\g^{(i)}_{-1}, \g^{(i)}_{-1}]\ ,
\qquad [\g^{(i)}_{1}, \g^{(j)}_{1}] = \g_2\ , \qquad
[\g^{(i)}_{-1}, \g^{(j)}_{-1}] = \g_{-2}\ ,\tag5.5$$
$$
[\g^{(i)}_1, \g_{-2}] = \g^{(j)}_{-1}\ ,
\qquad [\g^{(i)}_{-1}, \g_{2}] = \g^{(j)}_{1}\ ,
 \qquad \overline{\g^{(i)}_1}
= \g^{(i)}_{-1}
\qquad\quad (i \neq j)\ .\tag5.6$$
The moduli ${\g^{(i)}_{1}}$ and $\g^{(j)}_{-1}$ 
($i\neq j$)
 are isomorphic as $\g_0'$-moduli and, for both
values of $i$, $\dim_\C \g^{(i)}_{\pm1} = 
1/4 \dim_\C \n^\C$.\par
\medskip
When $\g^\C = G_2$
and $\alpha = \nu = \varepsilon_1$ 
is a short root, the eigenspace decomposition
of operator  $\operatorname{ad}_{H_\nu}$ defines
the following gradation of $\g^\C$:
$$\g^\C = \g_{-3} + \g_{-2} + \g_{-1} + \g_0 + \g_1 +
\g_2 + \g_3\ ,\tag5.7$$
where 
$$\g_0 =  \g'_0 + \C H_\nu\ ,
\quad \g_0' = C_{\g^\C}(\g(\nu))
= 
<E_{\pm(\varepsilon_2 - \varepsilon_3)}, 
H_{\varepsilon_2 - \varepsilon_3}>\ , $$
$$\g_2 = \C E_{\varepsilon_1} \quad, \quad
\g_{-2} = \C E_{-\varepsilon_1}\quad, \quad
\g^\C(\nu) = \g_2 + \g_{-2} + \C H_{\varepsilon_1}\ ,$$
$$\g_1 =  < E_{- \varepsilon_3},
E_{ - \varepsilon_2}>\quad,\quad
\g_3 = < E_{\varepsilon_1 - \varepsilon_3},
E_{\varepsilon_1 - \varepsilon_2}>\ ,$$
$$\g_{-i} = \overline{\g_i}\qquad \text{for} \ \ i = 1,3\tag5.8$$
(see Appendix for notation).\par
Note  that all subspaces $\g_i$   are irreducible $\g'_0$-moduli 
and that the moduli $\g_j$,
$j = \pm1, \pm3$, are  equivalent  $\g'_0$-moduli.
 Furthermore, $[\g_{\pm1}, \g_{\pm1}] =
\g_{\pm2}$ and $[\g_{\pm3}, \g_{\pm3}] =
\{0\}$. \par
\bigskip
The following Theorem gives
the complete classification 
of the invariant CR structures on  special contact manifolds.\par
\bigskip
\proclaim{Theorem 5.1} Let $(M = G/L, \D_Z)$ be a special contact 
manifold. Then:\par
\noindent{\rm a)}\  if $G\neq SU_{\ell+1}$,
 there exists (up to a  sign) a unique invariant CR structure 
$(\D_Z, J)$, and it is the  standard  one.
\par
\smallskip
\noindent{\rm b)}\  if $G = SU_2$ and hence $M = SU_2$, there
 exists a 1-1 correspondence
between the invariant CR structures (determined up 
to a sign)  and the points of the unit disc 
$$D = \{ t\in \C\ , \ |t|<1\ \} \ .
\tag5.9$$
Under  the identification  $Z = i H_\alpha$, a
point $t\in D$ corresponds to the 
CR structure $(\D, J_t)$ with the
holomorphic subspace 
$$
\m^{10} = \C(E_\alpha + t E_{-\alpha})\ .\tag5.10$$
The CR structure $(\D_Z, J_t)$ is standard if and only if 
$t = 0$.\par
\smallskip
\noindent{\rm c)}\  if $G = SU_{\ell}$, $\ell >2$,  and hence $M =
SU_{\ell}/U_{\ell-2}$, 
 the set of all 
invariant CR structures (determined up to a sign) 
 consists of:
\roster
\item"c.1)" the standard CR structure
 $(\D_Z, J^{(0)})$, induced by 
the  invariant complex structure $J^{(0)}$ on
$F_Z = SU_\ell/T^2\cdot SU_{\ell-2}$,  which is the natural
complex structure of the twistor space 
of the Wolf space $Gr_2(\C^{\ell}) = SU_\ell/S(U_2\cdot U_{\ell-2})$;
\item"c.2)" three families $(\D_Z, J_t)$, 
$(\D_Z, J'_t)$ and $(\D_Z, J^{(0)}_t)$ of invariant CR structures, 
parameterized by the points
of the unit disc $D $. 
 Under  the identification $Z = i H_\mu$, 
the  CR structures $(\D_Z,J_t)$, 
$(\D_Z,J'_t)$ and $(\D_Z, J^{(0)}_t)$
have the following holomorphic subspaces 
$$\text{(for \ }J_t\ )\quad
 \phantom{aaaaaaa}
\m^{10} = \C(E_\mu + t E_{-\mu})+ \g_1^{(1)}
+\g_{-1}^{(2)}\ ,
 \phantom{aaaaaaaaaaa}
\tag5.11$$
$$\text{(for \ }J'_t\ )\quad
 \phantom{aaaaaaa}
\m'{}^{10} = \C(E_\mu +  t E_{-\mu})+ \g_1^{(2)}
+\g_{-1}^{(1)}\ ,
 \phantom{aaaaaaaaaa}
\tag5.12$$
$$\text{(for \ }J^{(0)}_t\ )\quad 
\m''{}^{10} = \C(E_\mu + t^2 E_{-\mu}) + (\g_1^{(1)} + t \g_{-1}^{(2)})
+ (\g_1^{(2)} + t \g_{-1}^{(1)})\ ,
\tag5.13$$
 where 
$$\g^\C =  \l^\C + \C Z + \m^\C =
\g'_0  +\C(iH_{\mu}) + ( \g_{-2} + \g_{-1} + \g_1 + \g_2 )\ , $$
and where $(\g_1^{(i)} + t \g_{-1}^{(j)})$ denotes the unique $\g_0'$-invariant 
subspace of $\g_1 + \g_{-1}$, with highest weight vector $E^{(i)}_1 + t E^{(j)}_{-1}$, 
where $E^{(k)}_{\pm1}$, $k = 1,2$, are highest weight vectors of $\g^{(k)}_{\pm1}$.\par 
A CR structure
 $(\D_Z, J_t)$,  $(\D_Z, J'_t)$ or $(\D_Z, J^{(0)}_t)$ 
is standard if and only if $t = 0$.
\endroster
\endproclaim
\bigskip
\proclaim{Corollary 5.2} Let $(M = G/L, \D_Z)$ be a special 
contact manifold with $G = SU_\ell$.
\roster
\item if $M = SU_2$, then $(M, \D_Z)$ admits (up to sign)
only one standard CR structure and one family of non-standard
CR structures, parameterized by the punctured unit disc $D\setminus \{0\} \subset 
\C$;   any non-standard
CR structure is circular  and the  anti-canonical map
$\phi: M \longrightarrow \phi(M)$
is a finite covering;
\item if $M = SU_\ell/U_{\ell-2}$, $\ell >2$, then $(M, \D_Z)$ admits 
(up to a sign) exactly three standard CR structures (namely
$(\D_Z, J^{(0)})$, $(\D_Z, J_0)$ and $(\D_Z, J'_0)$) that are induced 
by three invariant complex structures of the corresponding
flag manifold $F_Z = SU_\ell/ T^2\cdot SU_{\ell-2}$, plus three families 
$(\D_Z, J^{(0)}_t)$, $(\D_Z, J_t)$ and $(\D_Z, J'_t)$ 
of non-standard CR structures, parameterized by the points of the 
punctured unit disc $t\in D\setminus \{0\}$; any non-standard CR structure
$(\D_Z, J^{(0)}_t)$ is primitive, while
the CR structures $(\D_Z, J_t)$ and $(\D_Z, J'_t)$ are circular; furthermore, 
 each  CR structure
$(\D_Z, J_t)$ or $(\D_Z, J'_t)$
 admits also 
a  CRF fibration 
$$\pi\: M = SU_\ell/U_{\ell-2} \longrightarrow Gr_2(\C^{\ell})
=SU_\ell/S(U_2\times U_{\ell-2})
$$
with  fiber $SO_3$ over the Wolf space $Gr_2(\C^{\ell})$ equipped 
with its (unique up to a sign) complex structure; finally,
for any non-standard CR structure, the  anti-canonical map
$\phi: M \longrightarrow \phi(M)$
is a finite covering.
\endroster
\endproclaim
\medskip
\remark{Remark 5.3} The  complex structures $J_0$ and $J_0'$
on $F_Z$ coincide on the fibers of the twistor  fibration 
$\pi: F_Z \to Gr_2(\C^\ell)$ but are projected into two opposite 
complex structures of $Gr_2(\C^\ell)$.
\endremark
\bigskip
\demo{Proof}
The proof of Theorem 5.1 reduces to classification of the  
decompositions (5.2), which correspond to an integrable CR structure,
for each special contact manifold $(G/L, \D_Z)$. For any decomposition
(5.2), the subspace $\a^{10}$ can be expressed as
$\a^{10} = \C Z'$ for some suitable $Z'\in \a^\C$.
Therefore we have to cases:
\roster 
\item"(1)"
 $Z'$ is a regular element of $\a^\C$;
\item"(2)" $Z'$ is a non-regular (hence nilpotent) element of
$\a^\C \simeq \goth{sl}_2(\C)$.
\endroster
\bigskip
\noindent{\it Case (1)}:\par
First of all we  show that this case may 
occur only if $\g = \goth{su}_\ell$.\par
Consider first    $\a = \g(\mu)$, with $\mu$ long root of 
the simple group $G$. Since $Z'$
is regular, we may 
assume that  $Z' = i H_\mu$ and we may consider the corresponding
graded decomposition 
(5.3). Recall that
$\l^\C = C_\g(\g(\mu)) = \g'_0$. \par
 Hence 
the subalgebra $\frak b =\l^\C + \m^{10}$  is contained in
$$\l^\C + \m^{10} = \g'_0 + \a^{10} + \n^{10} = 
\g'_0 + \C H_\mu + \n^{10} = \g_0 + \n^{10} \subset \g_0 + \g_1 + \g_{-1}$$
since  $\n^\C \subset \g_1 + \g_{-1}$, being
orthogonal to $\a^\C = \C H_\mu + \g_2 + \g_{-2}$.
In case
$\g^\C \neq A_{\ell}$, $\g_1$ and $\g_{-1}$ are irreducible
$\g_0$-modules  and hence either $\g_1$ or 
$\g_{-1}$ is included in $\n^{10}$. However $[\g_1, \g_1] = \g_2$
and $[\g_{-1}, \g_{-1}] = \g_{-2}$, and hence there is no subalgebra
$\frak b$ of $\g_0 + \g_1 + \g_{-1}$ which contains $\g_0$ properly.
Hence $\g^\C = A_{\ell}$.
In this  case,  each $\g_{\pm1}$ 
decomposes into two not equivalent
irreducible $\g_0$-moduli $\g^{(i)}_{\pm1}$, $i = 1,2$,
of dimension equal to 
$1/4 \dim_\C \n^\C$,
which verify (5.5) and (5.6). Like before it is easy to check
that the $\g_0'$-moduli decomposition of $\n^{10}$ has the form
$\n^{10} = \g_1^{(i)} + \g_{-1}^{(j)}$
for some choice of $i$ and $j$.
If $i=j = 1$, then
  $\n^{01} = \overline{\n^{10}} = \g_1^{(1)} + \g_{-1}^{(1)} = \n^{10}$
and this contradicts the condition $\n^{10} \cap \n^{01} = \{0\}$.
A similar contradiction arises when $i = j= 2$.\par
 In conclusion, if $\alpha = \mu$ is a long root, then
 $\g^\C = A_{\ell}$ and for any fixed $\a^{10}$ 
there exist at most two 
CR structures, i.e. those 
corresponding to  the
following two  possibilities for $\n^{10}$:
$$\n^{10} = \g_1^{(1)} + \g_{-1}^{(2)}\ ,\ \qquad 
\n^{10} = \g_1^{(2)} + \g_{-1}^{(1)}\ .\tag5.14$$
\medskip
It remains to consider the case in which $G= G_2$ and
 $\a = \g(\nu)$, with $\nu$  short root  of $\g^\C$. 
 We  
assume that  $Z' = i H_\nu$ and we  consider the  corresponding  
graded decomposition 
(5.7).\par
Then $\l^\C + \m^{10}$ is contained in
$$\l^\C + \m^{10} = \l^\C + \a^{10} +  \n^{10} 
= \g'_0 + \C H_\nu + \n^{10} \subset \g_0 + \g_1 + \g_{-1}
+ \g_{3} + \g_{-3} $$
because $\n^\C$ is orthogonal to $\a^\C = \C H_\nu + \g_2 + \g_{-2}$.
Since $\l^\C + \m^{10} = \g^0 + \n^{10}$ 
is a subalgebra and 
$\dim_\C \g_{\pm1} = 
\dim_\C \g_{\pm3}  = \frac{1}{2} \dim_\C \n^{10}$, 
$\n^{10}$  contains two of the four
irreducible $\g_0$-moduli $\g_{\pm1}$ and $\g_{\pm3}$. The only possibility
for $\n^{10}$, so that 
$\g_0 + \n^{10}$ is a subalgebra, is 
$\n^{10} = \g_{-3} + \g_{3}$.
This  implies that $\n^{01} = \overline{\n^{10}} = \g_{-3} + \g_{3} =
\n^{10}$
and it contradicts the condition 
$\m^{10} \cap \overline{\m^{10}} = \{0\}$.\par
\medskip
Now it remains to classify
the  invariant CR structures on
 $(SU_{\ell}/U_{\ell-2}, \D_Z)$.\par
For the following part of the proof, it is more convenient
to identify 
the contact element $Z$ (and no longer $Z'$)  with 
 $iH_\mu$. We also consider
the decomposition (5.3) determined by $Z = iH_\mu$. \par
Since $Z'$ is a regular element which is orthogonal to $Z = i H_\mu$,
it is (up to a factor) of the form 
$$Z' = E_\mu +  t E_{-\mu}\ ,\qquad |t| \neq 0\ .\tag5.15$$
Exchanging $\a^{10}$ with $\a^{01} = \overline{\a^{10}}$ if necessary
 (which
corresponds to changing sign to the complex structure),  we may  assume that
$0< |t| \leq 1$.\par
Since $\a^{10}\cap \a^{01}= \{0\}$ and hence $  E_\mu +  t E_{-\mu}$
and $\bar t E_\mu +  E_{-\mu}$ are linearly independent,
$t$   verifies  the condition
$$\det\left[\matrix 1 & t\\ \bar t & 1\endmatrix
\right] = 1 -|t|^2 \neq 0\tag5.16$$
and therefore $t \in D \setminus \{0\} = 
\{ 0< |t| < 1 \}$.\par
We claim that for any point $t\in D\setminus \{0\}$ 
there exist exactly three invariant CR structures, 
whose associated subspace 
$\a^{10}$ is equal to $\C(E_\mu + t E_{-\mu})$. In fact, one can check
that the only $\g'_0$-invariant subspaces $\m^{10}$ of $\C(E_\mu + t E_{- \mu}) + 
\g_1 + \g_{-1}$, which verify (i) and (ii) of Definition 4.1,
are
either  (5.11), (5.12) or a subspace of the form
$$\m^{10} = \C(E_\mu + t E_{- \mu}) +  (\g_1^{(1)} + s \g_{-1}^{(2)})
+ (\g_1^{(2)} + s \g_{-1}^{(1)})\tag 5.17$$
for some coefficient $s$. One can also check that  the subspaces 
(5.11) and (5.12) verify also the integrability condition, while (5.17) 
satisfies 
the integrability condition if and only if $t = s^2$. This proves 
that (5.11), (5.12) and (5.13) are the only holomorphic  subspaces of $\m^\C$
containing  $\C(E_\mu + t E_{- \mu}) $. In particular, they
define
three  distinct invariant CR structures, which we 
denote by $(\D_Z, J_t)$, $(\D_Z, J'_t)$ and $(\D_Z, J^{(0)}_t)$. \par
If $\ell = 2$ and hence
$M = SU_2$, then $\n = \{0\}$ and 
 the three CR structures  $(\D_Z, J_t)$, 
$(\D_Z, J'_t)$ and $(\D_Z, J^{(0)}_t)$ coincide for any $t$.\par
\medskip
Since for any $t \neq 0$ the holomorphic subspaces $\m^{10}$
and $\m'{}^{01}$ are  not $\operatorname{ad}_Z$-invariant, 
 any CR structure 
$(\D_Z, J_t)$, 
$(\D_Z, J'_t)$ or $(\D_Z, J^{(0)}_t)$ ($t \neq 0$) is non-standard
by  Lemma 4.5.\par
\bigskip
\noindent{\it Case (2)\/}:\par
Since $Z'$ is not regular, it
is a nilpotent element of $\a^\C = \goth{sl}_2(\C) = \g^\C(\alpha)$.
Then we may always choose
a Cartan subalgebra $\C H_\alpha$ of $\a$ so that $Z'\in \C E_\alpha$.
Furthermore,
since the contact element $Z$ is orthogonal to 
$\a^{10} + \a^{01} = \C E_\alpha + \overline{\C E_\alpha} = 
\C E_\alpha +  \C E_{-\alpha}$, we may 
assume (after rescaling) that $Z = iH_\alpha$.\par
Consider first that   $\alpha =\mu$ is a long root of $G$
and take 
the gradation (5.3) of $\g^\C$ determined with $H_\mu$. Then $\g_2
= \C Z' = \a^{10}$ and hence
$$\l^\C + \m^{10} = \g'_0 + \g_2 + \n^{10} \subset \g'_0 + 
\g_2 + \g_1 + \g_{-1}\ .$$
Assume that $\g^\C \neq A_{\ell}$. Then the $\g'_0$-moduli $\g_{\pm1}$
are
irreducible and $[\g_{\pm 1}, \g_{\pm 1}] = \g_{\pm 2}$. Hence
the only  subalgebra of $\g'_0 + 
\g_2 + \g_1 + \g_{-1}$, which properly contains 
$\g'_0 + 
\g_2$, is $\g'_0 + \g_1 + \g_2$. Hence $\m^{10}= \g_1 + \g_2$.\par
Vice versa, $\m^{10} = \g_1 + \g_2$ is a holomorphic subspace
of $\m^\C = (\l^\C + \C Z)^\perp = \g_0^\perp$  and
hence it corresponds to an
invariant CR structure on $(G/L, \D_Z)$. Since $Z = i H_\mu\in 
N_{\g}(\g'_0 + \g_{-1} + \g_{-2}) = N_\g(\l^\C + \m^{01})$,
this CR structure is standard.\par
\medskip
Assume now that $\g^\C =A_{\ell}$ and again
consider the decomposition (5.3) determined by $Z = iH_\mu$. 
Since
$\dim_\C \g^{(i)}_{\pm1} = 1/4 \dim_\C \n^\C$, the $\g'_0$-module
$\n^{10}$ can be written in one of the following five forms:
$$1)\ \n^{10} = (\g_1^{(1)})_\varphi + 
(\g_{-1}^{(1)})_\psi\ ,\qquad
2)\ \n^{10} = \g^{(1)}_1 + \g^{(2)}_{-1}\ ,\qquad
3)\ \n^{10} = \g^{(2)}_1 + \g^{(1)}_{-1}\ ,$$
$$4)\ \n^{10} = \g_1\ ,\qquad 5)\ \n^{10} = \g_{-1}\ ,$$
where $\varphi: \g^{(1)}_1 \to \g^{(2)}_{-1}$ and 
 $\psi: \g^{(1)}_{-1} \to \g^{(2)}_{1}$ are two $\g'_0$-equivariant
homomorphisms and 
where $(\g_1^{(1)})_\varphi$ and   
$(\g_{-1}^{(1)})_\psi$  denote the subspaces of the form
$$(\g_1^{(1)})_\varphi = \{ X + \varphi(X)\ : X\in \g^{(1)}_1\}\ , \qquad 
(\g_{-1}^{(1)})_\psi =  \{ X + \psi(X)\ : X\in \g^{(1)}_{-1}\}\ .$$
Case 5) cannot occur because in that case $[\n^{10}, \n^{10}] = \g_{-2}$
and this contradicts the fact that $\g'_0 + \n^{10} + \g_2$ is a subalgebra.\par
Also 
case 1) may not occur. In fact, $\varphi$ is either trivial
or an isomorphism. In case $\varphi$ is an isomorphism,  
for any $0\neq X \in \g^{(1)}_{1}$, it is possible 
to find an element $Y\in 
\g^{(1)}_{-1}$ so that $[\varphi(X), Y]$ is non-trivial and belongs
to $\g_{-2}$. Hence, 
$$ [X + \varphi(X), Y + \psi(Y)] \underset{\mod \g_0 + \g_1 + \g_2}\to\equiv
[\varphi(X), Y] \in \g_{-2} \ .$$
This contradicts the fact that
 $\l^\C + \m^{10}$ is a subalgebra
of $\g'_0 + \g_1 + \g_2$. We conclude that, if case 1) occurred, 
$\n^{10} = \g_1^{(1)} + (\g_{-1}^{(1)})_\psi$.
Now,  for any $X\in \g_1^{(1)}$
we may consider an element $Y + \psi(Y)\in   (\g_{-1}^{(1)})_\psi$ so that
$[X, Y] = \lambda H_\mu$ for some $\lambda \neq 0$. Hence 
$$[X, Y + \psi(Y)]  = \lambda H_\mu \mod \g'_0 + \g_2$$
This gives a contradiction with the fact that $\g'_0 + 
\n^{10} + \g_2$ is a subalgebra and the claim is proved. \par
For the cases 2), 3) and 4), 
$\m^{10}$ equals one of the following three subspaces
$$\g^{(1)}_1 + \g^{(2)}_{-1} + \g_2\ ,
\g^{(2)}_1 + \g^{(1)}_{-1} + \g_2\ ,
\g_1 + \g_2
\tag5.18$$
and one can check that any of them  is a holomorphic subspace.\par
By Proposition 4.2, they determine three  distinct CR structures denoted by
 $(\D, J)$, $(\D, J')$ 
and $(\D, J^{(0)})$, respectively. For any 
of the three subspaces (5.18), the normalizer
$N_\g(\l^\C + \m^{01})$ contains $\g_0 \cap \g
= \l + \R Z$ and hence the corresponding    CR structures are standard.  \par
Finally, observe that $(\D, J^{(0)})$
is induced by the invariant complex structure
$J_F$ on the flag manifold $F_Z = SU_\ell/ T^2\cdot SU_{\ell-2}$
which is associated to the following
black-white Dynkin graph
$$\dynkin\root{}\link\wroot{}\link\wroot{}\link\dots
\link\wroot{}\link\wroot{}
\link\root{}
\enddynkin$$
and which is the invariant complex structure of the twistor space of the Wolf
space
$Gr_2(\C^{\ell}) = SU_\ell/S(U_2\cdot U_{\ell-2})$; moreover, the subspace of $J^{(0)}$ 
coincides with the subspace given in (5.13) for $t = 0$; on the other hand,
the subspaces of $J$ and $J'$ are the subspaces given in 
(5.11) and (5.12) for $t = 0$.
All corresponding CR structures  coincide if $M = SU_2$.\par 
\medskip
It remains to consider the case in which $G= G_2$ and
 $\a =  \g(\nu)$, where $\nu$ 
is a  short root. Consider the decomposition (5.7)
determined by $H_\nu$ so that 
 $\C Z' = \C E_\nu = \g_2$.\par
 As before, we identify
 $Z$   with $iH_\nu$.  We have 
$$\l^\C + \m^{10} = \g'_0 + \a^{10} + \n^{10} \subset 
\g'_0 + \g_2 + \g_{-1} + \g_{1}
+ \g_{-3} + \g_{3} $$
because $\n^\C$ is orthogonal to $\a^\C = \C H_\nu + \g_{-2} + \g_{2}$.
We claim that $\g_3 \subset \n^{10}$. In fact, for any 
element $X\in \n^{10}$ consider the decomposition
$$X = X_{-3} + X_{-1} + X_1 + X_3\ ,\qquad  X_i \in \g_i\ .$$
Then,  one of the four 
vectors $X$, $X' =[E_\nu, X]$, 
$X'' = [E_\nu,[E_\nu, X]]$, $X''' = [E_\nu,[E_\nu,[E_\nu, X]]]$ is a 
non-trivial element of $\g_3$ and it belongs
to $\n^{10}$. Since $\g_3$ is $\g'_0$-irreducible,
the claim follows.\par
Similarly, we  claim that $\g_1 \subset \n^{10}$. To prove this, take any 
element $X\in \n^{10}$ which has a decomposition of the form
$$X = X_{-3} + X_{-1} + X_1\ ,\qquad  X_i \in \g_i\ .$$
Then, either $X$ or $X' =[E_\nu, X]$ or $X'' = [E_\nu,[E_\nu, X]]$ is a 
non-trivial element of $\g_1 + \g_3$, with non-vanishing projection
on $\g_1$. This implies that $\g_1 \cap \n^{10} \neq \{0\}$ and hence
that $\g_1 \subset \n^{10}$. Since $\dim_\C(\g_1 + \g_3) = 
\dim_\C \n^{10}$,
we conclude that
 $\n^{10} = \g_1 + \g_3$ and 
that 
$\m^{10} = \g_1 + \g_2 + \g_3$. Indeed, since $\l^\C + \g_1 + \g_2 + \g_3$
is always a subalgebra, there exists 
an integrable
CR structure whose associated holomorphic subspace
is $\m^{10} = \g_1 + \g_2 + \g_3$.
Furthermore,    
$N_\g(\l^\C + \m^{01})$ contains $Z = iH_\nu$
and hence  this CR structure is standard.
\qed 
\enddemo
\bigskip
\demo{Proof of Corollary 5.2} (1) By Theorem 5.1, 
it remains  only  need to check  that 
any non-standard CR structure on $M = SU_2$ is circular and 
that the associated anti-canonical map is a finite covering.\par 
By (5.10), the CR structure $(\D_Z, J)$ is 
non-standard if and only if the corresponding
 holomorphic subspace is of the form $\m^{10} = \C( E_{\alpha} + t E_{-\alpha})$ 
with 
$0 < |t| < 1$. Since $\l = \{0\}$
and the element $E_{\alpha} + t E_{-\alpha}$ is a regular
element of $\goth{sl}_2(\C)$, then $\m^{10}$ is a Cartan subalgebra of $\g^\C = 
\goth{sl}_2(\C)$ and any parabolic subalgebra $\p$ which contains $\m^{10}$
verifies the conditions a), b) and c) of Lemma 4.8. This implies that 
$M = SU_2$ admits a CRF fibration over $SU_2/T^1$, 
where $T^1$ is the 
1-dimensional subgroup generated by 
the subspace $\goth t = \p \cap \goth{su}_2$.\par
On the other hand, when $0 < |t| < 1$,
$$\eqalign{N_\g(\C(E_{\alpha} + t E_{-\alpha})) &= \cr
= 
\{\ X = a(i H_\alpha) + & b( E_{\alpha} + E_{-\alpha}) + 
i c( E_{\alpha} - E_{-\alpha})\in \goth{su}_2 \ : 
[X, E_{\alpha} + t E_{-\alpha}] \in \cr
& \phantom{aaaaaaaaaaaaaaaaaaaaaaaa}
\ \in \C (E_{\alpha} + t E_{-\alpha})\} = \{0\}\ .\cr}
$$
Then, by the remarks before Theorem 4.9,   the 
stabilizer $Q$ of the image of the anti-canonical map 
 $\phi(SU_2) = SU_2/Q$ 
is 0-dimensional  and 
the anti-canonical map is a covering map.\par
\medskip
(2) We first observe that each non-standard CR structure 
$(\D_Z, J^{(0)}_t)$ is primitive.  In fact, by Lemma 4.8, if one of such CR structures
is   non-primitive, then there exists a 
parabolic subalgebra $\p \subsetneq \g$, which verifies a), b) and c) 
of Lemma 4.8. On the other hand, one can check
that in this case, 
there is no proper subalgebra of $\g^\C$ which properly contains $\l^\C + 
\m^{10}$, with $\l^\C = \g'_0$ and $\m^{10}$ as in (5.13). \par
Now, we want to
prove that each 
 non-standard CR structure $(\D_Z, J_t)$ or $(\D_Z, J'_t)$ 
admits a CRF fibration onto $Gr_2(\C^{\ell}) = SU_\ell/S(U_2\cdot
U_{\ell-2})$.\par
Indeed, note  that, if we consider the decomposition 
(5.3) determined by the regular contact element $Z = i H_\mu$,  any 
CR structure $(\D_Z, J_{t})$ or $(\D_Z, J'_{t})$ 
corresponding to the holomorphic subspaces 
defined in (5.11) and (5.12)
verifies
$$\l^\C + \m^{01} \subset \p = \g_0 + \g^{(1)}_{-1} + 
\g^{(2)}_{1} +  \g_{-2} + \g_2 \ ,\tag5.19$$
$$\l^\C + \m'{}^{10} \subset \p' = \g_0 + \g^{(2)}_{-1} + 
\g^{(1)}_{1} + \g_{-2} + \g_2 \ ,\tag5.20$$
respectively. A reductive part for  both subalgebras
$\p$ and $\p'$
  is
$\r = \r' = (\l + \a)^\C$.
Therefore, by Lemma 4.8, the CR structures $(\D,J_{t})$ and 
 $(\D,J'_{t})$ are non-primitive and they admit a CRF fibration over 
the Wolf space $SU_{\ell+1}/ S(U_2\cdot U_{\ell-1})$ with typical fiber
$S(U_2\cdot U_{\ell-1})/U_{\ell-1} = SO_3$. \par
We now want to prove that any non-standard CR structure $(\D_Z, J_t)$ or $(\D_Z, J'_t)$ 
admits also a CRF fibration with standard fiber $S^1$. Let us use  
 the same notation
as before and observe that, for any complex holomorphic subspace 
$\m^{10}$ or $\m'{}^{10}$
defined in (5.11) or (5.12), the element 
$X = E_\mu + t E_{-\mu} \in \m^{10} \cap \m'{}^{10}$ is a regular element of 
$\g^\C(\mu)  \subset \g^\C$. 
Hence, if we denote by $\hat \p$ 
any  parabolic subalgebra $\hat \p(\mu) \subset \g^\C(\mu)$, which properly contains 
$E_\mu + t E_{-\mu}$ or $E_\mu + t^2 E_{-\mu}$, we get that
$$\l^\C + \m^{01} \subset  \p_\mu = \g'_0 + \g^{(1)}_{-1} + 
\g^{(2)}_{1} +  \hat \p(\mu)\ ,\tag 5.21$$
$$\l^\C + \m'{}^{10} \subset  \p'_\mu = \g_0 + \g^{(2)}_{-1} + 
\g^{(1)}_{1} + \hat \p(\mu) \ .\tag5.22$$
Note that $\p_\mu$ and $\p'_\mu$ are two parabolic subalgebras of 
$\g^\C$ which verify a), b) and c) of Lemma 4.8 and hence 
that the CR structures
$(\D, J_t)$ and $(\D, J'_t)$ admit CRF fibrations with 
1-dimensional fibers.\par
It remains to check that the anti-canonical map of any non-standard 
CR structure is a covering map. As in the proof of (1), this reduces to
checking  that for any 
holomorphic subspace defined in  (5.11) and (5.12),
$N_\g( \l^\C + \m^{10})=  N_\g( \l^\C + \m'{}^{10})  = \l^\C$  and hence that 
the image of the anti-canonical map has the same dimension as $G/L$. 
\qed
\enddemo
\bigskip
\bigskip
\subhead 6. Classification of non-standard CR structures
\endsubhead
\bigskip
\subsubhead 6.1 Notation
\endsubsubhead
\medskip
In all this section, 
\roster
\item"-" $(G/L, \D_Z)$ denotes a
simply connected non-special homogeneous contact  manifold
of a compact Lie group $G$; 
\item"-" $\k = C_\g(Z)
= \l + \R Z$ is the orthogonal decomposition of the centralizer $\k$ of
$Z$ and $\m$ is the orthogonal complement to $\k$ in $\g$;
\item"-" $\h \subset \k$ is a Cartan subalgebra of $\k$ and hence of $\g$;
\item"-" $\theta = \B\circ Z|_\h$ is the 1-form on $\h$ dual to $Z$ and 
$\vartheta = - i\theta = -i\B\circ Z|_\h$; we will refer to both 
of them as {\it contact forms\/};
\item"-" $R$ (resp. $R_o$) is the root system of $(\g^\C, \h^\C)$
(resp. of $(\k^\C, \h^\C)$) and $R' = R\setminus R_o$;
\item"-" $E_\alpha$ is the root vector with root $\alpha$ 
in the Chevalley normalization (see e.g. [7]); 
\item"-" a subset $S\subset R$ is called {\it closed subsystem\/} if 
$(S+S)\cap R \subset S$;
\item"-" if $S$ is a closed subsystem of roots, then $\g(S) \subset \g^\C$
is the subalgebra generated by the root vectors $E_\alpha$, $\alpha\in S$; 
\item"-" recall that the root vectors
$E_\alpha$, $\alpha \in R'$, span $\m^\C$;
\item"-" $\m(\alpha)$ denotes the irreducible $\k^\C$-submodules
of $\m^\C$, with  highest weight $\alpha\in R'$;
\item"-" if $\m(\alpha)$ and $\m(\beta)$ are equivalent
as $\l^\C$-moduli, we denote by $\m(\alpha) + t\m(\beta)$  the irreducible 
$\l^\C$-module with the highest weight vector $E_\alpha + 
t E_{\beta}$, $\alpha, \beta \in R'$, $t\in \C$; note that together with 
$\m(\beta)$,
these moduli exhaust all the irreducible $\l^\C$-submoduli of $\m(\alpha) 
+ \m(\beta)$ (see Lemma 6.1);
\item"-" by {\it Dynkin graph $\Gamma$\/}
 we will understand the Dynkin graph associated with a
root system $R$ of a compact semisimple Lie algebra $\g$; we associate
with the nodes of $\Gamma$ the simple roots of $R$ as in [7]
(see Table 4 in the Appendix).
\endroster
\bigskip
\subsubhead 6.2 Preliminaries
\endsubsubhead
\medskip
By the results in \S 5, 
the classification of  invariant CR structures
  reduces 
to the  classification of  non-standard CR structures 
on homogeneous contact manifolds of non-special type. This will be 
the contents of \S 6.3 and \S 6.4. \par
\medskip
In this section we give two important lemmata
that settle the main tools for the classification.
The first Lemma is an immediate corollary of
Proposition 3.9.\par
\bigskip
\proclaim{Lemma 6.1} Let $(M = G/L, \D_Z, J)$ be a homogeneous 
CR manifold associated with holomorphic subspace $\m^{10} \subset \m^\C$
and  $J$ the associated complex structure on $\m$.
Assume also that $G\neq G_2$ or that 
$G=G_2$  and that the contact form 
$\vartheta$ is not proportional to a short root of $R$.\par
Then a minimal $J$-invariant $\k^\C$-submodule $\n$ of $\m^\C$
 is either $\k^\C$-irreducible (and hence $\n = \m(\alpha)$
for some $\alpha \in R'$) or it is the sum
$\m(\alpha) + \m(\beta)$ of two such $\k^\C$-modules, where
the roots $\alpha$ and $\beta$ are $\vartheta$-congruent (i.e.
$\beta = \alpha + \lambda \vartheta$, for some $\lambda\in \R$).
\endproclaim
\demo{Proof}
Consider the decomposition $\m^\C = \sum \m(\gamma)$
into irreducible $\k$-submodules as in \S 3.3. 
The claim follows immediately from the fact that  any 
 $\operatorname{ad}_\l$-invariant complex structure $J$ on $\m$
preserves the $\l^\C$-{\it isotypic components\/} (i.e. the sum
of all mutually equivalent irreducible $\l^\C$-modules) and
that, under the hypotheses of Proposition 3.9, the multiplicity 
of any irreducible $\l$-module $\m(\gamma)$ is less or equal to 
2.\qed
\enddemo
\proclaim{Lemma 6.2} Let $(G/L,\D_Z, J)$ be a homogeneous CR manifold
with non-standard CR structure.
Then $G$ is either simple or of the form $G = G_1 \times G_2$,
where each $G_i$ is simple. \par
Moreover, if $G = G_1 \times G_2$ and $R = R_1 \cup R_2$ 
is the corresponding decomposition of the root system, then
there exist two roots $\mu_1 \in R_1$, $\mu_2 \in R_2$, such that 
the pairs of roots $(\mu_1, -\mu_2)$ and $(-\mu_1, \mu_2)$ are the only ones 
which are $\vartheta$-congruent; in particular, $\vartheta = \mu_1 + \mu_2$
is proportional to no root.
\endproclaim 
\demo{Proof} Since the CR structure $(\D_Z, J)$ is non-standard, 
the associated complex structure $J$
on $\m$ is not $\operatorname{ad}_\k$-invariant;
in particular there exists some minimal $J$-invariant $\k^\C$-module
in $\m^\C$, which is not $\k^\C$-irreducible.
By Lemma 6.1, there exist at least two roots
$\alpha, \beta$, which are $\vartheta$-congruent. Without
loss of generality, we may assume that $\vartheta = \alpha-\beta$.\par
 If $\vartheta$ is  proportional to some
root $\gamma$, then this root belongs to some summand
$\g_i$ of $\g$, $i = 1, \dots, r$. 
Hence, $\k = C_\g(Z)$ contains all other simple summands
of $\g$ and the same holds for $\l$. By effectivity, this implies
that $\g = \g_1$. \par
If $\vartheta = \alpha- \beta$ is not proportional to any
root and $\alpha$ and $\beta$ belong to the same summand
$\g_1$, then $\g = \g_1$ as before. Assume that they belong to two 
different summands $\g_1$ and $\g_2$. The same arguments of before 
show that $\g = \g_1 \oplus \g_2$
and that 
$\pm(\alpha, \beta)$ are the only pairs of roots which are
$\vartheta$-congruent.
\qed
\enddemo
\medskip
 We will perform the classification by considering separately 
two cases: when the contact form $\vartheta$ 
is proportional  to a root and when it is  not proportional  to any root.
Note that by Lemma 6.2, the first case may occur only when $G$ is simple.\par
\bigskip
\subsubhead 6.3 Case when the contact form is proportional to a root
\endsubsubhead
\medskip
Recall that the Weyl group of a simple Lie group
acts transitively on the set of roots
of the same length. In particular any long root can be considered as a maximal
root. Since we assume that the contact manifold $(M = G/L, \D_Z)$ is 
non-special and $G$ is simple, we may suppose that $\vartheta$ is proportional 
to a short root (i.e. strictly shorter
then a long root) and hence $G$ equals either $SO_{2n+1}$, $Sp_n$ or $F_4$.
Note that if $G = G_2$ then any contact manifold
$(G_2/L, \D_Z)$, with contact form $\vartheta$ proportional to a
short root, is special (see \S 3.2.2).\par
\medskip  
\proclaim{Proposition 6.3} Let
 $(G/L, \D_Z)$ be a homogeneous non-special contact  manifold
of a simple group $G$, such that the contact form
$\vartheta$ is proportional to a root. Then:\par
\roster 
\item  $G/L$ is    $SO_{2n+1}/SO_{2n-1}$,
$Sp_n/Sp_1\times Sp_{n-2}$ or $F_4/SO_7$ and $\vartheta$ is proportional to 
a short root of $G$;
\item there exists a 1-1 correspondence between the invariant 
CR structures on $(G/L, \D_Z)$ (determined up 
to a sign) and the points of the unit disc $D \subset 
\C$; 
\item more precisely,  any point $t \in D$ corresponds  to the CR structure 
 $(\D_Z, J_t)$ whose holomorphic subspace $\m^{10}$ is listed in the
following table (see \S 6.1 for notation):\par
\medskip
\moveright 0.2 cm
\vbox{\offinterlineskip
\halign {\strut\vrule\hfil\  $#$\ \hfil &\vrule\hfil\  $#$\ 
\hfil&\vrule\hfil\ $#$\ 
\hfil\vrule\cr
\noalign{\hrule}
{\phantom {\sum_1^N}}G/L
{}_{\phantom {\sum_1^N}}^{\phantom {\sum_1^N}} & 
 {\phantom {\sum_1^N}}\vartheta{}_{\phantom {\sum_1^N}}^{\phantom {\sum_1^N}} &
 \m^{10} 
 \cr \noalign{\hrule}
\underset{\phantom{A}}\to{
\overset{\phantom{A}}\to{
\frac{SO_{2n+1}}{SO_{2n-1}} = S(S^{2n})
}}
 & \varepsilon_1 &
\smallmatrix
\m(\varepsilon_1+\varepsilon_2) + t
 \m(- \varepsilon_1 +\varepsilon_2)  
\endsmallmatrix
\cr \noalign{\hrule}  
\underset{\phantom{A}}\to{
\overset{\phantom{A}}\to{
\frac{Sp_n}{Sp_1\times Sp_{n-2}} = 
S(\Bbb H P^{n-1}) 
}}
&  
\varepsilon_1 +\varepsilon_2 &
\smallmatrix
(\m(2\varepsilon_1) + t^2
 \m(-2\varepsilon_2))\oplus
(\m(\varepsilon_1+ \varepsilon_3) + t
 \m(-\varepsilon_2 +\varepsilon_3))
\endsmallmatrix
\cr \noalign{\hrule}
\underset{\phantom{A}}\to{
\overset{\phantom{A}}\to{
\frac{F_4}{Spin_7} = S(\Bbb O P^2)
}}
 &    \varepsilon_1 &
\smallmatrix
(\m(\varepsilon_1+\varepsilon_2) + t^2
 \m(-\varepsilon_1+\varepsilon_2))\oplus \\
(\m(1/2(\varepsilon_1 +\varepsilon_2+ 
\varepsilon_3+\varepsilon_4)) +
t
 \m(1/2(-\varepsilon_1+\varepsilon_2 + 
\varepsilon_3 +\varepsilon_4)))
\endsmallmatrix 
\cr \noalign{\hrule} 
}}
\medskip
\noindent
\item a CR structure $(\D_Z, J_t)$ is standard if and only if 
$t=0$; in all other cases it  
is primitive.
\endroster 
\endproclaim
\demo{Proof} For each group $G$ equal to $SO_{2\ell+1}$, $Sp_\ell$ or $F_4$
we may assume that $\vartheta$ is the short root 
$\vartheta = \varepsilon_1$, $\varepsilon_1 + \varepsilon_2$
 or $\varepsilon_1$, respectively. The associated decomposition
$\g = \l + \R Z + \m$ is given in
Table 2 of the Appendix. 
It is not difficult to determine the decomposition of $\m^\C$
into irreducible  submoduli. The result is given in 
 Table 2.
Then one has 
to find all  decompositions  $\m^\C = \m^{10} + \m^{01}$ into two
$\l^\C$-modules which satisfy the following  conditions:
 a) $\m^{01} = \overline{\m^{10}}$; b)
$[\m^{10}, \m^{10}] \subset \m^{10} + \l^\C$.
The moduli $\m^{10}$ which satisfy condition a) are of the following form:
$$\alignat 2
G &= SO_{2\ell + 1}: \ \ &
\m^{10}  &= \m_t^{10} = \m(\varepsilon_1 + \varepsilon_2) + 
t \m(-\varepsilon_1 + \varepsilon_2)\ ; \\
G &= Sp_{\ell}: 
&\m^{10} &= \m_{t,s}^{10} = 
(\m(2\varepsilon_1) + s
 \m(-2\varepsilon_2))\oplus
(\m(\varepsilon_1+ \varepsilon_3) + t
 \m(-\varepsilon_2 + \varepsilon_3))\ ;\\
G &= F_4:
&\m^{10}  &= \m_{t,s}^{10}{}' =
\endalignat$$
$$= (\m(\varepsilon_1+\varepsilon_2) + s
 \m(-\varepsilon_1+\varepsilon_2))\oplus 
(\m(1/2(\varepsilon_1 +\varepsilon_2+ 
\varepsilon_3+\varepsilon_4)) +
t
 \m(1/2(-\varepsilon_1+\varepsilon_2 + 
\varepsilon_3 +\varepsilon_4))) $$
for some $s, t \neq 0$. One can 
easily check that $\m^{10}_t$ verifies condition b) for every $t$. The module
$\m^{10}_{t,s}$ verifies condition b) if and only if $s=t^2$. To prove it one 
should observe that the only  brackets
between $\l^\C$-weight vectors in $[\m^{10}_{t,s}, \m^{10}_{t,s}]$, which are
non-trivial  modulo $\l^\C$, are 
$$\alignat 2
[E_{\varepsilon_1 + \varepsilon_i} + t E_{-\varepsilon_2 + 
\varepsilon_i}, & E_{\varepsilon_1 -\varepsilon_i} + t E_{-\varepsilon_2 -
\varepsilon_i}] \equiv \ &\\
\ &\equiv 
N_{\varepsilon_1 + \varepsilon_i, \varepsilon_1 -\varepsilon_i}
E_{2\varepsilon_1} + t^2
N_{-\varepsilon_2 + \varepsilon_i, -\varepsilon_2 +\varepsilon_i}
E_{-2\varepsilon_2}\ & \mod \l^\C\\
[E_{\varepsilon_1 + \varepsilon_i} + t E_{-\varepsilon_2 + 
\varepsilon_i}, & E_{\varepsilon_2 -\varepsilon_i} + t E_{-\varepsilon_1 -
\varepsilon_i}] \equiv\ &  \\
\ &\equiv
N_{\varepsilon_1 + \varepsilon_i, \varepsilon_2 -\varepsilon_i}
E_{\varepsilon_1 + \varepsilon_2} + t^2
N_{-\varepsilon_2 + \varepsilon_i, -\varepsilon_1 +\varepsilon_i}
E_{-\varepsilon_1 - \varepsilon_2}\ &\mod \l^\C
\endalignat$$
By a straightforward computation, it follows that 
these vectors are in $\m^{10}_{t,s}$ if and only if 
$s = t^2$.\par
A similar argument shows that also $\m^{10}_{t,s}{}'$ verifies condition b)
if and only if $s = t^2$. \par
Observe that
up to an exchange between $\m^{10}$ and $\m^{01}$ (which corresponds
to changing the sign of complex structure $J$), we may always assume that 
$|t| \leq 1$. 
It remains to check 
the condition 
$\m^{01}\cap \m^{10} = \{0\}$: in all cases, this implies  
$\det\left[\smallmatrix 1 & t \\ \bar t &
1\endsmallmatrix \right] \neq 0$ and hence that 
$|t| <1$. \par
\smallskip
To prove (4), note that, in all cases
listed in the table above, $N_\g(\l^\C + \m^{01})$ contains $Z$
only if $t = 0$ and hence, by Theorems 4.9 and 4.11, this is the only case
when the CR structure is standard. Moreover, in all  cases, if $t\neq 0$ 
there exists no proper
parabolic subalgebra $\p \supset \l^\C$ which verifies
the conditions of Lemma 4.8.\qed
\enddemo
\bigskip
\subsubhead 6.3 Case when the contact form is not proportional
to any root
\endsubsubhead
\medskip
In  this case we obtain the following classification.\par
\bigskip
\proclaim{Proposition 6.4}
Let $(M = G/L, \D_Z)$  be a contact manifold
with   contact form $\vartheta$  not proportional to any root.
If it admits a primitive invariant CR structure $(\D_Z, J)$, then it is
one of the following.\par
If $G$ is simple then
\roster
\item"a)"  $G/L = SO_{2n}/SO_{2n-2}$, $n>2$,  and 
$\vartheta$ is either $\varepsilon_1$ or, when $n = 4$,
  $\varepsilon_1 + 
\varepsilon_2 + \varepsilon_3 \pm \varepsilon_4$; moreover the 
holomorphic subspace of the CR structure $(\D_Z, J)$ is given by 
$$\m^{10} = \m(\varepsilon_1 + \varepsilon_2) + 
t \m(\beta)\tag6.1$$
where $\beta = - \varepsilon_1 + \varepsilon_2$, 
$ -\varepsilon_3 - \varepsilon_4$ or $-\varepsilon_3 + \varepsilon_4$
(the last two cases  occur only for $n=4$)
and $t$ belongs to the punctured unit disc $D\setminus \{0\} \subset \C$;
\item"b)" $G/L = Spin_7/SU_3 = S(S^7) = S^7 \times S^6$, $\vartheta = 
\varepsilon_1 + \varepsilon_2 + \varepsilon_3$ and the holomorphic 
subspace of $(\D_Z, J)$ is given by 
$$\m^{10} = \m(\varepsilon_1 + \varepsilon_2) + 
t \m( - \varepsilon_3)
+ \overline{\m(\varepsilon_1 + \varepsilon_2)} + 
\frac{1}{t} \overline{\m( - \varepsilon_3)}\tag6.2$$
for some $t\in D\setminus \{0\}$.
\endroster
If $G$ is not simple then
\roster
\item"c)" $G/L = SU_2\times SU_2/T^1 = S(S^3) = S^3 \times S^2$, $\vartheta = 
(\varepsilon_1 - \varepsilon_2) - (\varepsilon_1' - 
\varepsilon_2')$ and the holomorphic subspace of $(\D_Z, J)$ is 
$$\m^{10} = \C (E_{\varepsilon_1 - \varepsilon_2} + t 
E_{\varepsilon'_1 - \varepsilon'_2}) + 
\C( E_{-(\varepsilon_1 -\varepsilon_2)} + \frac{1}{t}
E_{-(\varepsilon'_1 - \varepsilon'_2)})\ .\tag6.3$$
\endroster
In all cases we considered $\vartheta$ up to a factor and up to 
a transformation from the Weyl group $W(R)$, and $J$ up to a sign.
\endproclaim
\bigskip
\proclaim{Proposition 6.5} 
A homogeneous
contact manifold $(G/L, \D_Z)$  with contact 
form $\vartheta$ not proportional to any root, admits a non-standard
non-primitive CR structure if and only if it is  $G$-contact diffeomorphic
 to  the contact manifold
$(M(\Gamma) = G/L, \D_{Z(\Gamma)})$ associated  with a non-special CR-graph
$(\Gamma,\vartheta(\Gamma))$ 
(see Definition 1.7).\par
For any invariant CR structure $(\D_{Z(\Gamma)}, J)$ on $M(\Gamma) = G/L$ 
the natural projection $\pi: M(\Gamma) = G/L \to F_2(\Gamma) = G/Q$
is holomorphic w.r.t. the complex structure $J_2(\Gamma)$ or 
$-J_2(\Gamma)$.\par
The CR structures for which $\pi$ is holomorphic w.r.t.
$J_2(\Gamma)$ are in 1-1 correspondence with the invariant CR structures
on the fiber $C = Q/L$ subordinated to the induced contact structure
$\D_{Z(\Gamma)} \cap TC$.\par
More precisely, if 
$$\q^\C = \l^\C + \C Z + \m^{10}_C + \m^{01}_C\ ,
\ \ \g^\C = \q^\C + \m^{10}_{J_2} + \m^{01}_{J_2}$$
are the two decompositions of $\q^\C$ and 
$\g^\C$ associated with an invariant CR structure 
on the fiber $C = Q/L$
 and with the complex structure $J_2(\Gamma)$ on 
$F_2(\Gamma)$, then 
$$\m^{10} = \m^{10}_C + \m^{10}_{J_2}\tag6.4$$
is the holomorphic subspace of the corresponding CR structure on 
$M(\Gamma)$. Moreover, this CR structure is non-standard if and only if 
the CR structure on $C$ is primitive.
\endproclaim 
 \par
\bigskip
 The rest part of the paper is devoted to the proof of Propositions
6.4 and 6.5.  We need some additional notations.\par
 For a fixed
CR structure $(\D_Z, J)$,
we set
$$R^\pm_J = \{ \alpha \in R' \: J(E_\alpha) = \pm i E_\alpha\}\ ,\qquad
R_J = R^+_J \cup R^-_J\ ,\qquad R_\e \= R' \setminus R_J\tag6.5$$
and we define the subspaces
$$\m_J^\pm = \sum_{\beta \in R^\pm_J} \C E_{\beta}\ ,\qquad
\m_J = \m_J^+ + \m_J^-\ ,\qquad \e \= \sum_{\beta \in R_\e} 
\C E_{\beta}\subset \m^\C\ .\tag6.6$$
Note that $J$ is standard if and only if $R_J = R'$.
We define also the closed subsystem 
$$\tilde R_\e = [R_\e] \=  R \cap span_\R( R_\e )\ , 
\quad \tilde R_o = R_o \cap \tilde R_\e\ ,$$
and we set $R'_o = R_o \setminus \tilde R_o$.\par
The following Lemma collects some basic properties of these
objects.\par
\medskip
\proclaim{Lemma 6.6}
\roster 
\item $R_J = - R_J$ and $R_\e = - R_\e$;
\item for any $\alpha \in R_\e$ there exists exactly one root 
$\beta \in R_\e$
which is $\vartheta$-congruent to $\alpha$;
\item for  any pair $\alpha, \beta \in R_\e$ of $\vartheta$-congruent
roots,  there exist two uniquely 
determined complex numbers $\lambda, \mu\neq 0$
 such that
$$e_{\alpha, \beta} = E_\alpha + \lambda E_\beta \in \m^{10}\ ,\qquad
f_{\alpha, \beta} = E_\alpha + \mu E_\beta\in \m^{01}\ .\tag6.7$$
\item $(R^\pm_J + R_o)\cap R \subset R^\pm_J$ and
$(R_\e + R_o)\cap R \subset R_\e$;
\item $(R^\pm_J +  R_\e) \cap R \subset R^\pm_J \cup R_\e \cup R_o$.
\endroster
\endproclaim
\demo{Proof} (1) is clear.  To see (2), (3) and (4), observe
that $\alpha\in R_J$ if and only if $E_\alpha$  belongs to an 
irreducible $\k^\C$-module which is also $J$-invariant;
hence (2), (3) and (4) follow from Lemma 6.1 and
Corollary 3.10.\par
The proof of (5) is the following. Let $\gamma \in R^+_J$
and $\alpha,\beta \in R_\e$ a pair of two $\vartheta$-congruent
roots. If $\gamma + \alpha \in R^-_J$, consider the element
$f_{-\alpha, - \beta} \in \m^{01}$ as defined in (6.7). 
Since $E_{\gamma + \alpha} \in \m^{01}$, by the
integrability condition
$$[E_{\gamma + \alpha}, f_{-\alpha, -\beta}] = C E_{\gamma} + X\in \m^{01} +
\l^\C$$
for some $C\neq 0$ and $X \notin \C E_\gamma$. This implies
that $\gamma \in R^-_J$: contradiction. \qed
\enddemo
For any $\alpha \in R$, 
a   root $\beta \in R$, which is $\vartheta$-congruent to
$\alpha$, 
is said to be {\it $\vartheta$-dual to $\alpha$\/}
and we say that $(\alpha, \beta)$ is a {\it $\vartheta$-dual pair\/}.
By Corollary 3.10  any root admits {\it at most one\/} $\vartheta$-dual
root;  by Lemma 6.6 (3), {\it any root in $R_\e$\/} has {\it
exactly one\/} $\vartheta$-dual
root.
\par
\bigskip
\proclaim{Lemma 6.7} Let 
 $(\alpha,\alpha')$ be a $\vartheta$-dual pair in $R_\e$.
Then the 
root subsystem $\tilde R = R\cap span_\R\{\alpha, \alpha'\}$ is of type
$A_1 + A_1$.
In particular  $\alpha \perp \alpha'$ and $\alpha \pm \alpha' \notin R$.
\endproclaim
\demo{Proof} 
Assume that
$\tilde R \neq A_1 + A_1$. Then 
$\tilde R$ is a root system of type  $A_2, B_2$ or $G_2$. 
Since by assumptions
$\vartheta = \alpha - \alpha'$ is proportional to no root,
looking at the corresponding root systems, we find that up to 
a transformation from the Weyl group 
there are the following  possibilities:
\roster
\item" " $R = A_2$\ :\qquad $\alpha = \varepsilon_0 - \varepsilon_2$\ ,\quad
$\alpha' = \varepsilon_2 - \varepsilon_1$\ ;
\item" " $\tilde R = B_2$\ :\qquad $\alpha =\varepsilon_1$\ ,\quad
$\alpha' = -\varepsilon_1 + \varepsilon_2$\ ;
\item" "
$\tilde R = G_2$\ :\qquad $\alpha = -\varepsilon_2$\ ,\quad
$\alpha' = -\varepsilon_1 + \varepsilon_2$\ .
\endroster
Note that in each of these three cases, $\alpha + \alpha' = \beta \in R$.\par
\medskip
\noindent{\it Case $\tilde R = A_2$}.\par
 In this case
$\vartheta = (\varepsilon_0 - \varepsilon_2) - (\varepsilon_2
-\varepsilon_1) = \varepsilon_0 + \varepsilon_1 - 2 \varepsilon_2 $ and
$\beta = \alpha + \alpha'$ is orthogonal to $\vartheta$
 and hence it belongs to $R_o$.
Moreover $\l^\C = C_{\g^\C}(Z)$
contains the subalgebra
$$\l' = \C H_{\varepsilon_0 - \varepsilon_1}
+ \C E_{\varepsilon_0- \varepsilon_1} + 
\C E_{\varepsilon_1- \varepsilon_0}\ .$$
At the same time, by Lemma 6.6 (3),  $\m^{01}$ contains the element 
$f_{\varepsilon_0 - \varepsilon_2, \varepsilon_2 - \varepsilon_1}
= E_{\varepsilon_0 - \varepsilon_2} + \mu E_{
\varepsilon_2 - \varepsilon_1}$,
with some fixed $\mu \neq 0$. Since $\m^{01}$ is $\l^\C$-invariant,
$\m^{01}$ contains also the subspace
$$[E_{\varepsilon_1 - \varepsilon_0}, 
\C f_{\varepsilon_0 - \varepsilon_2, \varepsilon_2 - \varepsilon_1}]
= 
\C (E_{\varepsilon_1 - \varepsilon_2} - \mu E_{
\varepsilon_2 - \varepsilon_0})\ .$$
By integrability condition, this implies that 
$$[E_{\varepsilon_0 - \varepsilon_2} + \mu E_{
\varepsilon_2 - \varepsilon_1}, 
E_{\varepsilon_1 - \varepsilon_2} - \mu E_{
\varepsilon_2 - \varepsilon_0}] = 
 \mu(- H_{\varepsilon_0 - \varepsilon_2} + 
 H_{\varepsilon_2 - \varepsilon_1}) \in \m^{01}+\l^\C$$
 and hence we conclude that
$- H_{\varepsilon_0 - \varepsilon_2} + 
 H_{\varepsilon_2 - \varepsilon_1} \in \l^\C$.
But this cannot be because $- H_{\varepsilon_0 - \varepsilon_2} + 
 H_{\varepsilon_2 - \varepsilon_1}$
is not orthogonal to $Z = i\B^{-1}\circ \vartheta$.\par
\medskip
\noindent {\it Case $\tilde R = B_2$ or $G_2$\/}.\par
Then $\beta = \alpha + \alpha'$ is not orthogonal 
to $\vartheta = \alpha - \alpha'$ and, moreover, 
$$(\beta + \R \vartheta) \cap R = \emptyset\ .$$
These two facts show that $\beta \in  R \setminus (R_\e \cup R_o) = R_J$. 
Changing
the sign of $\alpha$ and $\alpha'$, if necessary, we may assume
that $\beta \in R_J^+$.\par
Consider the vector
$f_{\alpha, \alpha'} = E_\alpha + \mu E_{\alpha'} \in \m^{01}$ 
which is defined
by
(6.7). Then  $\overline{E_\alpha + \mu E_{\alpha'}}
= E_{-\alpha} + \bar\mu E_{-\alpha'} \in \m^{10}$ and by 
integrability condition its commutator with $E_\beta$
is also in $\m^{10}+\l^\C$. Therefore 
$$[E_{-\alpha} + \bar\mu E_{-\alpha'}, E_{\beta}]
= N_{-\alpha, \beta} E_{\alpha'} +
\bar \mu N_{-\alpha', \beta} E_{\alpha} \in 
\m^{10}\ .$$
Hence the coefficient $\lambda$ of the vector
$e_{\alpha, \alpha'}$ defined by (6.7) is 
$$\lambda = \frac{N_{-\alpha, \beta}}
{\bar\mu N_{-\alpha', \beta}}\ .\tag6.8$$
Since we use the Chevalley normalization (see \S 6.1),
$N_{-\alpha, \beta}= \pm (p+1)$ for any two roots $\alpha$, $\beta$, 
where $p\geq 0$ is the maximal integer such that 
$\beta + p\alpha \in \tilde R$ (see e.g. [7]).
Using this formula,  
 we obtain from (6.8) that if $\tilde R = B_2$, $\lambda \bar\mu = \pm 2$,
while if  $\tilde R = G_2$, $\lambda \bar\mu = \pm 3$.\par
On the other hand,  by integrability condition 
$$[e_{\alpha, \alpha'}, \overline {f_{\alpha, \alpha'}}]
= [E_\alpha + \lambda E_\alpha', E_{-\alpha} + \bar \mu E_{-\alpha'}] 
= H_\alpha + \lambda \bar\mu H_{\alpha'}
\in   \l^\C\ .$$
This means that $\vartheta(H_\alpha + \lambda \bar\mu H_{\alpha'}) = 0$, 
i.e. that 
$$ <\vartheta | \alpha>
+ \lambda \bar \mu <\vartheta| \alpha'> = 0\ ,$$
where $<\vartheta|\alpha> = 2 (\vartheta,\alpha)/(\alpha,\alpha)$.
Hence for $\vartheta = \alpha - \alpha'$, we obtain
$$2\  -  <\alpha'|\alpha> + 
\lambda\bar \mu [-2\ + <\alpha|\alpha'>] = 0\ .$$
In case $\tilde R = B_2$, $<\alpha'|\alpha> = - 2$ and 
$ <\alpha|\alpha'> = - 1$ so that $\lambda\bar \mu = 4/3$; in 
case $\tilde R = G_2$, $<\alpha'|\alpha> = - 3$ and 
$ <\alpha|\alpha'> = - 1$ so that $\lambda\bar \mu = 5/3$. In both cases
we get a contradiction with the previously determined values for
$\lambda \bar \mu$.\qed
\enddemo
\bigskip
  Now we determine the possible types of the root subsystem
$\tilde R_\e = R \cap span_\R(R_\e)$.\par
\medskip
\proclaim{Lemma 6.8} If   $\tilde R_\e$ is not of the 
form $A_1\cup A_1$, then $\tilde R_\e$ and $R$ are both indecomposable
root systems.
\endproclaim
\demo{Proof} 
By Lemma 6.7,  we 
may assume that $\text{rank}\ \tilde R_\e >2$.
Suppose that $\tilde R_\e$ is decomposable into two
mutually orthogonal subsystems $R_1$ and $R_2$. Let $\alpha \in 
 R_1 \cap R_\e$, $\alpha' \in R_2\cap R_\e$
and  $\beta$, $\beta'$ the  $\vartheta$-dual roots of $\alpha$ and
$\alpha'$, respectively. Since 
$\vartheta$ cannot be  in the span of $R_1$, it is clear 
that  $\beta\in R_2$ and  that  $\beta' \in R_1$. Then the identity
$$\R \vartheta = \R (\alpha - \beta) = \R (\alpha' - \beta')$$
implies that $\alpha + \rho \beta' = \rho\alpha' + \beta = 0$ for some 
$\rho \neq 0$.  From this follows
 that $\beta' = - \alpha$, $\beta = - \alpha'$
and that $\text{rank}\tilde R_\e = 2$: contradiction.\par
A similar contradiction arises if we replace $\tilde R_\e$ by $R$.
\qed\enddemo
\medskip
Note that by Lemma 6.8, if $G = G_1 \times G_2$, then the only possibility 
for $\tilde R_\e$ is $A_1 \cup A_1$.\par
\bigskip
The following Lemma gives a more detailed description of the 
root subsystem $\tilde R_\e$.\par
\medskip
\proclaim{Lemma 6.9} The root subsystem $\tilde R_\e$ has type $D_\ell$, 
$\ell >1$ or $B_3$ and, up to a factor and a transformation from the Weyl
group $W = W(R)$, the contact form $\vartheta$ is one of the following:
\roster
\item
if $\tilde R_\e = D_2 = A_1 + A_1'$ and $\alpha, \alpha'$
are roots  of   the summands $A_1$ and $A_1'$,
then $\vartheta = \alpha - \alpha'$;
\item if $\tilde R_\e = D_3$ or $D_\ell$, with $\ell >4$,
then $\vartheta = 2\varepsilon_1$;
\item if $\tilde R_\e = D_4$ then $\vartheta = 2\varepsilon_1$
or  $\vartheta = \varepsilon_1 + \varepsilon_2 + \varepsilon_3 +
\varepsilon_4$ or $\vartheta = \varepsilon_1 + \varepsilon_2 + \varepsilon_3 -
\varepsilon_4$;
\item if $\tilde R_\e = B_3$ then $\vartheta = 
 \varepsilon_1 + \varepsilon_2 + \varepsilon_3 $.
\endroster
\endproclaim
\medskip
Note that in case  $\tilde R_\e = D_4$, all three contact forms $\vartheta$ 
in (3) are 
equivalent with respect to automorphisms of the root system.\par
\medskip 
\demo{Proof} From Lemma 6.8, it is sufficient to consider the case when 
rank $\tilde R_\e >2$ and $\tilde R_\e$ is indecomposable. 
For each indecomposable root system $\tilde R_\e$ we describe, up to a 
transformation from the Weyl group, all  pairs of roots $(\alpha, \alpha')$,
which are orthogonal and such that $\alpha \pm \alpha' \notin R$. By 
Lemma 6.7 such pairs are the only candidates for $\vartheta$-dual pairs in
$R_\e$. For each case, we consider the corresponding form $\vartheta = \alpha-
\alpha'$,
and  describe all $\vartheta$-dual pairs in $\tilde R_\e$ . Then,
assuming that $\alpha, \alpha' \in R_\e$, we check if the case is possible
looking if the $\vartheta$-dual pairs in $R_\e$ may 
generate $\tilde R_\e$.\par 
\medskip
\noindent{\bf Case (A)}: $\tilde R_\e = A_\ell$.\par
Up to a transformation from the Weyl group, the pair $(\alpha, \alpha')$
is equal to
$(\varepsilon_1 - \varepsilon_2, \varepsilon_3 - \varepsilon_4)$.
Then $\vartheta
= (\varepsilon_1 - \varepsilon_2) - (\varepsilon_3 - \varepsilon_4)$ and the
$\vartheta$-dual pairs are (up to sign)
$$(\varepsilon_1 - \varepsilon_2, \varepsilon_3 - \varepsilon_4)\ ;\quad
(\varepsilon_1 - \varepsilon_3, \varepsilon_2 - \varepsilon_4)\ .$$
Since $\beta = \varepsilon_2 - \varepsilon_3 \in 
R_o = (\vartheta)^\perp \cap R$, then
$\varepsilon_1 - \varepsilon_3 = \alpha + \beta \in R_\e$
and hence also the second $\vartheta$-dual pair is in $R_\e$. 
In particular
$\text{rank}\tilde R_\e = 3$ and $\tilde R_\e = A_3 = D_3$.\par
\medskip
\noindent{\bf Case (B)}:  $\tilde R_\e = B_\ell$.\par
We have three possibilities 
for  $(\alpha, \alpha')$
according to their lengths:
\roster 
\item"i)" $(\alpha, \alpha') = (\varepsilon_1 + \varepsilon_2,
-(\varepsilon_3 + \varepsilon_4))$;
\item"ii)" $(\alpha, \alpha') = (\varepsilon_1 + \varepsilon_2, -
\varepsilon_3)$;
\item"iii)" $(\alpha, \alpha') = (\varepsilon_1, -\varepsilon_2)$.
\endroster
\smallskip
The last case is not possible, since we assume that 
$\vartheta = \alpha -\alpha'$
is 
proportional to no root.\par
\smallskip
\noindent{\rm \ i)\/}\ 
$\vartheta = \varepsilon_1 + \varepsilon_2 + \varepsilon_3 + \varepsilon_4$
and the $\vartheta$-dual pairs are (up to sign)
$$(\varepsilon_1 + \varepsilon_2, -(\varepsilon_3 + \varepsilon_4))\ ;\quad
(\varepsilon_1 + \varepsilon_3, -(\varepsilon_2 + \varepsilon_4))
\ ; \quad (\varepsilon_1 + \varepsilon_4, -(\varepsilon_2 + \varepsilon_3))\ .
\tag 6.9$$
As in case (A), one can check that all these $\vartheta$-dual pairs
are in $R_\e$ and that they span a space of dimension $4$.
Since the $\vartheta$-dual pairs consist of long roots, they cannot
generate the root system $B_\ell$ and hence this case is impossible.\par
\smallskip
\noindent{\rm \ ii)\/}\
 $\vartheta = \varepsilon_1 + \varepsilon_2 + \varepsilon_3$ and 
the $\vartheta$-dual pairs are (up to sign)
$$(\alpha = \varepsilon_1 + \varepsilon_2, \alpha ' = -\varepsilon_3)\ ;
\quad(\beta = \varepsilon_2 + \varepsilon_3, \beta' = -\varepsilon_1)
\ ;\quad(\gamma = \varepsilon_3 + \varepsilon_1, \gamma' = -\varepsilon_2)
\ .\tag6.10$$
Again all pairs in (6.10) consist of roots in $R_\e$.
This implies that $\text{rank}\tilde R_\e = 3$.\par
\medskip
\noindent{\bf Case (C)}: $\tilde R_\e = C_\ell$.\par
As in case (B), we have three possibilities.
\roster 
\item"i)" $(\alpha, \alpha') = (\varepsilon_1 + \varepsilon_2,
 -(\varepsilon_3 +\varepsilon_4))$;
\item"ii)"  $(\alpha, \alpha') = (\varepsilon_1 + \varepsilon_2,
 -2\varepsilon_3 )$;
\item"iii)"  $(\alpha, \alpha') = (2\varepsilon_1,
 -2\varepsilon_2 )$.
\endroster
\smallskip
As in  (B),  the last case is not possible.\par
\smallskip
\noindent{\rm \ i)\/}\ 
$\vartheta = \varepsilon_1 + \varepsilon_2 + \varepsilon_3 + \varepsilon_4$
and
the $\vartheta$-dual pairs are given (up to sign) in 
(6.9). This implies that $\pm 2\varepsilon_i \in R_J$, $i =1,\dots, 4$,
 because it has no
$\vartheta$-dual root and it is not orthogonal to $\vartheta$.
Note also  that the roots $\varepsilon_i - \varepsilon_j$, 
$i, j = 1, \dots, 4$,
belong to $R_o$, because they are orthogonal to $\vartheta$. Therefore
$$R_\e \subset 
\{\pm(\varepsilon_i + \varepsilon_j),\ , \ 1\leq i, j \leq 4\}
\subset (R_o + \{\pm 2\varepsilon_i\})\cap R \subset R_J$$
and this is a contradiction.\par
\smallskip
\noindent{\rm \ ii)\/}\
 $\vartheta = \varepsilon_1 + \varepsilon_2 + 2\varepsilon_3$.\par
In this case, up to sign, there is only one $\vartheta$-dual pair, 
that is $(\varepsilon_1 + \varepsilon_2, - 2\varepsilon_3)$. 
On the other hand, $\varepsilon_1 - \varepsilon_2 \in R_o$ and hence
$2\varepsilon_1 = (\varepsilon_1 + \varepsilon_2) + 
(\varepsilon_1 - \varepsilon_2) \in R_\e$: contradiction.
\par
\medskip
\noindent{\bf Case (D)}: $\tilde R_\e = D_\ell$.\par
Since $D_3 = A_3$, we may assume that $\ell \geq 4$.
Then
we have three possibilities:
\roster 
\item"i)" $(\alpha, \alpha') = 
(\varepsilon_1 + \varepsilon_2 , -(\varepsilon_3 +\varepsilon_4))$;
\item"ii)" $(\alpha, \alpha') = 
(\varepsilon_1 + \varepsilon_2,  - (\varepsilon_3 -\varepsilon_4)$;
\item"iii)" $(\alpha, \alpha') = 
(\varepsilon_1 + \varepsilon_2, -(\varepsilon_1 -\varepsilon_2))$.
\endroster
\smallskip
\noindent{\rm \ i)\/}\ 
$\vartheta = \varepsilon_1 + \varepsilon_2 + \varepsilon_3 + \varepsilon_4$
and
the $\vartheta$-dual pairs are given (up to sign) in 
(6.9) and they all belong
to $R_\e$. Hence the rank of $\tilde R_\e$ is $4$. \par
\smallskip
A similar argument shows that $\text{rank} \tilde R_\e = 4$ in case ii), 
where $\vartheta = \varepsilon_1 + \varepsilon_2 + \varepsilon_3 -
\varepsilon_4$.
\smallskip
\noindent{\rm \ iii)\/}\ $\vartheta = 2\varepsilon_1$ and the 
$\vartheta$-dual pairs are $(\varepsilon_1 + \varepsilon_i, 
\varepsilon_1 - \varepsilon_i)$, with $i = 2, \dots \ell$, they 
are all in $R_\e$  and they span the whole system $D_\ell$.\par
\medskip
\noindent{\bf Case (E)}:  $\tilde R_\e = E_6, E_7$ or $E_8$.\par
Let $\alpha, \alpha'\in R_\e$ be a $\vartheta$-dual pair. 
Since $\alpha$ and $\alpha'$ are 
orthogonal, we may included them into a subsystem $\Pi$ of simple roots.
According to the type of 
$\tilde R_\e$, 
without loss of generality, we may assume that $\alpha'$ is 
one of the following:
$$\alignat 2
\tilde R_\e &= E_6: \qquad &
\alpha' &= \varepsilon_4 + \varepsilon_5
+\varepsilon_6 + \varepsilon\ ; \\
\tilde R_\e &= E_7: \qquad 
&\alpha' &= \varepsilon_5 + \varepsilon_6
+\varepsilon_7 + \varepsilon_8\ ;\\
\tilde R_\e &= E_8:
\qquad
&\alpha'  &= \varepsilon_6
+ \varepsilon_7 + \varepsilon_8\ .
\endalignat$$
For each case, it follows that
 $\alpha = \varepsilon_i - \varepsilon_{i+1}$ for some
 $i \neq  \ell -3$ where 
$\ell = \text{rank}\tilde R_\e$. It can be easily checked
that, using permutations of the vectors $\varepsilon_i$
which belong to the Weyl group of $E_\ell$ and  which preserve
$\alpha'$, we may assume that either $\alpha = 
\varepsilon_1-\varepsilon_2$ or $\alpha =
 \varepsilon_{\ell-1}-\varepsilon_\ell
= -\sum_{i=1}^{\ell-2}\varepsilon_i - 2\varepsilon_{\ell}$.
Therefore we have the following  possibilities:\par
\noindent if $\tilde R_\e = E_6$:
\roster
\item"i)" $\alpha = 
\varepsilon_1 - \varepsilon_2$ and $\vartheta = \alpha' - \alpha
= - \varepsilon_1 + 
\varepsilon_2 + \varepsilon_4 + \varepsilon_5
+\varepsilon_6 + \varepsilon$;
\item"ii)"$\alpha = \varepsilon_5 - \varepsilon_6$
and $\vartheta = \varepsilon_4 + 2 \varepsilon_6 + \varepsilon = 
\varepsilon_6 - \varepsilon_1 - \varepsilon_2 - \varepsilon_3
-\varepsilon_5 + \varepsilon$;
\endroster
\par
\noindent if $\tilde R_\e = E_7$:
\roster
\item"iii)" $\alpha = 
\varepsilon_1 - \varepsilon_2$ and $\vartheta = 
- \varepsilon_1 + 
\varepsilon_2 + \varepsilon_5 + \varepsilon_6
+\varepsilon_7 + \varepsilon_8$;
\item" iv)"$\alpha = \varepsilon_6 - \varepsilon_7$
and $\vartheta =
\varepsilon_5 + 2 \varepsilon_7 + \varepsilon_8 = 
\varepsilon_7 - \varepsilon_1 - \varepsilon_2 - \varepsilon_3
-\varepsilon_4 - \varepsilon_6$;
\endroster
\par
\noindent if $\tilde R_\e = E_8$:
\roster
\item"v)" $\alpha = 
\varepsilon_1 - \varepsilon_2$ and $\vartheta =
- \varepsilon_1 + 
\varepsilon_2 + \varepsilon_6
+ \varepsilon_7 + \varepsilon_8$;
\item"vi)"$\alpha = \varepsilon_7 - \varepsilon_8$
and $\vartheta = 
\varepsilon_6 + 2 \varepsilon_8 = 
\varepsilon_8 - \varepsilon_1 - \varepsilon_2 - \varepsilon_3
-\varepsilon_4 - \varepsilon_5 - \varepsilon_7$.
\endroster
\par
We claim that all $\vartheta$-dual pairs belong to $R_\e$ and
that the space they generate
has dimension $5$ for the cases i),  ii) and v);
it has dimension $6$ for the cases iii) and  iv)
and dimension $7$ for the case vi). 
Since in all cases the dimension 
 is strictly less then $\text{rank} \tilde R_\e = \ell$,
we conclude that the case $R_\e = E_\ell$ is impossible.\par
We prove  the claim in the cases v) and vi) which occur when
 $\tilde R_\e = E_8$; 
in all other cases the proof is  similar.\par
For case v), the $\vartheta$-dual pairs are (up to sign)
$(-\varepsilon_1 + \varepsilon_i, \vartheta +\varepsilon_1 - \varepsilon_i)$,
where $i = 2, 6, 7, 8$
and they all belong to $R_\e$. These vectors generate a 5-dimensional vector 
space.
In  case  vi)
the $\vartheta$-dual pairs are
$(-\varepsilon_8 + \varepsilon_i, \vartheta + \varepsilon_8 - \varepsilon_i)$,
where $ i = 1, \dots, 5$ or $7$,
and again they are all in $R_\e$. 
These vectors generate a 7-dimensional vector 
space.\par
\medskip
\noindent{\bf Case (F)}:
$\tilde R_\e = F_4$.\par
We have the following possibilities:
\roster
\item"i)" $(\alpha, \alpha') = 
 (\varepsilon_1 + \varepsilon_2, -(\varepsilon_3 + \varepsilon_4))$;
\item"ii)" $(\alpha, \alpha') = 
(\varepsilon_1 + \varepsilon_2, -\varepsilon_3)$;
\item"iii)" $(\alpha, \alpha') = (\varepsilon_1 + \varepsilon_2,
-( 1/2(\varepsilon_1 - \varepsilon_2 + \varepsilon_3 +
\varepsilon_4))$;
\item"iv)"  $(\alpha, \alpha') = (\varepsilon_1, -\varepsilon_2)$.
\endroster
Cases i) and iv) are impossible because $\vartheta = \alpha-
\alpha'$ should be proportional to no root.
The admissible $\vartheta$-dual pairs for case ii) are given 
by (6.10)
and they all belong to 
$R_\e$. They generate a 3-dimensional subspace
and this is impossible because $\text{rank} \tilde R_\e = 
\text{rank} F_4 = 4$. A similar argument is applied for case iii).
\qed
\enddemo
\bigskip
\proclaim{Corollary 6.10} If $G$ is simple, then the only possibilities
for the pair $(R, \tilde R_\e)$ are
$$(A_n, A_3)\ ,\quad (A_n, B_3)\ ,\quad
(B_n, A_3)\ ,\quad (B_n, B_3)\ ,$$
$$(B_n, D_4)\ ,\quad (D_n, D_4)\ ,\quad (D_n, D_n)\ ,\quad (E_6, D_5)\ ,$$
$$
(E_7, D_6)\ ,\quad
(E_8,D_5)\ ,\quad (E_8,D_7)\ ,\quad (F_4,A_3)\ ,\quad (F_4,B_3)\ . $$
\endproclaim
\demo{Proof} If $R$ is the root system of the simple Lie group $G$ and
$(\alpha, \alpha')$ is a $\vartheta$-dual pair in $R_\e$, then the 
arguments used in the proof of Lemma 6.9  give the result.\qed
\enddemo
\bigskip
\proclaim{Lemma 6.11} Let $\tilde R_o = R_o \cap \tilde
R_\e$, $R'_o = R_o \setminus \tilde R_o$
 and $\alpha,\alpha' \in  R_\e$ be a $\vartheta$-dual
pair. Then 
\roster
\item"a)" $\tilde R_\e = [(\{\pm \alpha, \pm\alpha'\} + \tilde R_o)\cap R] 
\cup \tilde R_o$;
\item"b)" $R_\e = (\{\pm \alpha, \pm\alpha'\} + \tilde R_o)\cap R$ 
and $R_J \cap \tilde R_\e = \emptyset$;
\item"c)" $R_Q = R_o \cup R_\e$ and $R_{P} = R_o \cup R_\e \cup R^+_J$ 
are closed subsystem of $R$; $R_Q$ is the maximal symmetric
subset in $R_P$ (i.e. the biggest subset such that $-R_Q = R_Q$), and
 $R_{P}$
is parabolic (i.e. for any root $\alpha$, either $\alpha$ or
$-\alpha$ belongs to it);
\item"d)" $(R_Q + \tilde R_\e)\cap R \subset \tilde R_\e$ and hence $R_Q =
R'_o \cup
\tilde R_\e$ is an orthogonal decomposition;
\item"e)" for any $\vartheta$-dual pair $(\alpha, \alpha')$ let
$R_o(\alpha) = (R_o + \{\alpha\})\cap R$ and $R_o(-\alpha') = 
(R_o +\{-\alpha'\})\cap R)$; then the
set of roots
$$S(\alpha, \alpha') = R_o \cup R_o(\alpha) \cup R_o(-\alpha')
\cup R_J^+\tag6.11$$
is a closed parabolic subsystem of $R$.
\endroster
\endproclaim
\demo{Proof} a) When $\text{rank} \tilde R_\e = 2$
the claim is trivial.  \par
If $\tilde R_\e = B_3$, we may assume that $\alpha = \varepsilon_1 +
\varepsilon_2$,
 $\alpha' = - \varepsilon_3$ and $\vartheta = \varepsilon_1 + \varepsilon_2
+ \varepsilon_3$. Hence
$$\tilde R_o = \tilde R_\e \cap (\vartheta)^\perp = 
\{\varepsilon_i - \varepsilon_j
\ , \ i,j = 1, \dots, 3\}\ .$$
By Lemma 6.6 (4), 
$$(\{\pm\alpha, \pm\alpha'\} + \tilde R_o) \cap R = 
\{\pm(\varepsilon_i + \varepsilon_j), \pm \varepsilon_i
\ , \ i,j = 1, \dots, 3\} \subset R_\e\ .$$
Since $\tilde R_\e = \tilde R_o \cup 
\{\pm(\varepsilon_i + \varepsilon_j), \pm \varepsilon_i
\ , \ i,j = 1, \dots, 3\}$, the claim is proved for this case.\par
If 
$\tilde R_\e = D_\ell$, the argument is similar.
In particular, if $\vartheta = 2\varepsilon_1$,  one 
obtains that  $\tilde R_o = D_{\ell - 1} =
\{\pm \varepsilon_i \pm \varepsilon_j\ ,\ i,j >1\}$ and 
$R_\e = \{\pm \varepsilon_1 \pm \varepsilon_i\}$.\par
b) follows directly from a).\par
c)  The closeness of $R_Q$ and $R_P$ follows from Lemma 6.6 (4) and (5) and 
from point b). The last
statement is obvious.\par
d) The first claim follows from the facts that 
$\tilde R_\e = span_\R(R_\e) \cap
R$
and $(R_o + R_\e) \cap R \subset R_\e$. This implies
that $\g(\tilde R_\e)$ is an ideal of the semisimple Lie algebra
$\g(R_Q)$ and from this also the second claim follows.\par
e) By point b), $R_o(\alpha) \cup 
R_o(-\alpha') \subset R_\e$ and hence 
$R_o \cup R_o(\alpha) \cup 
R_o(-\alpha') \subset R_Q$ and $S(\alpha, \alpha') \subset R_P = R_Q \cup 
R_J^+$. Since $R_Q$ corresponds to a reductive part of the 
parabolic subalgebra $\g(R_P)$ and $R^+_J$ corresponds to 
the nilradical, it follows that $(S(\alpha, \alpha') + R_J^+) \cap R
\subset R_J^+$. By  d),  it remains to check that $R_o(\alpha)
\cup R_o(-\alpha)$ is a closed subsystem. \par
In case $\tilde R_\e = 2 A_1  = D_2$, we have that $R_o(\alpha) 
\cup R_o(-\alpha) = \{\alpha, -\alpha'\}$ and hence the claim is trivial.\par
In case $\tilde R_\e = 
D_\ell$, $\ell>2$, we may assume that $\vartheta = 2\varepsilon_1$, 
$\alpha = \varepsilon_1 + \varepsilon_2$, 
$\alpha' = -(\varepsilon_1 - \varepsilon_2)$. Then $\tilde R_o
= \{\pm \varepsilon_i \pm \varepsilon_j\ ,\ 1 <i, j\}$ and
$$R_o(\alpha) = \{ \varepsilon_1 \pm \varepsilon_i\ ,\ 1< i\} =
R_o(-\alpha') \tag6.12$$
and the conclusion follows. In case
$\tilde R_\e = 
B_3$, then $\vartheta = \varepsilon_1 + \varepsilon_2 + \varepsilon_3$, 
$\alpha = \varepsilon_1 + \varepsilon_2$ and $\alpha' = -\varepsilon_3$. Then
$\tilde R_o = \{\pm (\varepsilon_i - \varepsilon_j)\ \}$ and
$$R_o(\alpha) = \{ \varepsilon_i + \varepsilon_j\}\ ,\qquad
R_o(-\alpha')  = \{ \varepsilon_i \}\tag6.13$$
and again the conclusion follows.\qed
\enddemo
\medskip
Since $\g(R^+_J)$ is the nilradical of the parabolic subalgebra
$\g(R_P)$, we may choose an ordering of the roots such that 
the positive root system $R^+$ contains $R^+_J$. In the following
{\it $\alpha$ denotes the maximal root in $R_\e$ w.r.t. this ordering\/} 
and $\alpha'$ is its associated $\vartheta$-dual root.\par
\medskip
\proclaim{Proposition 6.12} The orthogonal complement $\m^\C$ to 
$\k^\C$ in $\g^\C$ admits the following $\k^\C$-invariant 
decomposition:
\roster
\item if $\tilde R_\e = B_3$ or $D_2 = A_1 + A_1$, then 
$$\m^\C = \e + \m^+_J + \m^-_J = (\m(\alpha) + \m(\alpha') + 
\overline{\m(\alpha)} + \overline{\m(\alpha')}) + \m^+_J + \m^-_J
\ ,\tag6.14$$
\item if $\tilde R_\e = D_\ell$, then 
$$\m^\C = \e + \m^+_J + \m^-_J = (\m(\alpha) + \m(\alpha'))+ \m^+_J + \m^-_J
\tag6.15$$
\endroster
where $\m(\alpha)$ and $\m(\alpha')$ are irreducible $\k^\C$-moduli
with highest weights $\alpha, \alpha'$, which are equivalent
and irreducible as $\l^\C$-moduli.\par
In terms of this decomposition, the holomorphic subspace $\m^{10}$
of the CR structure $(\D_Z, J)$  (up to sign) is  of the form
\roster
\item if $\tilde R_\e = B_3$ or $D_2 = A_1 + A'_1$
$$\m^{10} = (\m(\alpha) + 
t\m(\alpha') ) + (\overline{\m(\alpha)} + \frac{1}{t}\overline{\m(\alpha')})
+ \m^+_J\ ,
\tag6.16$$
\item if $\tilde R_\e = D_\ell$
$$\m^{10} = (\m(\alpha) + 
t\m(\alpha') ) 
+ \m^+_J
\tag6.17$$
\endroster
for some $t\in \{ x\in \C\ : 0 < |x| < 1\} = D\setminus\{0\}$.
\endproclaim
\demo{Proof} From $(R_o + \alpha)\cap R \subset R_\e$ and the definition of 
$\alpha$, the root $\alpha$
is the maximal weight of the $\k^\C$-module in $\m^\C$ which contains
$E_\alpha$. Moreover since $\alpha'$ is $\vartheta$-congruent to $\alpha$, 
then also $\alpha'$ is the maximal weight of an $\l^\C$-  and hence
$\k^\C$-module, and the $\l^\C$-moduli $\m(\alpha)$ and $\m(\alpha')$ are 
equivalent. By Lemma 6.11 b), it follows that the subspace $\e$,
spanned by the root vectors $E_\gamma$, $\gamma \in R_\e$, is given by
$$\e = \m(\alpha) + \m(\alpha') + 
\overline{\m(\alpha)} + \overline{\m(\alpha')}\ .$$
Moreover if  $\tilde R_\e = D_\ell$, $\ell>2$,   $R_o(\alpha) = R_o(-\alpha')$ 
(see (6.12)) and hence $\m(\alpha) = \overline{\m(\alpha')}$
(see also Table 3 in the Appendix).\par
>From Lemma 6.1 and the  remark in the second to the last point of 
\S 6.1, we obtain that 
the holomorphic subspace $\m^{10}$ is of the form
$$\m^{10} = (\m(\alpha) + 
t\m(\alpha') ) + \m^+_J$$
when $\tilde R_\e = D_\ell$, $\ell>2$, and of the form
$$\m^{10} = (\m(\alpha) + 
t\m(\alpha') ) + (\overline{\m(\alpha)} + s\overline{\m(\alpha')})
+ \m^+_J$$
when $\tilde R_\e = B_3$ or $D_2 = A_1 + A_1$, 
for some $t, s \neq 0$. By exchanging $\m^{10}$ with $\m^{01}$ (which 
corresponds to changing the sign of $J$) we may assume that 
$|t| \leq 1$. Using the integrability condition and the assumption
that $\vartheta = \alpha - \alpha' \notin R$, we have
$$[E_\alpha + t E_{\alpha'}, E_{-\alpha} + s E_{-\alpha'}]
= H_\alpha + ts H_{\alpha'} \in \m^{10} + \l^\C$$
and therefore $H_\alpha + ts H_{\alpha'}\in \l^\C$. Using (3.1) we get
$$0 = \vartheta(H_\alpha + ts H_{\alpha'}) = <\vartheta|\alpha> + 
ts <\vartheta|\alpha'>\ .$$
So 
$$ s = - \frac{1}{t} \frac{<\vartheta|\alpha> }{<\vartheta|\alpha'>}\ .$$
 If $\tilde R_\e = 2 A_1$, it is immediate to check that 
$<\vartheta|\alpha> = 2 = - <\vartheta|\alpha'>$. In case $\tilde R_\e = 
B_3$, we may assume that $\vartheta = \varepsilon_1 + \varepsilon_2
+ \varepsilon_3$, $\alpha = 
\varepsilon_1 + \varepsilon_2$ and
  $\alpha' = -\varepsilon_3$. Hence again 
$<\vartheta|\alpha> = - <\vartheta|\alpha'>$ 
and this shows that in both cases $s = 1/t$.\par
Finally, the condition $\m^{10} \cap \m^{01} = \{0\}$ implies
that the vectors $E_\alpha + t E_{\alpha'}$ and $\overline{E_{-\alpha}
+ \frac{1}{t} E_{-\alpha'}} = E_\alpha + \frac{1}{\bar t} E_{\alpha'}$
are linearly independent, and hence  $|t|\neq 1$.\qed
\enddemo
\bigskip
\proclaim{Lemma 6.13}
\roster
\item  Let $\q^\C = \k^\C + \e$ and $\p = \q^\C + \m^+_J$. Then
$\p$ is a parabolic subalgebra of $\g^\C$, with reductive part $\q^\C$ and 
nilradical $\m^+_J$. Moreover, if $Q$ is the connected subgroup of $G$ 
with Lie algebra $\q = \q^\C \cap \g$, then $F_2 = G/Q$ is a flag 
manifold and $\m^+_J$ is the holomorphic subspace of an invariant 
complex structure $J_2$ on $F_2 = G/Q$.
\item The subspace $\m^{10}_{J_1} = \m(\alpha) + \m(-\alpha') + \m^+_J$
is the holomorphic subspace of an invariant complex structure $J_1$
of $F_Z = G/K$.
\item The natural $G$-equivariant projections 
$$\pi: G/L \longrightarrow  G/Q\ ,\qquad \pi': G/K \longrightarrow G/Q$$
are holomorphic fibrations w.r.t. the CR structure $(\D_Z, J)$ on $G/L$,
the complex structure $J_1$ on $F_Z = G/K$ and the complex 
structure $J_2$ on 
$F_2 = G/Q$, respectively. Moreover, the typical fiber $C = Q/L$ of 
$\pi$ is either $Spin_7/SU_3 = S^7\times S^6$ or
$SO_{2\ell}/SO_{2 \ell -2}$, $\ell > 1$ and the induced invariant 
CR structure is primitive. 
\item The typical fiber $C = Q/L$ of $\pi$ may be equal  
$SO_4/SO_2 = S^3\times S^2$ only if
$G = G_1\times G_2$, with each $G_i$ simple.
\endroster
\endproclaim
\demo{Proof} (1) The proof follows from Lemma 6.11 c) and the remark that 
$\p = \g(R_P) + \h^\C$ and $\q^\C = \g(R_Q) + \h^\C$.\par
(2) We have to check the conditions a) and b) of (4.2). Condition a) is 
obvious. Condition b) means that $\k^\C + \m(\alpha) + \m(-\alpha') =
\g(S(\alpha, \alpha')) + h^\C$ is a subalgebra. This follows from Lemma 6.11 e).
\par
(3) The first claim follows from Lemma 4.8. \par
For the second claim, we recall that 
 we have the following decompositions of the Lie algebras
$\q^\C$ and $\l^\C$:
$$\l^\C = \g(R'_o) \oplus
\left( \g(\tilde R_o) + Z(\l^\C)\right)\ ,$$
$$\q^\C = \k^\C + \e = \g(R'_o)
\oplus \left( \g(\tilde R_\e) + Z(\q^\C)\right)\ .$$
Since the fiber $Q/L$ has a non-standard CR structure, the group 
$Q' = Q/N$, where $N$ is its kernel of non-effectivity, is semisimple
by Corollary 3.2 and Proposition 4.6. Therefore it has Lie algebra
$\q'{}^\C = \g(\tilde R_\e) = B_3$ or $D_\ell$. The corresponding 
stability subalgebra $\l'{}^\C = \l^\C/\n^\C$ has rank equal to 
$\text{rank}(\q'{}^\C) - 1$ and his semisimple part is $\g(\tilde R_o) = 
A_2$ or $D_{\ell-1}$. Hence the fiber $Q/L = Q'/L'$, considered as homogeneous
manifold of the effective group $Q'$, is 
either $Spin_7/SU_3$ or $SO_{2\ell}/SO_{2\ell - 2}$ (note that $SO_7$
does not contains $SU_3$). The manifold $Spin_7/SU_3$ can be identified
with the unit sphere bundle $S(Spin_7/G_2) = S(S^7) = S^7 \times S^6$.\par
The holomorphic subspace $\m^{10}(Q/L)$ of the CR structure of
the fiber $Q/L$ is of the form 
$(\m(\alpha) + t \m(\alpha')) + (\overline{\m(\alpha)}
+ 1/t \overline{\m(\alpha')})$ for some $t\neq 0$ and the minimal 
$\k^\C$-module generated by $\m^{10}(Q/L)$ is $\e$. By Lemma 4.8, 
this implies  that the CR structure on $Q/L$ is primitive.\par
(4) It is sufficient to observe that if $G$ is simple, the case
$\tilde R_\e =  A_1 \cup A_1$ cannot occur by Corollary 6.10.
\qed
\enddemo
\bigskip
Lemma 6.13 (3) and 
Proposition 6.12
directly imply
Proposition 6.4.\par
\bigskip
Now it remains to prove Proposition 6.5. Let $(M = G/L, \D_Z, J)$
be a non-standard non-primitive 
CR manifold with contact form $\vartheta$
not proportional to any root. We recall that in Lemma 6.13 (3)
we  defined a complex structure  $J_1$
on the flag manifold $F_Z = G/K$, associated with the decomposition
$\g^\C = \k^\C + \m^{10}_{J_1} +  \m^{01}_{J_1}$.
 We  also defined  another flag manifold $F_2 = G/Q$, 
with $\q^\C = \k^\C + \e$, with invariant complex structure $J_2$ 
associated with the
decomposition $\g^\C = \q^\C + \m^+_J + \m^-_J$ and such that the 
projection $\pi: (F_Z = G/K, J_1) \to (F_2 = G/Q, J_2)$ is holomorphic.
Moreover the CR structure $(\D_Z, J)$ on $G/L$ has the holomorphic subspace 
defined in (6.16) and (6.17).\par
The subalgebra $\k^\C$ corresponds to the root subsystem $R_o$, which has the
orthogonal decomposition $R_o = R_o' \cup \tilde R_o$, and $\q^\C$ corresponds
to the root subsystem with the orthogonal decomposition $R_Q = R_o' \cup
\tilde R_\e = R_o' \cup (\tilde R_o \cup R_\e)$ (see Lemma 6.11). Moreover
there are only three possibilities for the pair of subsystems 
$(\tilde R_\e, \tilde R_o)$
, namely $(D_2 = 2A_1,\emptyset,)$, $( D_\ell),D_{\ell-1} ,$, $\ell >2$,
or $(B_3, A_2 )$. However, the following lemma
shows that this last case cannot occur.\par
\bigskip
\proclaim{Lemma 6.14} If $R_J\neq \emptyset$, then $\tilde R_\e \neq B_3$.
\endproclaim
In other words, the fiber $C = Q/L$ of the CRF fibration 
$\pi: G/L \to G/Q$ described in Lemma 6.13 (3) cannot be $Spin_7/ SU_3$
if the base is not trivial.\par
\demo{Proof} Assume that $\tilde R_\e = B_3$. Then $G$ is simple
and $R$ is indecomposable by Lemma 6.8. So $R$ has type either $B_n$ or
$F_4$, because these are the only connected Dynkin graphs
which contain a subgraph of type $B_3$. \par
If $R = F_4$, using
the notation of the Appendix, we may assume that 
$(\alpha = \varepsilon_2 + \varepsilon_3, \alpha' = -\varepsilon_4)$
is a $\vartheta$-dual pair in $R_\e$.
Since $\vartheta =  \varepsilon_2 + \varepsilon_3
+ \varepsilon_4$, then $-\varepsilon_4 + \varepsilon_1 \in R_J$, because
it is not orthogonal to $\vartheta$ nor has a $\vartheta$-dual root; moreover
$-\varepsilon_1 \in R_o = R \cap (\vartheta)^\perp$ and hence
$-\varepsilon_4 =  (-\varepsilon_4 + \varepsilon_1)  -\varepsilon_1\in 
R_J$, by Lemma 6.6 (4): contradiction.\par
Assume now that $R = B_n$, $n>3$. Then we may assume that 
$(\alpha, \alpha') = (\varepsilon_1 + \varepsilon_2, - \varepsilon_3)$
is a $\vartheta$-dual pair in $R_\e$ and hence that $\vartheta = 
\varepsilon_1 + \varepsilon_2 + \varepsilon_3$. Then, as before, we get that 
$-\varepsilon_3 + \varepsilon_4 \in R_J$, $-\varepsilon_4 \in R_o$ and 
hence that $-\varepsilon_3 = (-\varepsilon_3 + \varepsilon_4) +
(-\varepsilon_4)  \in R_J$: contradiction.\qed
\enddemo
 \bigskip
Now we construct some special basis $\Pi$ for $R$, 
which we will call {\it good\/}. For any basis $\Pi$ let
$$\Pi_o = \Pi \cap R_o\ ,\quad \tilde \Pi_o = \Pi_o \cap \tilde R_\e \ ,\quad
\tilde \Pi_\e = \Pi \cap \tilde R_\e
\ ,\quad \Pi_\e = \Pi \cap R_{\e}\ .$$
 Then $$\tilde \Pi_\e = \Pi_\e \cup \tilde \Pi_o\ .$$
A basis $\Pi$ is called {\it good\/} if
$$\tilde R_o = [\tilde \Pi_o]\ ,\quad \tilde R_\e = [\tilde \Pi_\e]\ ,
\quad R_o = [\Pi_o]\ ,$$
where for any subset $A\subset \Pi$ we denote $[A] = span(A) \cap R$.\par
\medskip
A good basis exists because $R_o \cup \tilde R_\e = R_o' \cup
\tilde R_\e$ is a closed subset of roots, $R'_o$ is orthogonal to 
$\tilde R_\e$ and $R_o = R'_o \cup (R_o \cap \tilde R_\e) = R'_o \cup
\tilde R_o$. In fact, we may take a basis $\tilde \Pi_o$ for $\tilde R_o$, 
extend it to a basis $\tilde \Pi_\e$ for $\tilde R_\e$, add to it 
a basis for $R'_o$ and finally extend everything to a basis $\Pi$ for
$R$.\par
By the remarks before Lemma 6.14, the pair $(\tilde \Pi_\e, \tilde \Pi_o)$
is of type $(D_\ell, D_{\ell-1})$, $\ell >2$, or $(2 A_1, \emptyset)$
and it can be represented by the following two graphs
$$\dynkin\xroot{\ \ 2}\link
\wroot{\ \ 2}\link\dots\link\wroot{\ \ 2}
\link\wroot{\ \ 2}\wrootupright{1}
\wrootdownright{1}\enddynkin \tag6.18$$ 
$$\dynkin\xroot{\ \ 1}\enddynkin \quad
\dynkin\xroot{\ \ -1}\enddynkin
\tag6.19$$
where the subdiagram of $\tilde \Pi_o$ is obtained by deleting the grey nodes.
Moreover, by Lemma 6.9, the contact form $\vartheta$ is the linear
combination of the simple roots associated with the nodes of (6.18)
and (6.19) with the indicated coefficients. For example, if 
$(\tilde \Pi_\e, \tilde \Pi_o) = (D_\ell, D_{\ell-1})$ and if we use
the standard correspondence between nodes and roots, we get
$$
\vartheta =
2(\varepsilon_1 - \varepsilon_2) + \dots + 
2(\varepsilon_{\ell - 2} - \varepsilon_{\ell-1}) + 
(\varepsilon_{\ell-1} - \varepsilon_\ell) + 
(\varepsilon_{\ell-1} + \varepsilon_\ell) = 2\varepsilon_1\ .$$
Note that if $\ell = 4$, using two permutations of the simple roots 
corresponding
to the end nodes, one gets the  other 
two possible contact forms,
namely
$\vartheta = \varepsilon_1 + \varepsilon_2 + 
\varepsilon_3 + \varepsilon_4$ and  $\vartheta = 
\varepsilon_1 + \varepsilon_2 + 
\varepsilon_3 -\varepsilon_4$.\par
Remark that a good basis $\Pi$ together with the subsets $\Pi_o$ and $\Pi_\e$
completely determines the homogeneous CR manifold $M = G/L$ and the 
flag manifolds $(F_Z = G/K, J_1)$ and $(F_2 = G/Q, J_2)$. In fact the root
systems $R_o = R(K)$ of $K$ and $R(Q)$ of $Q$ are given by 
$R_o = [\Pi_o]$ and $R(Q) = [\tilde \Pi_e = \Pi_o \cup \Pi_\e]$ and 
$\l= \k \cap (ker \vartheta)$, where $\vartheta$ is 
defined by (6.18)-(6.19).\par
Notice also that by definition of  good basis
$$R\cap (\vartheta)^\perp = R_o = [\Pi_o]\tag6.20$$
and hence that 
$$\Pi_o = \Pi \cap (\vartheta)^\perp\ . \tag6.21$$
Any good basis $\Pi$ together with the subsets $\Pi_o$ and $\Pi_\e$ can be 
represented by a painted Dynkin graph $\Gamma = \Gamma(\Pi)$ if we paint
the nodes corresponding to the roots of $\Pi_\e$ in grey, the nodes
of $\Pi_o$ in white and all  others in black.\par
 We call such graph 
$\Gamma$ {\it a painted Dynkin graph
associated with the CR manifold $(M = G/L, \D_Z, J)$\/}.\par
Any associated painted Dynkin graph has the following two properties.\par
\roster
\item It contains a unique proper subgraph $\Gamma_\e$ of type
(6.18), if it is connected, or of type (6.19), if it is not
connected; moreover in this second case, $\Gamma = \Gamma_1 \cup \Gamma_2$
has two connected components and each of them contains exactly 
one grey node.
\item The black nodes are exactly the nodes which are linked to 
$\Gamma_\e$.
\endroster
Indeed, (1) follows from definition of good basis,  Lemma 6.8 and Lemma 6.2.
(2) follows from (6.21).\par
\medskip
A painted Dynkin graph which verifies (1) and (2) is called {\it 
admissible graph\/}.\par
\medskip
Let $\Gamma$ be an admissible graph and $\Gamma_\e$ the corresponding 
subgraph of type (6.18) or (6.19). We denote by $\vartheta(\Gamma)$ the linear
combination of roots associated with the nodes of $\Gamma_\e$ as prescribed
in (6.18)-(6.19).\par
An admissible graph $\Gamma$ is called {\it good\/} if 
$$[\Pi_o] = R\cap (\vartheta)^\perp\ ,\tag6.22$$
where $\Pi_o$ is the set of simple roots associated with the white nodes of 
$\Gamma$. Remark that by (6.20) any graph associated with 
$(M = G/L, \D_Z, J)$ is a good graph. The converse of this statement
is also true.\par
\bigskip
\proclaim{Lemma 6.15} Any good graph
is a painted Dynkin graph associated to 
a homogeneous CR  manifolds $(G/L, \D_Z,J)$, which have a
contact form $\vartheta$ parallel to no roots and 
where $(\D_Z, J)$  is non-standard and
non-primitive.
\endproclaim
\demo{Proof} Let $\Gamma$ be a good graph 
and $\vartheta(\Gamma)$ the corresponding 
contact form. As described in the Introduction,
$\Gamma$ defines two flag manifolds
$F_1(\Gamma) = G/K$ and $F_2(\Gamma) = G/Q$, with invariant complex
structures $J_1(\Gamma)$ and $J_2(\Gamma)$, respectively. Denote by
$$\g^\C = \k^\C + \m^{10}_{J_1} + \m^{01}_{J_1}\ ,\qquad
\g^\C = \q^\C + \m^{10}_{J_2} + \m^{01}_{J_2}$$
the corresponding associated decompositions. Consider also
the element $Z = i\B^{-1}\circ \vartheta(\Gamma)$.
Since the 1-parametric subgroup generated by $Z$ is closed, by Proposition 3.3
it defines a contact manifold 
$(M = G/L, \D_Z)$ with $\l = 
 \k \cap (Z)^\perp$. Moreover the fiber $C = Q/L$ of the 
fibration $\pi:  G/L \to G/Q$,  together with the contact structure 
induced  on $C$ by $Z$, is one of the contact manifolds 
described in Proposition
6.4
admitting a primitive CR structure. \par
If $\m^{10}_C$ is the holomorphic
subspace of such CR structure, then $\m^{10} = 
\m^{10}_C + \m^{10}_{J_2}$ is the holomorphic subspace of 
a non-standard CR structure
$(\D_Z, J)$ on $G/L$ and  the associated painted Dynkin graph is exactly 
$(\Gamma, \vartheta(\Gamma))$.
In fact, the conditions i) and ii) of Definition 4.1 are immediate. 
The integrability condition follows from the fact that $\m^{10}_C$ is a
holomorphic
subspace for a CR structure on $C = Q/L$ (and hence that $\l^\C + \m^{10}_C$
is a subalgebra), 
that $\m^{10}_{J_2}$ is the nilradical of the parabolic subalgebra
$\q^\C + \m^{10}_{J_2}$, and that $\m^{10}_C \subset \q^\C$.
\qed
\enddemo
\medskip
Now the classification of  homogeneous CR manifolds of the considered
type reduces to the classification of good graphs  $\Gamma$. 
\medskip
\noindent{\it Case 1\/}.\ $\Gamma$ is not connected.\par
In this case $\Gamma_\e = A_1 \cup A_1$, $\Gamma = \Gamma_1 \cup \Gamma_2$, 
where each $\Gamma_i$ is a connected component which corresponds to  a
root system $R_i$, and $R = R_1 \cup R_2$. Moreover $\vartheta = 
\alpha_1 - \alpha_2$, where $\alpha_i \in R_i$.\par
We prove that if $\Gamma$ is good then $R = A_p \cup A_q$, with 
$p+q >1$ and that $\Gamma$
is a CR-graph of type II.\par
First of all, one can easily check that if one of the connected components 
$\Gamma_i$ is 
not of type $A_q$, then $\Gamma$ is not good, that is that there exists
a root $\beta \in R\cap (\vartheta)^\perp$ which is not in $[\Pi_o]$. 
For example, if $R_1 = D_q$, we may assume that $\vartheta = \alpha_1 - 
\alpha_2$, where $\alpha_1 = \varepsilon_1 - \varepsilon_2$. Then
$\beta = \varepsilon_1 + \varepsilon_2 \in R\cap (\vartheta)^\perp$
but it is not in $[\Pi_o]$.\par
Assume now that  $R = A_p \cup A_q$. Without loss of generality
we may assume that $\alpha_1 = \varepsilon_{k} - \varepsilon_{k+1}$, 
$\alpha_2 = \varepsilon'_{r} - \varepsilon_{r+1}'$ are the roots associated 
with the grey nodes of $\Gamma_1$ and $\Gamma_2$, respectively. Then
$R\cap (\vartheta)^\perp = A_{p-2} \cup A_{q-2}$ and  it coincides
with $[\Pi_o]$ if and only if the nodes of the roots $\alpha_i$ are end nodes.
This proves that $\Gamma$ is good if and only if it is  a CR-graph 
of type II ( see Definition 1.7). \par
\medskip
\noindent{\it Case 2\/}.\ $\Gamma$ is  connected.\par
In this case, $\Gamma$ is a good graph only if the type of the pair
$(\Gamma, \Gamma_\e)$ is one of the following
$$(A_n, A_3)\ ,\quad (B_n,A_3)\ ,\quad (D_n, D_4)\ ,\quad (E_6, D_5)\ ,$$
$$(E_7, D_6)\ ,\quad (E_8, D_5)\ ,\quad (E_8, D_7)\ .$$
This follows from Corollary 6.10 and the fact that 
$(A_n,A_3)$, $(B_n, B_3)$, $(B_n, D_4)$, 
$(D_n,D_n)$, $(F_4,A_3)$ and $(F_4, B_3)$ do not correspond
to any admissible graph.\par
\medskip
We first prove that the cases $(B_n, A_3)$, $(E_7,D_6)$, 
$(E_8, D_5)$ and $(E_8, D_7)$
are not possible.\par
\roster
\item"i)" $(\Gamma, \Gamma_\e) = (B_n, A_3)$. In this case $\Gamma$
is of the form
$$\dynkin
\wroot{}\link \dots \link\wroot{}\link\root{}
\link\wroot{\ \alpha_1}\link\xroot{\ \alpha_2}\link\wroot{\ \alpha_3}
\link\root{}\link \dots \link\wroot{}\llink>\wroot{}
\enddynkin$$
where $\alpha_1 = \varepsilon_{k} - \varepsilon_{k+1}$, 
$\alpha_2 = \varepsilon_{k+1} - \varepsilon_{k+2}$ and 
$\alpha_3 = \varepsilon_{k+2} - \varepsilon_{k+3}$. Then 
$\vartheta(\Gamma) = \alpha_1 + 
2 \alpha_2 + \alpha_3 = \varepsilon_k + \varepsilon_{k+1} - \varepsilon_{k+2}
- \varepsilon_{k+3}$ and 
$$ \Pi_o= \{ \varepsilon_i -\varepsilon_{i+1}, i=1,\dots,k-2;k;k+2;k+4,\dots ,
n-1;\varepsilon_n \}\ . $$
However the root $\beta = \varepsilon_{k+1} + \varepsilon_{k+2}
\in (\vartheta(\Gamma))^\perp \cap R$ but it does not belong to $[\Pi_o]$:
contradiction.
\item"ii)" $(\Gamma, \Gamma_\e) = (E_7, D_6)$. In this case $\Gamma$
and $\vartheta(\Gamma)$ are 
$$
\dynkin\xroot{\ \alpha_1}\link\wroot{\ \alpha_2}\link\wroot{\ \alpha_3}
\link\wroot{\ \alpha_4}
\wrootdown{\alpha_7}\link\wroot{\ \alpha_5}
\link\root{\ \alpha_6}
\enddynkin\ ,\quad\vartheta(\Gamma) 
 = 2\varepsilon_1 +
\varepsilon_7 + \varepsilon_8\ .$$
However, this situation corresponds to no good graph,  
because  the root $\beta = \varepsilon_7 - \varepsilon_8$ is in 
 $\vartheta(\Gamma)^\perp\cap R$, but it does not belong to 
$$[\Pi_o]  = [\{ \alpha_2, 
\alpha_3, \alpha_4, \alpha_5, \alpha_7\}] = 
\{\ \varepsilon_a - \varepsilon_b\ ,
\ \pm( \varepsilon_a + \varepsilon_b + \varepsilon_7 + \varepsilon_8)\ ,
\ 1 \leq a,b \leq 6\ \}\ .$$
\par
\medskip
\item"iii)" $(\Gamma, \Gamma_\e) = (E_8, D_5)$. Then $\Gamma$ and 
$\vartheta(\Gamma)$ are 
$$
\dynkin\wroot{}\link\wroot{}\link\root{}\link\wroot{}\link\wroot{}
\wrootdown{}\link\wroot{}
\link\xroot{}
\enddynkin
\ ,\quad\vartheta(\Gamma) = \varepsilon_4 +  \varepsilon_5 + 
\varepsilon_6 + \varepsilon_7 -
 \varepsilon_8\ .$$
One can easily check that the root $\beta = \varepsilon_1 + 
\varepsilon_2 + \varepsilon_4 $ is orthogonal to $\vartheta(\Gamma)$,
but it doesn't belong to the subsystem 
$$[\Pi_o] = [\{ \alpha_1, \alpha_2, \alpha_4, \alpha_5, \alpha_6, \alpha_8\}]$$
generated by  white roots: contradiction.
\item"iv)" $(\Gamma, \Gamma_\e) = (E_8, D_7)$.Then $\Gamma$ and 
$\vartheta(\Gamma)$ are 
$$\dynkin\xroot{}\link\wroot{}\link
\wroot{}\link\wroot{}\link
\wroot{}\wrootdown{}\link
\wroot{}\link\root{}\enddynkin
\ ,
\quad
\vartheta(\Gamma) =  2\varepsilon_1 +  \varepsilon_8\ .$$
Also this case is not possible because $\varepsilon_7 - \varepsilon_9
\in R\cap (\vartheta(\Gamma))^\perp$ but it is not in 
$[\Pi_o] = [\{ \alpha_2, \alpha_3, \alpha_4, \alpha_5, \alpha_6, \alpha_8\}]$.\par
\endroster
\medskip
It remains to describe the good graphs of the following types
$$1) \ (A_n, A_3)\ ,\qquad 2)\ (D_n, D_4)\ ,\qquad
3)\ (E_6, D_5)\ .$$
\par
\noindent{\it (1) $(\Gamma, \Gamma_\e) = (A_n, A_3)$\/}.\par
Assume that $\Gamma_\e$ is not at an end of $\Gamma$, that is 
$$\dynkin\wroot{}\link\dots\link\wroot{}\link\root{}
\link\wroot{\ \alpha_1}\link\xroot{\ \alpha_2}\link\wroot{
\ \alpha_3}\link\root{}\link\wroot{}\link
\dots\link\wroot{}
\enddynkin$$
Then we may assume that $\alpha_i = \varepsilon_{p+i} - \varepsilon_{p+i+1}$. 
Then $\vartheta(\Gamma) = \alpha_1 + 2 \alpha_2 + \alpha_3 = 
\varepsilon_{p+1} + \varepsilon_{p+2} -
\varepsilon_{p+3} - \varepsilon_{p+4}$ and the root $\beta = 
\varepsilon_p - \varepsilon_{p+5} \in R \cap (\vartheta(\Gamma))^\perp$
but it is not in the span of $\Pi_o$; hence the graph is not good. On the 
other hand one can easily check that the graph
$$\dynkin \wroot{}\link\xroot{}\link\wroot{}\link\root{}\link
\wroot{}\link\dots\link\wroot{}\enddynkin$$
is good.\par
\medskip
\noindent{\it (2) $(\Gamma, \Gamma_\e) = (D_n, D_4)$\/}.\par
In this case we have two admissible graphs:
$$\dynkin
\wroot{}\link\dots\link\wroot{}\link\root{}\link\wroot{\ \alpha_{n-3}}
\link\wroot{\ \alpha_{n-2}}\xrootupright{\alpha_{n-1}}
\wrootdownright{\alpha_{n}}
\enddynkin\tag6.23$$
$$\dynkin
\wroot{}\link\dots\link\wroot{}\link\root{}\link\xroot{\ \alpha_{n-3}}
\link\wroot{\ \alpha_{n-2}}\wrootupright{\ \alpha_{n-1}}\wrootdownright{
\ \alpha_{n}}
\enddynkin\tag6.24$$
Using the standard equipment, we have that if $\Gamma$ is given by (6.23), 
then 
$$\vartheta(\Gamma) = \alpha_{n-3} + 
2 \alpha_{n-2} + \alpha_{n-1} + 2\alpha_n 
 = \varepsilon_{n-3} + \varepsilon_{n-2} + \varepsilon_{n-1} + 
\varepsilon_n$$
and $R\cap (\vartheta(\Gamma))^\perp = D_{n-4} \cup A_3$. If 
$\Gamma$ is given by (6.24), then 
$$\vartheta(\Gamma) = 2\alpha_{n-3} + 
2 \alpha_{n-2} +  \alpha_{n-1} + \alpha_n = 2\varepsilon_{n-3}$$
and $R\cap (\vartheta(\Gamma))^\perp = D_{n-1}$. \par
Since in both cases $[\Pi_o] = A_{n-5} \cup A_3$, the graph (6.24) is not 
good, while the graph (6.23) is good only when $n = 5$. \par
\medskip
\noindent{\it (3) $(\Gamma, \Gamma_\e) = (E_6, D_5)$\/}.\par
Up to isomorphism, we have only one admissible graph
$$\dynkin\xroot{\ \alpha_1}\link\wroot{\ \alpha_2}\link
\wroot{\ \alpha_3}\wrootdown{\alpha_6}\link
\wroot{\ \alpha_4}\link\root{\ \alpha_5}\enddynkin$$
Using the standard equipment, we get 
$$\vartheta(\Gamma) = 2(\alpha_1 + \alpha_2 + \alpha_3) 
+ \alpha_4 + \alpha_6 = 2\varepsilon_1 + \varepsilon_6 + \varepsilon\ .$$
Then 
$$R\cap (\vartheta(\Gamma))^\perp = 
\{\ \varepsilon_a - \varepsilon_b\ ,
\ \pm(\varepsilon_a + \varepsilon_b + \varepsilon_6 + \varepsilon)\ ,
\ a,b = 2,3, 4, 5\ \} = [\Pi_o] = D_4$$
and hence the graph is good.\par
\medskip
This concludes the classification of good graphs. We proved that 
the pairs $(\Gamma, \vartheta(\gamma))$ given by all good graphs are 
exactly the non-special CR-graphs of Definition 1.7. By the remarks 
before Lemma 6.15 and the Lemmata 6.13 and 6.15, Proposition 
6.5 follows.\par
\bigskip
\bigskip
\newpage
\head APPENDIX
\endhead
\bigskip
The notation used in the following Tables is the same 
of [7]. We recall that
the weights of the groups $B_\ell, C_\ell, D_\ell$ and $F_4$
are expressed in terms of an 
orthonormal basis $(\varepsilon_1, \dots, \varepsilon_\ell)$
of $\goth{h}(\Bbb Q)^*$. The weights of the groups $A_\ell,
E_7, E_8$ and $G_2$ are expressed
using vectors $\varepsilon_1, \dots, \varepsilon_{\ell+1}\in 
\goth{h}(\Bbb Q)^*$
such that 
$$
\sum \varepsilon_i = 0\ ,\qquad
(\varepsilon_i, \varepsilon_j) = 
\left\{\matrix \frac{\ell}{\ell+1} & i= j\\
-\frac{1}{\ell+1} & i\neq j
\endmatrix \right.\tag A.1$$
It is useful to recall that  if $\sum a_i = 0$,
then $(\sum a_i \varepsilon_i, \sum b_j \varepsilon_j) =
\sum a_i b_i$.
For $E_6$, the weights are expressed 
by  vectors $\varepsilon_1, \dots, \varepsilon_6$,
which verify (A.1) with $\ell = 5$, and by an auxiliary vector $\varepsilon$
which  is orthogonal
to all $\varepsilon_i$ and verifies $(\varepsilon, \varepsilon) = 1/2$. \par
\bigskip
In Table 1,  for any simple complex Lie group 
$\g^\C$, we give the corresponding root system $R$, 
the longest root $\mu$ (unique up to inner automorphisms), 
the 
 subalgebra $\g'_0 = C_{\g^\C}(\g(\mu))$,
 the  subsystem of roots
$R_o$ corresponding to $\g'_0$, the decomposition
into irreducible submodules of the $\g_0$-module $\g_1$ which appear
in the decomposition (3.2), 
and the set of roots $R_1 = R^+ \setminus ({\mu}\cup R_o)$.\par
For a set of simple roots of $\g'_0$, we denote by 
$\{\pi_1, \dots, \pi_\ell\}$ the corresponding system of fundamental 
weights and, for any weight $\lambda = \sum a_i \pi_i$, we denote
by $V(\lambda)$ the irreducible $\g'_0$-module 
with highest weight $\lambda$.\par
\bigskip
In Table 2, we give the information
needed to determine  the holomorphic subspaces
$\m^{10}$  
when $\g^\C$ is a simple Lie algebra and the contact form
$\vartheta= -i\B \circ Z|_\h$ is parallel to 
a short root. \par
In Table 3
we give the same information for the cases
$\g^\C = B_3$ or $D_\ell$ and $\vartheta$   proportional to 
no root and associated with a primitive CR structure. \par
In both tables
we give the root
systems $R$, the  contact form $\vartheta$,
the subalgebra $\l^\C =  C_{\g^\C}(Z)\cap (Z)^\perp$, the root subsystem $R_o$ 
of  $\l^\C$ and  the list of 
the highest weights for the irreducible
$\k^\C$-moduli in $\m^\C$ ($\k^\C = C_{\g^\C}(Z)$). We group
the  highest weights
corresponding to equivalent $\l^\C$-moduli
 with curly
brackets.\par
\bigskip
In Table 4 we recall the Dynkin graphs associated with   indecomposable
root systems  and the correspondence  used in [7] between  
nodes   and  simple roots.
\par
\newpage
\centerline{\bf Table 1}
\medskip
\vbox{\offinterlineskip
\halign {\strut\vrule\hfil\ $#$\ \hfil 
 &\vrule\hfil\ $#$\ 
\hfil&\vrule\hfil\  $#$\ 
\hfil
&\vrule\hfil\  $#$\ 
\hfil
&\vrule\hfil\  $#$\ 
\hfil
&\vrule\hfil $#$\ 
\hfil
&\vrule\hfil\  $#$\ 
\hfil
\vrule\cr
\noalign{\hrule}
\phantom{\frac{\frac{1}{1}}{\frac{1}{1}}}\g
\ \ & 
R
& 
\mu
&
\g_0'
&
R_o
&
\g_1
&
R_1
\cr \noalign{\hrule}
A_\ell 
&
\smallmatrix
\varepsilon_i - \varepsilon_j\\
\phantom{a}\\
1 \leq i,j \leq \ell+1\\
\phantom{a}
\endsmallmatrix
&
\smallmatrix
\varepsilon_1 - \varepsilon_{\ell+1}
\endsmallmatrix
&
\smallmatrix
A_{\ell-2} + \R
\endsmallmatrix
&
\smallmatrix
\varepsilon_a - \varepsilon_b\\
\phantom{a}\\
2\leq a,b \leq\ell\\
\phantom{a}
\endsmallmatrix
&
\smallmatrix V(\pi_1) + \\
V(\pi_{\ell-2})
\endsmallmatrix
&
\smallmatrix
\varepsilon_1 - \varepsilon_a,
\ \varepsilon_a - \varepsilon_{\ell+1}\\
\phantom{a}\\
 2\leq a\leq \ell\\
\phantom{a}
\endsmallmatrix
 \cr \noalign{\hrule}
B_\ell 
&
\smallmatrix
\pm\varepsilon_i \pm \varepsilon_j, 
\ \pm\varepsilon_i\\
\phantom{a}\\
1 \leq i,j \leq \ell
\endsmallmatrix
&
\smallmatrix
\varepsilon_1 + \varepsilon_2
\endsmallmatrix
&
\smallmatrix
A_1 + B_{\ell-2}
\endsmallmatrix
&
\smallmatrix
\pm(\varepsilon_1 -\varepsilon_2),
\ \pm\varepsilon_a \pm \varepsilon_b\\
\pm\varepsilon_a\\
\phantom{a}\\
3\leq a,b\leq \ell\\
\phantom{a}
\endsmallmatrix
&
\smallmatrix 
V(\pi_1)\otimes
\\
V(\pi'_1)
\endsmallmatrix
&
\smallmatrix
\varepsilon_1,\ \varepsilon_2\\
\varepsilon_1 \pm \varepsilon_a,\  \varepsilon_2 \pm \varepsilon_a\\
\phantom{a}\\
3\leq a \leq \ell
\endsmallmatrix
 \cr \noalign{\hrule}
C_\ell 
&
\smallmatrix
\pm\varepsilon_i \pm \varepsilon_j\ ,
\ \pm2\varepsilon_i\\
\phantom{a}\\
1\leq i,j \leq \ell\\
\phantom{a}
\endsmallmatrix
&
\smallmatrix
2\varepsilon_1
\endsmallmatrix 
&
\smallmatrix
C_{\ell-1}
\endsmallmatrix
&
\smallmatrix
\pm\varepsilon_a \pm \varepsilon_b,
\ \pm2\varepsilon_a\\
\phantom{a}\\
2\leq a,b \leq  \ell\\
\phantom{a}
\endsmallmatrix
&
\smallmatrix 
V(\pi_1)
\endsmallmatrix
&
\smallmatrix
\varepsilon_1 \pm \varepsilon_a\\
\phantom{a}\\
 2\leq a\leq \ell\\
\phantom{a}
\endsmallmatrix
 \cr \noalign{\hrule}
D_\ell 
&
\smallmatrix
\pm\varepsilon_i \pm \varepsilon_j\\
\phantom{a}\\
1\leq i,j\leq \ell\\
\phantom{a}
\endsmallmatrix
&
\smallmatrix
\varepsilon_1 + \varepsilon_2
\endsmallmatrix
&
\smallmatrix A_1 + 
D_{\ell-2}  
\endsmallmatrix
&
\smallmatrix
\pm(\varepsilon_1 -\varepsilon_2),
\ \pm\varepsilon_a \pm \varepsilon_b\\
\phantom{a}\\
3 \leq a,b \leq  \ell\\
\phantom{a}
\endsmallmatrix
&
\smallmatrix 
V(\pi_1)\otimes
\\
V(\pi'_1)
\endsmallmatrix
&
\smallmatrix
\varepsilon_1 \pm \varepsilon_a,
\  \varepsilon_2 \pm \varepsilon_a\\
\phantom{a}\\
3\leq a\leq \ell\\
\phantom{a}
\endsmallmatrix
 \cr \noalign{\hrule}
E_6 
&
\smallmatrix
\varepsilon_i - \varepsilon_j,
\  \pm2\varepsilon\\
\varepsilon_i + \varepsilon_j +
\varepsilon_k \pm\varepsilon\\
\phantom{a}\\
1\leq i,j,k \leq 6\\
\phantom{a}
\endsmallmatrix
&
\smallmatrix
2\varepsilon
\endsmallmatrix
&
\smallmatrix
A_{5}
\endsmallmatrix
&
\smallmatrix
\varepsilon_i - \varepsilon_j
\endsmallmatrix
&
\smallmatrix 
V(\pi_1)
\endsmallmatrix
&
\smallmatrix
\varepsilon_i + \varepsilon_j +
\varepsilon_k \pm\varepsilon\\
\phantom{a}
\endsmallmatrix
 \cr \noalign{\hrule}
E_7 
&
\smallmatrix
\varepsilon_i - \varepsilon_j\\
\varepsilon_i + \varepsilon_j +
\varepsilon_k +\varepsilon_\ell\\
\phantom{a}\\
1\leq i,j,k,\ell\leq 8
\\
\phantom{a}
\endsmallmatrix
&
\smallmatrix
-\varepsilon_7 + \varepsilon_8
\endsmallmatrix
&
\smallmatrix
D_{6}
\endsmallmatrix
&
\smallmatrix
\varepsilon_a - \varepsilon_b\\
\varepsilon_7 + \varepsilon_8 +
\varepsilon_a +\varepsilon_b\\
\varepsilon_a + \varepsilon_b +
\varepsilon_c +\varepsilon_d\\
\phantom{a}\\
1\leq a,b, c, d\leq 6\\
\phantom{a}
\endsmallmatrix
&
\smallmatrix 
V(\pi_1)
\endsmallmatrix
&
\smallmatrix
-\varepsilon_7 + \varepsilon_a,
\ \varepsilon_8 - \varepsilon_a\\
\varepsilon_8 + \varepsilon_a +
\varepsilon_b +\varepsilon_c\\
\phantom{a}\\
1\leq a,b,c \leq 6\\
\phantom{a}
\endsmallmatrix
 \cr \noalign{\hrule}
E_8 
&
\smallmatrix
\varepsilon_i - \varepsilon_j\\
\pm(\varepsilon_i + \varepsilon_j +
\varepsilon_k)\\
\phantom{a}\\
1\leq i,j,k \leq 9\\
\phantom{a}
\endsmallmatrix
&
\smallmatrix
\varepsilon_1 - \varepsilon_9
\endsmallmatrix
&
\smallmatrix
E_{7}
\endsmallmatrix
&
\smallmatrix
\varepsilon_a - \varepsilon_b\\
\pm(\varepsilon_1 + \varepsilon_9 +
\varepsilon_a)\\
\pm(\varepsilon_a + \varepsilon_b +
\varepsilon_c)\\
\phantom{a}\\
2\leq a, b,c\leq 8\\
\phantom{a}
\endsmallmatrix
&
\smallmatrix 
V(\pi_1)
\endsmallmatrix
&
\smallmatrix
\varepsilon_1 - \varepsilon_a,
\ -\varepsilon_9 + \varepsilon_a\\
\varepsilon_1 + \varepsilon_a +
\varepsilon_b\\
\phantom{a}\\
2\leq a,b\leq 8\\
\phantom{a}
\endsmallmatrix
 \cr \noalign{\hrule}
F_4 
&
\smallmatrix
\frac{\pm\varepsilon_1 \pm \varepsilon_2
\pm\varepsilon_3 \pm \varepsilon_4}{2}\\
\pm\varepsilon_i \pm \varepsilon_j,
\ \pm\varepsilon_i\\
\phantom{a}\\
1\leq i, j \leq 4\\
\phantom{a}
\endsmallmatrix
&
\smallmatrix
\varepsilon_1 + \varepsilon_2
\endsmallmatrix
&
\smallmatrix
C_3
\endsmallmatrix
&
\smallmatrix
\pm(\varepsilon_1 - \varepsilon_2)\\
\pm\frac{\varepsilon_1 - \varepsilon_2
\pm\varepsilon_3 \pm \varepsilon_4}{2}\\
\pm\varepsilon_a,
\ \pm\varepsilon_a \pm\varepsilon_b\\
\phantom{a}\\
3\leq a,  b\leq 4\\
\phantom{a}
\endsmallmatrix
&
\smallmatrix V(\pi_1)
\endsmallmatrix
&
\smallmatrix
\varepsilon_1,
\  \varepsilon_2\\
\frac{\varepsilon_1 + \varepsilon_2
\pm\varepsilon_3 \pm \varepsilon_4}{2}
\endsmallmatrix
 \cr \noalign{\hrule}
G_2 
&
\smallmatrix
\varepsilon_i - \varepsilon_j\ ,\ \pm\varepsilon_i\\
\phantom{a}\\
1\leq i,j\leq 3\\
\phantom{a}
\endsmallmatrix
&
\smallmatrix
\varepsilon_1 - \varepsilon_2
\endsmallmatrix
&
\smallmatrix
A_1
\endsmallmatrix
&
\smallmatrix
\pm\varepsilon_3
\endsmallmatrix
&
\smallmatrix 
V(\pi_1)
\endsmallmatrix
&
\smallmatrix
\varepsilon_1 - \varepsilon_3\ ,\  \varepsilon_1- \varepsilon_2\\
\varepsilon_3 - \varepsilon_2
\endsmallmatrix
 \cr \noalign{\hrule}
}}
\bigskip
\bigskip
\centerline{\bf Table 2}
\medskip
\vbox{\offinterlineskip
\halign {\strut\vrule\hfil\ $#$\ \hfil 
 &\vrule\hfil\ $#$\ 
\hfil&\vrule\hfil\  $#$\ 
\hfil
&\vrule\hfil\  $#$\ 
\hfil
&\vrule\hfil $#$
\hfil
&\vrule\hfil $#$\
\hfil
\vrule\cr
\noalign{\hrule}
\phantom{\frac{\frac{1}{1}}{\frac{1}{1}}}\g
\ \ & 
R
& 
\underset{\underset{=i\B\circ Z}\to
{\phantom{A}}}\to\vartheta
&
\underset{\underset{=C_{\g^\C}(Z)\cap (Z)^\perp}\to
{\phantom{A}}}\to
{\l^\C} 
&
R_o
&
\smallmatrix
\text{highest weights
for}\ \m^\C\\
\text{grouped into  sets of} \\
\text{equivalent}\ \l^\C-\text{moduli}
\endsmallmatrix
\cr \noalign{\hrule}
B_\ell 
&
\smallmatrix
\pm\varepsilon_i \pm \varepsilon_j, 
\ \pm\varepsilon_i\\
\phantom{a}\\
1\leq i,j\leq  \ell\\
\phantom{a}
\endsmallmatrix
&
\smallmatrix
\varepsilon_1
\endsmallmatrix
&
\smallmatrix
B_{\ell-1}
\endsmallmatrix
&
\smallmatrix
\pm\varepsilon_a \pm \varepsilon_b,
\ \pm\varepsilon_a\\
\phantom{a}\\
2\leq a,b  \leq \ell\\
\phantom{a}
\endsmallmatrix
&
\smallmatrix
\{\varepsilon_1+\varepsilon_2, -\varepsilon_1+\varepsilon_2\}
\endsmallmatrix
 \cr \noalign{\hrule}
C_\ell 
&
\smallmatrix
\pm\varepsilon_i \pm \varepsilon_j,
\ \pm 2 \varepsilon_i\\
\phantom{a}\\
1\leq 
i,j \leq \ell\\
\phantom{a}
\endsmallmatrix
&
\smallmatrix
\varepsilon_1 + \varepsilon_2
\endsmallmatrix 
&
\smallmatrix
A_1 + C_{\ell-2}
\endsmallmatrix
&
\smallmatrix
\pm(\varepsilon_1 -\varepsilon_2),
\ \pm2\varepsilon_a\\
\pm\varepsilon_a \pm \varepsilon_b\\
\phantom{a}\\
3\leq a,b\leq  \ell\\
\phantom{a}
\endsmallmatrix
&
\smallmatrix
\{2\varepsilon_1, - 2\varepsilon_2\}\\
\{\varepsilon_1 + \varepsilon_3,
\ -\varepsilon_2 +\varepsilon_3\} 
\endsmallmatrix
 \cr \noalign{\hrule}
F_4 
&
\smallmatrix
\phantom{a}\\
\frac{\pm\varepsilon_1 \pm \varepsilon_2
\pm\varepsilon_3 \pm \varepsilon_4}{2}\\
\pm\varepsilon_i \pm \varepsilon_j,
\ \pm\varepsilon_i\\
\phantom{a}\\
1\leq i,j \leq 4\\
\phantom{a}
\endsmallmatrix
&
\smallmatrix
\varepsilon_1
\endsmallmatrix
&
\smallmatrix
B_3
\endsmallmatrix
&
\smallmatrix
\pm\varepsilon_a \pm \varepsilon_b,
\ \pm \varepsilon_a\\
\phantom{a}\\
2\leq a,b\leq 4\\
\phantom{a}
\endsmallmatrix
&
\smallmatrix
\{\varepsilon_1+\varepsilon_2, -\varepsilon_1+\varepsilon_2\}\\
\{\frac{\varepsilon_1 + \varepsilon_2+
\varepsilon_3 + \varepsilon_4}{2}, 
\frac{-\varepsilon_1 +\varepsilon_2 +
\varepsilon_3 + \varepsilon_4}{2}\}
\endsmallmatrix
 \cr \noalign{\hrule}
G_2 
&
\smallmatrix
\phantom{a}\\
\varepsilon_i - \varepsilon_j\ , \ 
\pm\varepsilon_i\\
\phantom{a}\\
1\leq i,j\leq  3\\
\phantom{a}
\endsmallmatrix
&
\smallmatrix
\varepsilon_1
\endsmallmatrix
&
\smallmatrix
A_1
\endsmallmatrix
&
\smallmatrix
\pm(\varepsilon_2 - \varepsilon_3)
\endsmallmatrix
&
\smallmatrix
\{\varepsilon_1 , - \varepsilon_1\}\\
\{\varepsilon_3 -  \varepsilon_1,\ 
\varepsilon_3 ,\ 
-\varepsilon_2 ,\  
\varepsilon_1 - \varepsilon_2\}
\endsmallmatrix
\cr \noalign{\hrule}
}}
\bigskip
\bigskip
\newpage
\centerline{\bf Table 3}
\medskip
\moveright 0.4 cm
\vbox{\offinterlineskip
\halign {\strut\vrule\hfil\ $#$\ \hfil 
 &\vrule\hfil\ $#$\ 
\hfil&\vrule\hfil\  $#$\ 
\hfil
&\vrule\hfil\  $#$\ 
\hfil
&\vrule\hfil $#$
\hfil
&\vrule\hfil $#$\
\hfil
\vrule\cr
\noalign{\hrule}
\phantom{\frac{\frac{1}{1}}{\frac{1}{1}}}\g
\ \ & 
R
& 
\underset{\underset{=i\B\circ Z}\to
{\phantom{A}}}\to\vartheta
&
\underset{\underset{=C_{\g^\C}(Z)\cap (Z)^\perp}\to
{\phantom{A}}}\to
{\l^\C} 
&
R_o
&
\smallmatrix
\text{highest weights
for}\ \m^\C\\
\text{grouped into  sets of} \\
\text{equivalent}\ \l^\C-\text{moduli}
\endsmallmatrix
\cr \noalign{\hrule}
B_3 
&
\smallmatrix
\pm\varepsilon_i \pm \varepsilon_j, 
\ \pm\varepsilon_i\\
\phantom{a}\\
1\leq i,j\leq  3\\
\phantom{a}
\endsmallmatrix
&
\smallmatrix
\varepsilon_1 + \varepsilon_2 + 
\varepsilon_3
\endsmallmatrix
&
\smallmatrix
A_2
\endsmallmatrix
&
\smallmatrix
\pm(\varepsilon_a - \varepsilon_b),
\\
\phantom{a}\\
1\leq a,b  \leq 3\\
\phantom{a}
\endsmallmatrix
&
\smallmatrix
\{\varepsilon_1+\varepsilon_2, -\varepsilon_3\}\quad
\{-\varepsilon_2 -\varepsilon_3, \ \varepsilon_1\}
\endsmallmatrix
 \cr \noalign{\hrule}
D_\ell 
&
\smallmatrix
\pm\varepsilon_i \pm \varepsilon_j\\
\phantom{a}\\
1\leq i, j\leq  \ell\\
\phantom{a}
\endsmallmatrix
&
\smallmatrix
\varepsilon_1
\endsmallmatrix
&
\smallmatrix
D_{\ell -1}
\endsmallmatrix
&
\smallmatrix
\pm\varepsilon_i \pm \varepsilon_j\\
\phantom{a}\\
2\leq i, j \leq \ell\\
\phantom{a}
\endsmallmatrix
&
\smallmatrix
\{ \varepsilon_1 + \varepsilon_2, -\varepsilon_1 + \varepsilon_2\}
\endsmallmatrix
 \cr \noalign{\hrule}
}}
\bigskip
\bigskip
\centerline{\bf Table 4}
\medskip
\moveright 0.5 cm
\vbox{\offinterlineskip
\halign {\strut\vrule\hfil\ $#$\ \hfil 
 &\vrule\hfil\ $#$\ 
\hfil&\vrule\hfil\  $#$\ 
\hfil\vrule\cr
\noalign{\hrule}
\phantom{\frac{\frac{1}{1}}{\frac{1}{1}}}
\text{Type of}\ G
& 
\text{Dynkin graphs}
&
\text{Simple roots}
\cr \noalign{\hrule}
A_\ell
& 
\underset{\phantom{a}}\to{\overset{\phantom{a}}\to{
\dynkin\wroot{\ 1}\link\wroot{\ 2}\link
\dots\link\wroot{\ \ell-1}\link\wroot{\ \ell}
\enddynkin
}}
&
\smallmatrix
\alpha_i = \varepsilon_i - \varepsilon_{i+1}
\endsmallmatrix
\cr \noalign{\hrule}
B_\ell
& 
\underset{\phantom{a}}\to{\overset{\phantom{a}}\to{
\dynkin\wroot{\ 1}\link\wroot{\ 2}\link
\dots\link\wroot{\ \ell-1}\llink>\wroot{\ \ell}
\enddynkin
}}
&
\smallmatrix\alpha_i = \varepsilon_i - \varepsilon_{i+1}\ (i<\ell)\ ,\\ 
\phantom{a}\\
\alpha_\ell = \varepsilon_\ell
\endsmallmatrix
\cr \noalign{\hrule}
C_\ell
& 
\underset{\phantom{a}}\to{\overset{\phantom{a}}\to{
\dynkin\wroot{\ 1}\link\wroot{\ 2}\link
\dots\link\wroot{\ \ell-1}\llink<\wroot{\ \ell}
\enddynkin
}}
&
\smallmatrix
\alpha_i = \varepsilon_i - \varepsilon_{i+1}\ (i<\ell)\ ,\\ 
\phantom{a}\\
\alpha_\ell = 2\varepsilon_\ell
\endsmallmatrix
\cr \noalign{\hrule}
D_\ell
& 
\underset{\phantom{a}}\to{\overset{\phantom{a}}\to{
\dynkin\wroot{\ 1}\link\wroot{\ 2}\link
\dots\link\wroot{\ \ell-2}\wrootupright{\ \ell-1}
\wrootdownright{\ \ell}\enddynkin
}}
&
\smallmatrix
\alpha_i = \varepsilon_i - \varepsilon_{i+1}\ (i<\ell)\ ,\\
\phantom{a}\\
\alpha_\ell = \varepsilon_{\ell-1} + \varepsilon_\ell
\endsmallmatrix
\cr \noalign{\hrule}
E_6
&
\underset{\phantom{a}}\to{\overset{\phantom{a}}\to{
\dynkin\wroot{\ 1}\link\wroot{\ 2}\link\wroot{\ 3}
\wrootdown{6}\link\wroot{\ 4}
\link\wroot{\ 5}
\enddynkin
}}
& \smallmatrix
\alpha_i = \varepsilon_i - \varepsilon_{i+1}\ (i<6)\ ,\\
\phantom{a}\\
\alpha_6 = \varepsilon_4 + \varepsilon_5 + \varepsilon_6 + \varepsilon
\endsmallmatrix
\cr \noalign{\hrule}
E_7
&
\underset{\phantom{a}}\to{\overset{\phantom{a}}\to{
\dynkin\wroot{\ 1}\link\wroot{\ 2}\link\wroot{\ 3}\link\wroot{\ 4}
\wrootdown{7}\link\wroot{\ 5}
\link\wroot{\ 6}
\enddynkin
}}
& \smallmatrix
\alpha_i = \varepsilon_i - \varepsilon_{i+1}\ (i<7)\ ,\\
\phantom{a}\\
\alpha_7 = \varepsilon_5 + \varepsilon_6 + \varepsilon_7 + \varepsilon_8
\endsmallmatrix
\cr \noalign{\hrule}
E_8
&
\underset{\phantom{a}}\to{\overset{\phantom{a}}\to{
\dynkin\wroot{\ 1}\link\wroot{\ 2}\link\wroot{\ 3}\link
\wroot{\ 4}\link\wroot{\
5}
\wrootdown{8}\link\wroot{\ 6}
\link\wroot{\ 7}
\enddynkin
}}
& \smallmatrix
\alpha_i = \varepsilon_i - \varepsilon_{i+1}\ (i<8)\ ,\\
\phantom{a}\\
\alpha_8 = \varepsilon_6 +  \varepsilon_7 + \varepsilon_8
\endsmallmatrix
\cr \noalign{\hrule}
F_4
&
\underset{\phantom{a}}\to{\overset{\phantom{a}}\to{
\dynkin\wroot{\ 1}\link\wroot{\ 2}\llink<\wroot{\ 3}
\link\wroot{\ 4}\enddynkin
}}
&
\smallmatrix
\alpha_1 = (\varepsilon_1 - \varepsilon_2 - \varepsilon_3 - 
\varepsilon_4)/2\ , \
\alpha_2 = \varepsilon_4\\
\phantom{a}\\
\alpha_3 = \varepsilon_3 - \varepsilon_4\ ,\ 
\alpha_4 = \varepsilon_2 - \varepsilon_3
\endsmallmatrix
\cr \noalign{\hrule}
G_2
&
\underset{\phantom{a}}\to{\overset{\phantom{a}}\to{
\dynkin\wroot{\ 1}\lllink<\wroot{\ 2}\enddynkin
}}
&
\smallmatrix
\alpha_1 = - \varepsilon_2\ ,
\ \alpha_2 = \varepsilon_2 - \varepsilon_3
\endsmallmatrix
\cr \noalign{\hrule}
}}


\bigskip
\bigskip
\Refs
\widestnumber\key{12}

\ref
\key 1
\by H. Azad, A. Huckleberry and W. Richthofer
\paper Homogeneous CR manifolds 
\jour J. Reine und Angew. Math.
\vol 358
\yr 1985
\pages 125--154
\endref

\ref
\key 2
\by D. V. Alekseevsky
\paper Contact homogeneous spaces
\jour Funktsional. Anal. i Prilozhen.
\vol 24
\yr 1990 
\issue 4
\transl\nofrills Engl. transl. in
\jour Funct. Anal. Appl.  
\vol 24
\yr 1991
\issue 4
\pages 324--325
\endref

\ref
\key 3
\by D. V. Alekseevsky
\paper Flag Manifolds
\inbook Sbornik Radova, 11 Jugoslav. Seminr.
\vol 6
\issue 14
\yr 1997
\publ Beograd
\pages 3--35
\endref

\ref
\key 4
\by D. V. Alekseevsky and A. M. Perelomov
\paper Invariant Kahler-Einstein metrics
on compact homogeneous spaces
\jour Funktsional. Anal. i Prilozhen.
\vol 20
\yr 1986 
\issue 3
\transl\nofrills Engl. transl. in
\jour Funct. Anal. Appl.  
\vol 20
\yr 1986
\issue 3
\pages 171--182
\endref

\ref
\key 5
\by M. Bordermann, M. Forger and H. R\"omer
\paper Homogeneous K\"ahler Manifolds: paving the way
towards new supersymmetric Sigma Models
\jour Comm. Math. Phys.
\vol 102
\yr 1986
\pages 605--647
\endref

\ref
\key 6
\by C. P. Boyer, K. Galicki and B. M. Mann
\paper The geometry and topology of 3-Sasakian manifolds
\jour J. Reine und Angew. Math.
\yr 1994
\vol 144
\pages 183--220
\endref

\ref
\key 7
\by V. V. Gorbatsevic, A. L. Onishchik and E. B. Vinberg
\paper Structure of Lie Groups and Lie Algebras
\inbook in Encyclopoedia of Mathematical Sciences -
Lie Groups and Lie Algebras III
\ed A. L. Onishchik and E. B. Vinberg
\publ Springer-Verlag -- VINITI
\publaddr Berlin
\yr 1993 (Russian edition: VINITI, Moscow,1990) 
\endref

\ref
\key 8
\by A. Huckleberry and W. Richthofer
\paper Recent Developments in Homogeneous CR Hypersurfaces
\inbook in Contributions to Several Complex Variables
\bookinfo Aspects of Math. 
\vol E9 
\publ Vieweg, Brauschweig 
\yr 1986
\pages 149--177
\endref

\ref
\key 9
\by A. Morimoto and T. Nagano
\paper On pseudo-conformal transformations of
hypersurfaces
\jour J. Math. Soc. Japan
\vol 15
\yr 1963
\pages 289--300
\endref

\ref
\key 10
\by M. Nishiyama
\paper Classification of invariant complex
structures on irreducible compact simply connected coset spaces
\jour Osaka J. Math.
\vol 21
\yr 1984
\pages 39--58
\endref

\ref
\key 11
\by W. Richthofer
\paper Homogene CR - Mannigfaltigkeiten
\paperinfo Dissertation zur Erlangung des Doktorgrades
der Abteilung f\"ur Mathematik
der Ruhr - Universit\"at 
\publaddr Bochum
\yr 1985
\endref

\ref
\key 12
\by A. Spiro
\paper Groups acting transitively on
compact CR manifolds of hypersurface type
\jour Proc. A.M.S.
\vol 128
\issue 4
\yr 1999
\pages 1141--1145
\endref

\endRefs
\enddocument